\documentclass{amsart}

\usepackage[english]{babel}
\usepackage{array,booktabs,float}
\usepackage{amsmath,amssymb,amsfonts,amscd}
\usepackage{amsthm}
\usepackage{bm}
\usepackage{amscd}
\usepackage{mathtools}
\usepackage{mathrsfs}
\usepackage{pb-diagram}
\usepackage{booktabs}
\usepackage{cleveref}
 
\usepackage{todonotes}
\usepackage{comment}
\usepackage{tikz-cd}

\theoremstyle{plain}
\newtheorem{theorem}{Theorem}[section]
\crefname{theorem}{Theorem}{Theorems}
\Crefname{theorem}{Theorem}{Theorems}
\newtheorem{proposition}[theorem]{Proposition}
\crefname{proposition}{Proposition}{Propositions}
\Crefname{proposition}{Proposition}{Propositions}
\newtheorem{lemma}[theorem]{Lemma}
\crefname{lemma}{Lemma}{Lemmas}
\Crefname{lemma}{Lemma}{Lemmas}
\newtheorem{corollary}[theorem]{Corollary}
\crefname{corollary}{Corollary}{Corollaries}
\Crefname{corollary}{Corollary}{Corollaries}

\crefname{claim}{Claim}{Claims}
\Crefname{claim}{Claim}{Claims}

\crefname{property}{Property}{Properties}
\Crefname{property}{Property}{Properties}

\crefname{problem}{Problem}{Problems}
\Crefname{problem}{Problem}{Problems}

\DeclareMathOperator{\coker}{coker}

\theoremstyle{definition}
\newtheorem{definition}[theorem]{Definition}
\crefname{definition}{Definition}{Definitions}
\Crefname{definition}{Definition}{Definitions}

\crefname{notation}{Notation}{Notations}
\Crefname{notation}{Notation}{Notations}

\crefname{convention}{Convention}{Conventions}
\Crefname{convention}{Convention}{Conventions}

\crefname{condition}{Condition}{Conditions}
\Crefname{condition}{Condition}{Conditions}

\crefname{assumption}{Assumption}{Assumptions}
\Crefname{assumption}{Assumption}{Assumptions}

\theoremstyle{remark}
\newtheorem{remark}[theorem]{Remark}
\crefname{remark}{Remark}{Remarks}
\Crefname{remark}{Remark}{Remarks}

\crefname{example}{Example}{Examples}
\Crefname{example}{Example}{Examples}

\crefname{section}{Section}{Sections}
\Crefname{section}{Section}{Sections}
\crefname{subsection}{Subsection}{Subsections}
\Crefname{subsection}{Subsection}{Subsections}
\crefname{figure}{Figure}{Figures}
\Crefname{figure}{Figure}{Figures}

\newtheorem*{acknowledgement}{Acknowledgement}

\newcommand{\N}{\mathbb{N}}
\newcommand{\Z}{\mathbb{Z}}
\newcommand{\Q}{\mathbb{Q}}
\newcommand{\R}{\mathbb{R}}
\newcommand{\C}{\mathbb{C}}

\newcommand{\ind}{\mathrm{ind}}
\newcommand{\Slash}[1]{{\ooalign{\hfil/\hfil\crcr$#1$}}}
\newcommand{\Spinc}{\mathrm{Spin}^c}

\begin{document}

\title[Exotic $P^2$-knots]{A gauge theoretic invariant of embedded surfaces in $4$-manifolds and exotic $P^2$-knots}
\author{Jin Miyazawa}
\address{Graduate School of Mathematical Sciences, the University of Tokyo, 3-8-1 Komaba, Meguro, Tokyo 153-8914, Japan}
\email{miyazawa@ms.u-tokyo.ac.jp}
 
\begin{abstract}
    We give an infinite family of embeddings of $\R P^2$ to $S^4$ such that they are mutually topologically isotopic however are not smoothly isotopic to each other. Moreover, they are topologically isotopic to the standard $P^2$-knot.  
    To prove that these $P^2$-knots are not smoothly isotopic to each other, we construct a gauge theoretic invariant of embedded surfaces in $4$-manifolds using a variant of the Seiberg--Witten theory, which is called the Real Seiberg--Witten theory. 
\end{abstract}
\maketitle
\tableofcontents
\section{Introduction}
We explain our results and backgrounds. In this paper, we assume embeddings of surfaces in $4$-manifolds are smooth otherwise stated. When we consider topological embeddings of surfaces, we always assume that the embeddings are locally flat. 
\subsection{Exotic surfaces}

Let $\Sigma$ be a compact surface. Two embedding of surfaces \[j \colon \Sigma \hookrightarrow X \;\text{and}\; j' \colon \Sigma \hookrightarrow X\] in a smooth $4$-manifold $X$ are said to be smoothly (resp. topologically) isotopic if there is a diffeomorphism (resp. homeomorphism) $\beta \colon X \to X$ such that $\beta$ is smoothly (resp. topologically) isotopic the identity and $j'=\beta \circ j$. We call $j$ and $j'$ are exotic pair if they are topologically isotopic and they are not smoothly isotopic. We sometimes identify the map $j$ with its image $j(\Sigma)$.   
   The problem of the existence of exotic surfaces attracts a great deal of attention. 
However, there are some open problems of the existence of the exotic surfaces. 

\subsubsection{Non-orientable surfaces in $S^4$}
There are only few examples of exotic pair of non-orientable surfaces in $S^4$ while there are a great deal of the research of the non-orientable surface knots in $S^4$. 
The first example of the exotic pair of non-orientable surfaces are given by Finashin--Kreck--Viro \cite{MR0970078}. 
They construct infinite family of embeddings of $\#_{10} \R P^2$ to $S^4$ 
such that they are topologically isotopic however they are not smoothly isotopic each other. 
Finashin \cite{MR2528682} give an example of an infinite family of exotic embeddings of $\#_{6} \R P^2$ to $S^4$. 
Among all known examples, this has the minimum non-orientable genus. 
Other examples are given by Havens \cite{havens2021equivariant} and Levine--Lidman--Piccirillo \cite[Theorem 1.17.]{levine2023new}. 

On the other hand, in the topological counter part, there is a strong unknotting result of non-orientable surface knots in $S^4$ given by Conway--Orson--Powell \cite{conway2023unknotting}. 

A reason why there are no known examples of exotic embeddings of non-orientable surfaces with small genus is that it is difficult to prove that such embeddings are not smoothly isotopic. 
In the examples above, they prove that the double branched cover along the embedded non-orientable surfaces are not diffeomorphic to each other. If we want to use this strategy to give examples of exotic embeddings of non-orientable surfaces with small genus, we have to give examples of small exotic $4$-manifolds which is known to be difficult to detect.

In this paper, we construct a new smooth isotopy invariant of embeddings of $\R P^2$ to $S^4$, which are called \textit{$P^2$-knot}, and we give an infinite family of $P^2$-knot 
 such that they are mutually topologically isotopic 
 however they are not smoothly isotopic to each other. 
Moreover, they are topologically isotopic to the \textit{standard $P^2$-knot}, or called $P^2$-unknot. 
There are two standard $P^2$-knots distinguished by normal Euler number, a simple topological invariant (see \cref{def_of_normal_euler_number}). 
The standard $P^2$-knots are given by the image of the map $\R P^2 \to \pm \C P^2 \to \pm \C P^2/(\text{complex conjugate}) \cong S^4$, where the first map is the embedding to the fixed point part of the complex conjugate. We denote $P_+ \subset S^4$ be the standard $P^2$-knot with normal Euler number is $2$. Note that the double branched cover of $S^4$ along $P_+$ is $-\C P^2$. 

Suppose a $P^2$-knot $P$ is smoothly (resp. topologically) isotopic to the standard $P^2$-knot, we call $P$ is \textit{smoothly (resp. topologically) unknotted}. Our main result \cref{main0} can be expressed as follows: The $P^2$-knot version of the smooth unknotting conjecture is false. 

\begin{theorem}[see \cref{main_theorem_2}]\label{main0}
    There is an infinite family of $P^2$-knots $\{ P_n\}_{n=0}^{\infty}$ such that $P_n$ is topologically isotopic to $P_+$ for all $n$ however the family of the pairs $\{(S^4, P_n)\}$ are not mutually pairwise diffeomorphic.  
\end{theorem}

The definition of $P_n$ is as follows: Let $K$ be the $(-2, 3, 7)$-pretzel knot and let $\tau_{0, 1}(K)$ be the $1$-roll spun $2$-knot of $K$ in $S^4$ (see \cref{Twist and rolling spun knots}). Then $P_n$ is given by the inner connected sum $P_+ \# n\tau_{0, 1}(K)$, where $n\tau_{0, 1}(K)$ is a $2$-knot given by $n$ times inner connected sum of the copies of $\tau_{0, 1}(K)$.

To prove these $P^2$-knots are not smoothly isotopic, we construct an invariant of embedded surfaces in $4$-manifold using the Real Seiberg--witten theory. We explain the invariant in \cref{intro_Real_SW}. To confirm that they are topologically isotopic, we use the result of Conway--Orson--Powell \cite[Theorem A.]{conway2023unknotting} that state two embeddings of non-orientable surface whose knot group is $\Z/2$ are topologically unknotted under some condition on the normal Euler number and genus. 

By using our invariant, we deduce that the complement of the exotic $P^2$-knots are exotic rational homology balls. We prove the following theorem. 
\begin{theorem}[see \cref{complement}]\label{complement_}
   Let $P_n$ be a $P^2$-knot in the statement in \cref{main0}. Let $\nu_n$ be a tubular neighborhood of $P_n$. Then $X_n:=S^4\setminus \nu_n$ are homeomorphic to $S^4\setminus \nu(P_+)$ however they are not diffeomorphic to each other. 
\end{theorem}

We also solving a negative answer to the following question by Kamada \cite[Epilogue]{MR3588325}: Let $K$ be a knot in $S^3$ and let  $\tau_{k, 0}(K)$ be the $k$-twist spun $2$-knot of $K$. Then the question is that whether or not $P_+ \# \tau_{k, 0}(K)$ and $P_+ \# \tau_{k', 0}(K)$ are smoothly isotopic when $k \equiv k' \mod 2$. 
One can easily check that the fundamental group of $S^4\setminus P_+ \# \tau_{k, 0}(K)$ and $S^4\setminus P_+ \# \tau_{k', 0}(K)$ are the same if $k \equiv k' \mod 2$. 
If $k$ is odd, then $\pi_1(S^4\setminus P_+ \# \tau_{k', 0}(K))\cong \Z/2$. In this case the problem above asked by Melvin in the Kirby list \cite[Problem 4.58]{kirby1997problems}. 

\begin{theorem}[see \cref{even_twist}]\label{main00}
    Let $K$ be the $(-2, 3, 7)$-pretzel knot and $k$ and $k'$ be two even integers and let us denote the $k$-twist spun $2$-knot of $K$ by $\tau_{k, 0}(K)$. If $k - k'$ is not divided by $4$, then two $P^2$-knots $P_+ \# \tau_{k, 0}(K)$ and $P_+ \# \tau_{k', 0}(K)$ are not smoothly isotopic. Moreover, $(S^4, P_+ \# \tau_{k, 0}(K))$ and $(S^4, P_+ \# \tau_{k', 0}(K))$ are not pairwise diffeomorphic. 
\end{theorem}

We do not know whether $P_+ \# \tau_{k, 0}(K)$ and $P_+ \# \tau_{k', 0}(K)$ are topologically isotopic or not. 

\subsubsection{Exotic surfaces in other closed $4$-manifolds}
We can apply our invariant to prove that two embedded orientable surfaces are not smoothly isotopic. To the best of our knowledge, the following example is new when $N \ge 5$. 

\begin{theorem}\label{main1}
    Let $X= \# M \overline{\C P^2}\# {N} \C P^2$. If $M \ge 20$ and $N \ge 4$, there are exotic surfaces $S$ and $S'$ in $X$ given by embeddings of $S^2$ which represent the homology class 
    \[
        \sum_{i=1}^{N-3} 2[\C P^1]_i 
    \] 
where $[\C P^1]_i$ is the homology class represented by the $\C P^1$ in the $i$-th connected component of $\# {N} \C P^2$. Moreover, $(X, S)$ and $(X, S')$ are not pairwise diffeomorphic.
\end{theorem}

There are many examples of orientable exotic surfaces in closed $4$-manifolds. For example, Fintushel--Stern \cite{MR1492129}, Finashin \cite{finashin2002knotting},Auckly--Kim--Melvin--Ruberman \cite{auckly2019isotopy}, Kim--Ruberman \cite{kim2008topological,kim2008smooth,kim2011double}, Kim \cite{kim2006modifying}, Mark \cite{mark2013knotted}, and Hoffman--Sunukjian \cite{hoffman2013surfaces} proved that the complement of the surfaces are not diffeomorphic using Seiberg--Witten invariant or Ozsv\'{a}th--Szabo invariant. Recently, Baraglia \cite{baraglia2020adjunction} gives examples of the exotic surfaces in closed $4$-manifolds whose complements are diffeomorphic using Family Seiberg--witten invariant. 
 
In fact, when $N=4$ in \cref{main1}, we can prove that the surface in \cref{main1} is exotic by proving the complement of the tubular neighborhood of the surfaces are not diffeormorphic using the relative Seiberg--Witten invariants. 
However, if $N\ge 5$, the author does not know how to prove \cref{main1} by using the Seiberg--Witten invariant or the Ozsv\'{a}th--Szabo mixed invariant. 

We have a result on an infinite family of embedding of $S^2$ into $\C P^2$ whose homology class is $2[\C P^2]$. Let $C \subset \C P^2$ be a holomorphic embedding of $S^2$ whose homology class is $2[\C P^1] \in H_2(\C P^2, \Z)$. Let $K$ be the $(-2, 3, 7)$-pretzel knot and $\tau_{0, 1}(K)$ is $1$-roll spun knot in $S^4$. Let $C_n:=C \# n \tau_{0, 1}(K)$ be an embedding of $S^2$ in $\C P^2 \cong \C P^2 \# S^4$ given by inner connected sum of $C$ and $n \tau_{0, 1}(K)$. 
\begin{theorem}[see \cref{conic}]
    We have $\pi_1(\C P^2 \setminus C_n)\cong \pi_1(\C P^2\setminus C) \cong \Z/2$ however these embeddings $\{(\C P^2, C_n)\}_{n=0}^{\infty}$ are not mutually pairwise diffeomorphic.  
\end{theorem}
We do not know that they are topologically isotopic or not. If they are topologically isotopic, then this is the first example of the exotic embedding of $S^2$ in $\C P^2$. 
Related result is given by Finashin \cite{finashin2002knotting}. He gives an infinite family of surfaces for each $k\ge 3$ whose homology class, genus, and the fundamental group of the complement of the surfaces are the same as that of algebraic surfaces with degree $2k$. 
Kim and Ruberman prove that they are topologically isotopic in \cite{kim2008topological}. 
To the best of our knowledge, the case that $k=3$ is the example of exotic orientable surfaces in $\C P^2$ with minimal genus.

\subsection{Exotic group action}
In fact our invariant is capable of detecting exotic $\Z/2$-action. If a pair of the $4$-manifolds with involutions $(\tilde{X}, \iota)$ and $(\tilde{X}', \iota')$ are $\Z/2$-equivariant diffeomorhic, then our invariant coincides. As applications, we have the following theorems:

\begin{theorem}[\cref{Z2_action_CP2}]\label{exotic_diffeo_cp2}
    At least one of the following two statements holds:
    \begin{itemize}
        \item There exist an of exotic $\C P^2$. 
        \item There exist infinite family of exotic $\Z/2$-actions on $\C P^2$. 
    \end{itemize}
\end{theorem}

\begin{theorem}[\cref{Z2_action_S4}]\label{exotic_diffeo}
At least one of the following two statements holds:
\begin{itemize}
    \item There exist a counter example to the $4$-dimensional smooth Poincar\'{e} conjecture. 
    \item There exist infinite family of $\Z/2$-actions on $S^4$ that do not preserve any positive scalar curvature metric. 
\end{itemize}
\end{theorem}

We do not know whether or not the double branched covers of $S^4$ along the $P^2$-knots and $2$-knots that we give in this paper are diffeomorphic to $\C P^2$ and $S^4$ respectively.


\subsection{Background of the Real Seiberg--Witten theory}\label{intro_Real_SW}
Several authors have used various invariants to obstruct smooth isotopies between embedded surfaces in $4$-manifolds. For example, Finashin--Kreck--Viro used the Donaldson invariant of the double branched cover of surfaces in \cite{MR0970078}. Fintushel--Stern \cite{MR1492129} used the Seiberg--Witten invariant.  Baraglia \cite{baraglia2020adjunction} used the family Seiberg--Witten invariant. To detect exotic surfaces embedded in the $4$-manifold with boundary, Juh\'{a}sz--Miller--Zemke used the Heegaard Floer homology in \cite{MR4347309}, Heyden--Sundberg used the Khovanov homology in \cite{hayden2021khovanov}, Hayden used symplectic geometry and contact geometry in \cite{hayden2020exotically}, and Lin--Mukherjee \cite{lin2021family} used the family Bauer--Furuta invariant. 

In this paper, we construct an invariant of diffeomorphism classes of the pair $(X, S)$, where $X$ is an oriented closed compact $4$-manifold and $S$ is an embedded surface in $X$ such that $H^1(X, \Z/2)=0$ and the homology class $[S]_2$ in the $\Z/2$-coefficient is $0$,i.e. $[S]_2=0 \in H_2(X, \Z/2)$. (See \cref{nice_emb_closed} and \cref{invariant_closed} for detail.)  

To construct the invariant, we use a variant of the Seiberg--Witten theory, which we call the Real Seiberg--Witten theory. A brief setup of the Real Seiberg--Witten theory as follows.  Suppose that a $3$- or $4$-manifold $X$ with an involution $\iota$ whose fixed point set is a codimension $2$-submanifold, and a spin$^c$ structure $\mathfrak{s}$ on $X$ satisfies $\iota^*c_1(\mathfrak{s})=-c_1(\mathfrak{s})$. Suppose further that there is an anti-linear lift $I$ of the involution $\iota$ to the spinor bundle which is compatible with the Clifford multiplication. We can easily seen that the Seiberg--Witten equations are equivariant under the involution $(\nabla_A, \phi) \mapsto (I \circ \nabla_A \circ I^{-1}, I\phi)$, where $A$ is a spin$^c$ connection and $\phi$ is a positive spinor. We consider solutions of the Seiberg--Witten equations that are fixed by the involution, and we call them the solution of the Real Seiberg--Witten equations. 

This setting is introduced by Tian--Wang \cite{tian2009orientability} for hermitian almost complex manifolds. Kato \cite{MR4365044} introduced Real Seiberg--Witten theory for general $4$-manifolds with an involtuion with non-empty fixed points. Real Seiberg--Witten Theory is used by Konno, Taniguchi and the author \cite{konno2021involutions,konno2023involutions}, and Li \cite{li2022monopole,li2023real}.  
When the involution is a free action, the Real theory is called $\mathrm{Pin}^-(2)$-monopole that is introduced by Nakamura \cite{MR3118618,MR3415770}. 

For $2$-knots and $P^2$-knots, we define a numerical invariant in \cref{def_degree} and \cref{def_degree_P2} that is given by counting solutions in the framed moduli space of the Real Seiberg--Witten equations. This invariant has significant property: The invariant of a surface knot given by connected sum of $2$-knots or $P^2$-knots is multiplication of the invariants of the connected components. 

In deed, we define the numerical invariant for the pair $(\Sigma, c)$, where $\Sigma$ is a $2$-knot or $P^2$-knot and $c$ is a Real $\mathrm{spin}^c$ structure on the double branched cover of $S^4$ along $\Sigma$. 
For simplicity, here we state our result in the case that we can take canonical Real $\mathrm{spin}^c$ structure $c$. 
If $\Sigma$ is a $P^2$-knot, to get the canonical Real $\mathrm{spin}^c$ structure, we assume the double branched cover along $\Sigma$ is an integer homology $\pm \C P^2$ here. 

\begin{theorem}[see \cref{multiplicative} and \cref{invariant_Unknot}]\label{main2}
    Let $S$ be a $2$-knot in $S^4$. Then we have a $\Z_{\ge 0}$-valued invariant $\lvert deg(S) \rvert$ with the following properties: 
    \begin{itemize}
        \item For every $2$-knot $S$ and $S'$, if $(S^4, S)$ and $(S^4, S')$ are pairwise diffeomorphic, then $\lvert deg(S) \rvert =\lvert deg(S') \rvert$. 
        \item For every $2$-knot $S$ and $S'$, we have $\lvert deg(S\#S') \rvert=\lvert deg(S) \rvert \lvert deg(S') \rvert$. 
        \item For the unknot $U$, we have $\lvert deg(U) \rvert=1$. 
    \end{itemize}
\end{theorem}
\begin{theorem}[see \cref{conn_sum_S_P}, and \cref{standard_P2_knot}]
    Let $P, P' \subset S^4$ be $P^2$-knots with 
    the double branched covers $\Sigma_2(S^4, P)$ and $\Sigma_2(S^4, P')$ of $S^4$ along $P$ and $P'$ are integer homology $\pm \C P^2$. Then we have a $\Z_{\ge 0}$-valued invariant $\lvert deg(P) \rvert$ that satisfies the following properties:
    \begin{itemize}
        \item 
        If $(S^4, P)$ and $(S^4, P')$ are pairwise diffeomorphic, then $\lvert deg(P) \rvert =\lvert deg(P') \rvert$.
        \item If $S$ is a $2$-knot that satisfies the double branched cover $\Sigma_2(S^4, S)$ of $S^4$ along $S$ is an integer homology $S^4$, then 
        \[
            \lvert deg(P \# S) \rvert=\lvert deg(P) \rvert \lvert deg(S) \rvert. 
        \]
        \item Let $P_{\pm}$ be the standard $P^2$-knot with $P_{\pm} \circ P_{\pm}=\pm 2$. Then $\lvert deg(P_+) \rvert=1$ and $\lvert deg(P_-) \rvert=0$. 
    \end{itemize}
\end{theorem} 

\subsection{Organization of this paper}
We explain the outline of this paper. In \cref{Involution_SW}, we recall the construction of the Real version of the Bauer--Furuta invariant and Floer homotopy type. Moreover, we prove a gluing theorem of the Real version of relative Bauer--Furuta invariants. In \cref{def_inv_surfaces}, we give the definition of the invariant of embedded surfaces which take value in non-equivariant stable homotopy group of the spheres. For $2$-knots and $P^2$-knots, we have numerical invariants. 
In \cref{Examples_of_surfaces}, we calculate the invariants of examples of surfaces in $\# M \overline{\C P^2}\# {N} \C P^2$ ($M\ge 20$ and $N \ge 4$), $2$-knots, and $P^2$-knots. We also give proofs of \cref{main0,complement_,main00,main1,exotic_diffeo_cp2,exotic_diffeo} in \cref{Examples_of_surfaces}. 

\section{Involution and Seiberg-Witten theory}\label{Involution_SW}
In this section, we recall the Real Seiberg-Witten theory. 
\subsection{Definition and classification of the real structure}
We recall the definition and the classification of the real structure. 
\begin{definition}\label{real}
Let $M$ be a $3$- or $4$-manifold and let $\iota \colon M \to M$ be an involution. Let $\mathfrak{s}$ be a $\mathrm{spin}^c$ structure and $S$ be the spinor bundle of $\mathfrak{s}$. Let $I \colon S \to S$ be an anti-linear map which covers $\iota$ and preserves $\Z/2$-grading of the spinor bundle if $\dim(M)=4$. 
If $I$ satisfies the following compatibility with the Clifford action $\rho$, we call $I$ a real structure on $\mathfrak{s}$ : 
$I(\rho(a)\phi)=\rho(\iota^* a)I(\phi)$ for all $a \in \Omega^1(M)$ and $\phi \in \Gamma(S)$. 
\end{definition}

From \cref{real}, it is easy to see that if there is a real structure which covers an involution $\iota$, we have $\iota^* \mathfrak{s}\cong \overline{\mathfrak{s}}$, where $\overline{\mathfrak{s}}$ is the reversed $\mathrm{spin}^c$ structure. 
The following classification theorem is given by Li \cite[Subsection 3.1]{li2022monopole}. 
\begin{theorem}[Li\cite{li2022monopole}]\label{classification_of_real_spinc}
    Let $\mathfrak{s}$ is a $\mathrm{spin}^c$ structure on $M$. Suppose there is a real structure on $\mathfrak{s}$ which covers an involution $\iota$. Then the isomorphism classes of Real structures on the $\mathrm{spin}^c$ structure $\mathfrak{s}$ which cover $\iota$ are identified with the set
\[
 \frac{H^1(M, \Z)^{\iota^*}}{(1+\iota^*)(H^1(M, \Z))}. 
\]
\end{theorem}
\begin{remark}
    If $H^1(M, \Z)^{\iota^*}=0$, the number of the isomorphism class of the Real structure on a $\mathrm{spin}^c$ structure $\mathfrak{s}$ is at most $1$. 
\end{remark}

\subsection{Real Seiberg--Witten Floer homotopy type}

Let us recall the construction the real Seiberg--Witten Floer homotopy type and introduce an obvious generalization to the procedure. In \cite{konno2021involutions} and \cite{konno2023involutions}, we assume the three-dimensional manifold is a rational homology sphere. However, if $\iota^*$ acts trivially on the first cohomology of the three manifold, then it is possible to define the Floer homotopy type in a same manner even if $b_1\neq 0$. 

Let $Y$ be a closed oriented $3$-manifold and $\mathfrak{t}$ be a $\mathrm{spin}^c$ structure on $Y$. Let $\iota$ be an orientation preserving involution on $Y$. Let $I$ be a real structure on $\mathfrak{t}$ which covers $\iota$. We fix a reference $\Spinc$ connection $A_0$ which satisfies $I \circ A_0 \circ I^{-1}=A_0$ and $F_{A_0^{\tau}}=0$. The existence of a connection which satisfies $I \circ A_0 \circ I^{-1}=A_0$ is shown in \cite[Proposition 2.12]{konno2023involutions} by taking the mean of the connection. From the assumption $H^1(Y, \R)^{-\iota^*}=0$ we have $H^2(Y, \R)^{-\iota^*}=0$ by Poincare duality. Then from the assumption $\iota^*c_1(\mathfrak{t})=-c_1(\mathfrak{t})$ and $H^2(Y, \R)^{-\iota^*}=0$, we have 
$c_1(\mathfrak{t})=0$ in $\R$ coefficient. We assume $F_{A_0^{\tau}}=0$ on the $\Spinc$ connection $A_0$. 

Then the space of $\Spinc$ connections $\mathcal{A}$ can be identified with $i\Omega^1(M)$ and the space of $I$-invariant $\Spinc$ connections are identified $i\Omega(M)^{-\iota^*}$, which is the space of imaginary valued $1$-forms with $\iota^*a=-a$. 

\begin{proposition}[\cite{konno2021involutions} and \cite{konno2023involutions}]\label{RSWF}
    Let $\mathfrak{t}$ be a $\mathrm{spin}^c$ structure on $Y$ and $\iota$ be an involution. We assume that $H^1(Y, \R)^{-\iota^*}=0$. Let $I$ be a real structure on $\mathfrak{t}$ which covers $\iota$. 
    Let $Coul(Y)^I$ be a slice of the configuation space fixed by the involution $(-\iota^*, I)$:
    \[Coul(Y)^I=L^2_{k+1/2}(\ker(d^*\colon i\Omega^1(Y)\to i\Omega^0(Y))^{-\iota^*} \oplus \Gamma(S)^I). \]
Let $l=(-\ast d, \Slash{D}_{A_0})$ be a Fredholm operator on $Coul(Y)^I$, where $A_0$ be a reference connection and $\Slash{D}_{A_0}$ be a $\Spinc$ Dirac operator on $S$. Let us denote $V^{\nu}_{-\nu}$ by the sum of the eigenspaces of $l$ such that the eigenvalue $\lambda$ satisfies $ \lambda \in (-\nu,  \nu]$. For a positive real number $R$, we denote $B(R)$ by the set $\{(a, \phi) \in Coul(Y)^I \mid \lVert (a, \phi) \rVert _{L^2_{k+1/2}}<R \}$, where $k>100$. 
Then there is a positive real number $R$ such that if we take $\nu$ sufficiently large, $B(R) \cap V^{\nu}_{-\nu}$ is an isolating neighborhood of the finite dimensional approximated Seiberg--Witten Flow. 
\end{proposition}
\begin{remark}
    Form the assumption that $-\iota^*$ fixed part of $H^1(Y, \R)$ is trivial, we see that 
    \[
        \ker(d^*\colon i\Omega^1(Y)\to i\Omega^0(Y))^{-\iota^*} = d^*(i\Omega^+(Y))^{-\iota^*}. 
    \]
    Therefore, when we have $b^1(Y)=b^1_{\iota}(Y)$, we do not have to consider parametrized spectrums on the Picard torus even if the first betti number $b^1(Y)$ is not $0$. 
\end{remark}

Before we prove \cref{RSWF}, we recall the Manolescu's finite dimensional approximation of the Seiberg--Witten Flow with the Real involution. 

Let us denote $\xi \colon i\Omega^1(Y) \to i\Omega^0(Y)$ be the composition of $d^*$ and the inverse of the Laplacian $(d^* d)^{-1} \colon i\Omega^0_0(Y) \to i\Omega^0_0(Y)$, where $\Omega^0_0(Y)=\{ f \in \Omega^0(Y) \mid \int_Y f d\mu_g=0\}$. 
Let $Q_{\mu} \colon Coul(Y) \to Coul(Y)$ be a non-linear term of the Seiberg--Witten flow on the slice $Coul(Y)$:
\[
    Q_{\mu}(a, \phi)=(q(\phi)-d\xi(q(\phi))+d^*\mu, \rho(a)\phi-\xi(q(\phi))\phi)
\]
where $q$ is a bilinear form $\Gamma(S)\otimes_{\R} \Gamma(S) \to i\Omega^1(Y)$ which satisfies $\langle q(\phi, \psi), a \rangle=\text{Re}(\langle \rho(a)\phi, \psi \rangle)$ and $\mu \in i\Omega^2(Y)$. 

Note that $\xi$ and $q$ are equivariant under the real involution $-\iota^*$ and $I$. If we have $-\iota^*\mu =\mu$, we see that $Q_{\mu}$ is equivariant under the real involution. 

Let $\tilde{p}^{\nu}_{-\nu}$ be the $L^2$ orthogonal projection to the finite dimensional subspace $V^{\nu}_{-\nu}$ and $\beta \colon [0, 1] \to [0, 1]$ be a smooth function such that $\beta(\tau)=0$ if $\tau$ is in the neighborhood of $\{0, 1\}$ and $\int_{0}^{1} \beta(\tau) d\tau=1$. 
We set 
\[
    {p}^{\nu}_{-\nu}:=\int_{0}^{1} \beta(\tau) \tilde{p}^{\nu-\tau}_{-\nu+\tau} d\tau. 
\]
We can easily check that $p^{\nu}_{-\nu}$ is $I$-equivariant. 

Let us consider the flow on $V_{-\nu}^{\nu}$ that is given by

\[
    \frac{d\gamma}{dt}=-l\gamma - p_{-\nu}^{\nu} \circ Q_{\mu}(\gamma)
\]
for $\gamma \colon \R \to V_{-\nu}^{\nu}$. 
We call this flow by the finite dimensional approximated Seiberg--Witten flow. 
\begin{proof}[Proof of \cref{RSWF}]
From the assumption that $H^1(Y, \R)^{-\iota*}=0$, we have $H^2(Y, \R)^{-\iota^*}=0$ from the Hodge decomposition and the fact that Hodge star operator commutes with $\iota^*$. Thus $c_1(\mathfrak{t})$ is a torsion in $H^2(Y, \Z)$ since $\iota^*c_1(\mathfrak{t})=-c_1(\mathfrak{t})$ and the Chern--Simons--Dirac functional is balanced. 
From the proof of \cite[Proposition 10]{MR3211454}, it is shown that there exist a positive Real number $R$ such that for all Seiberg--Witten trajectory $\gamma \colon \R \to Coul(Y)$ with finite energy and for all $t \in \R'$, there is a gauge transformation $u \colon Y \to U(1)$ which satisfies 
\begin{align}\label{finiteness}
    \lVert u \cdot \gamma(t) \rVert _{L^2_{A_0, k}(Coul(Y))}\le R. 
\end{align}
Note that we do not have to have the transversality of the trajectories of flows when the Chern--Simons--Dirac functional is balanced, whereas the transversality is assumed in \cite[Proposition 10]{MR3211454}. 

From the assumption that $H^1(Y, \Z)^{-\iota^*}=0$, if $\gamma$ is invariant under the action of $(-\iota^*, I)$, the trivial gauge transformation $1$ satisfies \eqref{finiteness}. 

From \cite[Corollary 3.7]{MR3784516}, we see that if $\nu >>1$ and $R'<R$, $B(R)\cap V^{\nu}_{-\nu}$ is an isolated neighborhood. 
\end{proof}
Since the $(-\iota^*, I)$-invariant part of the finite dimensional approximated Seiberg--Witten flow is equivariant under the action of constant gauge transformation in Real theory, that is, $\pm 1$. 
Moreover, if the $\mathrm{spin}^c$ structure $\mathfrak{t}$ comes from a spin structure on $Y$, we have an additional symmetry. In that cese we have a specific Real structure $I_0=\tilde{\iota} \circ j$, where $\tilde{\iota}$ is a lift to the spinor bundle of the spin structure of order $4$. We see that $(-\iota^*, I_0)$ commutes $(-1, j)$-action on $Coul(Y)$. 

Thus we can define the stable $G$-stable spectrum valued invariant $SWF_G(Y, \mathfrak{t}, I)$ defined as in \cite[Definition 3.11.]{konno2023involutions}. 

\begin{proposition}\label{conley}
    Let $G=\Z/2$ or $\Z/4$ and $\mathbb{K}_G$ be a representation space 
    \[
        \mathbb{K}_G= \begin{cases}
            \tilde{\R} & \textnormal{if} \; G=\Z/2 \\
            \R^2_1 & \textnormal{if} \; G=\Z/4 
        \end{cases}
    \]
   , where $\tilde{\R}$ is the sign representation of $\Z/2$ and, where $\R^2_1$ is given by regarding the rank $1$ complex representation of $\Z/4$ with weight $1$ as a rank $2$ real representation of $\Z/4$.
    Let $\mathfrak{C}_{G}$ be the set of $G$-stable equivalence class that defined in \cite[Subsection 3.4]{konno2023involutions}. 
    Let $I^{\nu}_{-\nu}(Y, \mathfrak{t}, \iota, I, g)$ be a Conley index of the isolated neighborhood $B(R) \cap V^{\nu}_{-\nu}$ of the finite dimensional approximated Seiberg--Witten flow. 
    We denote $n(Y, \mathfrak{t}, g)$ by the collection term given in \cite[Subsection 3.4, (16)]{konno2023involutions}. 
    Then the triple 
    \[
        (I^{\nu}_{-\nu}(Y, \mathfrak{t}, \iota, I, g), \dim(V^{0}_{-\nu}((i\ker d^*)^{-\iota^*})), \dim_{\mathbb{K}_G}(V^{0}_{-\nu}(\Gamma(S)^I))+\frac{n(Y, \mathfrak{t}, g)}{2})
    \]
    defines an element of $\mathfrak{C}_G$. 
    Then the $G$-stable equivalence class of the above triple is independent of the choice of the $\iota$-invariant metric $g$, the number $\nu$, $R$, and the index pair of the isolated neighborhood of $B(R) \cap V^{\nu}_{-\nu}$. 
\end{proposition}
Let us denote $V^{\lambda}_{-\nu}((i\ker d^*)^{-\iota^*})$ by $V^{\lambda}_{-\nu}$ and denote $V^{\lambda}_{-\nu}(\Gamma(S)^I)$ by $W^{\lambda}_{-\nu}$. 
\begin{definition}[see also \cite{konno2021involutions}, \cite{konno2023involutions}]\label{real_SWF}
   Let $Y$ be a closed oriented $3$-manifold and $\iota$ be an involution. Let $(\mathfrak{t}, I)$ be a Real $\mathrm{spin}^c$ structure that covers $\iota$. Then we define a $\mathfrak{C}_G$ valued invariant
    \[
        SWF_G(Y, \mathfrak{t}, I) : = [(I^{\nu}_{-\nu}(Y, \mathfrak{t}, \iota, I, g), \dim(V^{0}_{-\nu}), \dim_{\mathbb{K}_G}(W^{0}_{-\nu})+\frac{n(Y, \mathfrak{t}, g)}{2})]
    \]
   , where $G=\Z/4$ if $\mathfrak{t}$ comes from a spin structure and $I$ commutes with $j$, and $G=\Z/2$ otherwise. 
\end{definition}
\begin{remark}
    For simplicity, we sometimes write 
    \[
        SWF_G(Y, \mathfrak{t}, I):=\Sigma^{-V^{0}_{-\nu}}\Sigma^{-W^{0}_{-\nu}-\mathbb{K}_G^{n(Y, \mathfrak{t}, g)/2}}I^{\nu}_{-\nu}(Y, \mathfrak{t}, \iota, I, g). 
    \]
    Note that if $H_1(Y, \Z)^{-\iota^*} \neq 0$, $n(Y, \mathfrak{t}, g)$ may be in $\Q \setminus \Z$. This expression is just a formal notation. 
\end{remark}

\subsection{Real Bauer--Furuta invariant and its gluing}
Here we recall the construction of the Real Bauer--Furuta invariant given in \cite[Section 3.4-Section 3.7]{konno2023involutions}. 
\subsubsection{Definition}
    Let $W$ be a compact oriented four manifold with boundary $Y$ with involution $\iota$. 
    We assume that $H^1(Y, \R)^{-\iota^*}=H^1(W, \R)^{-\iota^*}=0$. Let $\mathfrak s$ be a $\mathrm{spin}^c$ structure on $W$ and let $I$ be a real structure on $\mathfrak{s}$. Then we have 
    a finite dimensional approximation of the relative $I$-invariant Seiberg--Witten map 
    \begin{align*}
        (\text{SW}^I, \Pi^{\nu} \circ r) \colon & L^2_{k+1}(i\Omega^1_{CC}(W))^I \times L^2_{k+1}(\Gamma(S^+))^I  \\
       & \to L^2_k(i\Omega^+(W))^I \times L^2_{k}(\Gamma(S^-))^I \times (\ker d^*_Y)^I \times (H^I)^{\nu}_{-\infty} 
    \end{align*}
   , where $r \colon i\Omega^1_{CC}(W) \oplus \Gamma(S^+) \to Coul^I$ is a restriction map and $\Pi^{\nu} \colon Coul(Y)^I \to (\ker d^*_Y)^I \times (H^I)^{\nu}_{-\infty}$ be a $L^2$ orthogonal projection.  
    
    Let us define the infinite representation space $\mathcal{V}$ and $\mathcal{W}$ given as follows: We set $\mathcal{V}=\oplus_{\N}\R$ if $G=\Z/2$ and $\mathcal{V}=\oplus_{\N}\tilde{\R}$ if $G=\Z/4$, and we set $\mathcal{W}=\oplus_{\N} \mathbb{K}_G$.  
We identify the Hilbert space $L^2_{k+1}(i\Omega^1_{CC}(W))^I \oplus L^2_{k+1}(\Gamma(S^+))^I$ with the Hilbert space $\bar{\mathcal{V}} \oplus \bar{\mathcal{W}}$, the completion of $\mathcal{V} \oplus \mathcal{W}$ by $l^2$ norm. Note that the choice of the identification is contractible. We also identify $L^2_k(i\Omega^+(W))^I \times L^2_{k}(\Gamma(S^-))^I$ with $\bar{\mathcal{V}} \oplus \bar{\mathcal{W}}$ in the same manner. 

Let us take the set $\mathcal{U}$ of the finite dimensional subspace given by direct sum of finite dimensional subspaces $V \subset \mathcal{V}$ and $W \subset \mathcal{W}$. For $U \in \mathcal{U}$ we have a functor $\Sigma^{\pm U} \colon \mathfrak{C}_G \to \mathfrak{C}_G$ given as in the same manner in \cite[Section 6.]{manolescu2003seiberg} and \cite[Section 4.1, Lemma 4.1, Proposition 4.2]{MR3784516}. Note that the functor $\Sigma^{\pm U}$ is well defined up to homotopy even if $U$ is not a standard Euclidean space. 
By taking the finite dimensional approximation of the map $(\text{SW}^I, \Pi^{\nu} \circ r)$, we have a sequence of a map
\[
    f_n \colon (\ker{(D SW^I_0, \Pi^0 \circ r)} \oplus U_n)^+ \to \Sigma^{\coker{(D SW^I_0, \Pi^0 \circ r)}  \oplus U_n}\Sigma^{-V^0_{\nu}}I^{\nu}_{-\nu}(Y, \mathfrak{t}, \iota, I, g)
\]
for $n \in \Z_{>0}$ and $\nu \in \R_{>0}$. The object in $\mathfrak{C}_G$ given by 
\[
    V_W^+ \colon = [(\ker{\Slash{D}_{\tilde{g}}}^+, 0, \dim_{\mathbb{K}_G}(\coker{\Slash{D}_{\tilde{g}}^+})+\frac{n(Y, \mathfrak{t}, g)}{2})]
\]
here $\tilde{g}$ be a metric on $W$ with $\tilde{g}|_Y=g$ and $\Slash{D}_{\tilde{g}}^+$ is a Real $\Spinc$ Dirac operator with Atiyah--Patodi--Singer boundary condition with respect to the metric $\tilde{g}$ and the boundary operator $\Slash{D}_{A_0}$, the reference Real $\Spinc$ connection on $Y$. From \cite[Lemma 3.22.]{konno2023involutions} the object $V_W^+$ is independent of the choice of the metrics $\tilde{g}$ and $g$ and the reference connection $A_0$ since $\ind_{\R}(\Slash{D}_{\tilde{g}}^+)-\dim_{\R}(\mathbb{K}_G)n(Y, \mathfrak{t}, g)/2 =\frac{c_1(\mathfrak{s})^2-\sigma(W)}{8}$ here we use the term $\frac{c_1(\mathfrak{s})^2-\sigma(W)}{8}$ in the same meaning of Manolescu's paper \cite[Section 6.]{manolescu2003seiberg}. Then $\{f_n\}_n$ defines a sequence of the morphism in the category $\mathfrak{C}_G$. Hence we have the following relative invariant. 

\begin{definition}\label{definition_of_the_relative_bauer_furuta}
    Let $W$ be a compact oriented four manifold with boundary $Y$ with involution $\iota$. Let $(\mathfrak{s}, I)$ be a Real $\mathrm{spin}^c$ structure on $W$. We assume $(W, \iota, \mathfrak{s}, I)$ satisfies the assumptions as above. 
Then we have tha stable $G$-equivariant homotopy invariant:
\begin{align*}
    \Psi(W, \mathfrak{s}, I) \in \text{colim}_{U \in \mathcal{U}}[\Sigma^U (V_W)^+ \to \Sigma^{H^+(W)^{-\iota^*} \oplus U}(\text{SWF}(Y, \mathfrak{t}, I)) ]_G .
\end{align*}
where $(V_W)^+$ be an object in $\mathfrak{C}_G$ given as above. We call this invariant by Real Bauer--Furuta invariant of the Real structure $(\mathfrak{s}, I)$ on $X$. 
\end{definition}
Let us denote \[\text{colim}_{U \in \mathcal{U}}[\Sigma^U (V_W)^+ \to \Sigma^{H^+(W)^{-\iota^*} \oplus U}(\text{SWF}(Y, \mathfrak{t}, I)) ]_G\] by $\pi_{V_W}^{st}(\Sigma^{H^+(W)^{-\iota^*}}\text{SWF}(Y, \mathfrak{t}, I))$. 
\begin{remark}
    To define the relative Real Bauer--Furuta invariant, we identify the Hilbert spaces of the domain and target of the Real Seiberg--Witten map to the standard Hilbert spaces with $G$-action. Form the Kuiper's theorem, the choice of the identification is contractible so that we can define the invariant as above. 
    However, there is no canonical trivialization of the virtual vector space $V_W$. In the ordinary Seiberg--Witten theory, we can take a  trivialization of $V_W$ using the complex structure of the spinor bundle, this is unique up to homotopy since $\mathrm{U}(m)$ are connected for $m \ge 0$. In our case, $\mathrm{O}(n)$ has two connected components so that there are two trivializations of $V_W$ up to homotopy. 
    In this paper, we only consider the numerical invariants defined by using $\Psi(W, \mathfrak{s}, I)$ that does not depends on the choice of the trivialization of $V_W$. 
\end{remark}
\begin{remark}
    In this paper we do not use the group action of $G$ when we apply the Real Bauer--Furuta invariant to detect exoticness on the embedded surfaces. 
\end{remark}

\subsubsection{Statement and proof of the Gluing theorems}

We show that the following gluing theorem:

\begin{theorem}\label{gluing_general}
    Let $X_0$ and $X_1$ be compact oriented $4$-manifolds with boundary $\partial X_0=Y$ and $\partial X_1=-Y$. ($Y$ may be disconnected. )
    Let $\iota_0$, $\iota_1$ and $\iota'$ are involutions on $X_0$, $X_1$ and $Y$ respectively. We assume that $\iota_0|_Y=\iota_1|_Y=\iota'$. 
    The $4$-manifold that is given by $X_0 \cup_Y X_1$ is denoted by $X$. Let $\iota=\iota_{X_0} $
    We assume $H^1(X, \R)^{-\iota^*} \cong H^1(X_0, \R)^{-\iota_0^*} \cong H^1(X_1, \R)^{-\iota_1^*} \cong H^1(Y, \R)^{-\iota^*}=0$, $H^0(X, \Z)^{-\iota^*}=0$, and $H^0(X_0, \Z)^{-\iota^*}=0$. We also assume that $X/\iota$, $X_0/\iota_0$, $X_1/\iota_1$, and $Y/\iota'$ are connected. 
    Then following statements holds:
    \begin{itemize}
        \item For each Real $\mathrm{spin}^c$ structure $(\mathfrak{s}, I)$ on $X$ that covers $\iota$, there is a unique Real $\mathrm{spin}^c$ structure $(\mathfrak{s}_i, I_i)$ on $X_i$ that covers $\iota_i$ for $i=0, 1$ such that $(\mathfrak{s}, I)|_{X_i}=(\mathfrak{s}_i, I_i)$. 
        \item Let $\epsilon'$ be the Spaniel--Whitehead duality map in $\mathfrak{C}_G$. Then there is a map 
        \[
          \epsilon \colon  \pi_{V_{X_0}}^{st}(\Sigma^{H^+(X_0)^{-\iota^*}}SWF(Y, \mathfrak{t}, I))\times \pi_{V_{X_1}}^{st}(\Sigma^{H^+(X_1)^{-\iota^*}}SWF(-Y, \mathfrak{t}, I)) \to \pi_{V_{X}}^{st}(S^{H^+(X)^{-\iota^*}})
        \]
        given by composing smash product and the map $\epsilon$. Then
         we have the gluing formula for the Real Bauer--Furuta invariant: 
        \begin{align}\label{gluing_formula_general}
            \Psi(X, \mathfrak{s}, I)=\epsilon(\Psi(X_0, \mathfrak{s}_0, I_0)\wedge \Psi(X_1, \mathfrak{s}_1, I_1))
        \end{align}
    \end{itemize}
\end{theorem}
\begin{remark}
    For simplicity, we state the statement when the resulting $4$-manifold $X$ is closed. However if there are boundaries in $X_0$ and $X_1$ other that $\pm Y$, we have similar gluing result under some assumption. 
\end{remark}
From \cref{gluing_general}, we have two corollaries. 
\begin{corollary}\label{gluing_conn}
    Let $W_0$ and $W_1$ be connected compact oriented $4$-manifolds with boundary $\partial{W_0}= Y$ and $\partial{W_1}=-Y$. 
We assume $Y$ is connected and $H^1(W_i, \Z)^{-\iota^*}=0$ for $i=0, 1$. 
Let $\iota_i$ be an involution on $W_i$ and $(\mathfrak{s_i}, I_i)$ be a Real $\mathrm{spin}^c$ structure on $W_i$. We assume that $(\mathfrak{s}_0, I_0) |_{Y} \cong (\mathfrak{s}_1, I_1)|_Y$ and $H^1(\partial W_i, \Z)^{-\iota_i^*}=0$ for $i=0, 1$. 
    Then we have the gluing formula of the Real relative Bauer--Furuta invariant: 
    \begin{align}\label{gluing_formula_conn}
        \Psi(W_0 \cup_Y W_1, \mathfrak{s}_0 \cup_\mathfrak{t} \mathfrak{s}_1, I_0 \cup_{I'} I_1)=\epsilon(\Psi(W_0, \mathfrak{s}_0, I_0)\wedge \Psi(W_1, \mathfrak{s}_1, I_1))
    \end{align}    
   , where $\epsilon$ is the Spaniel--Whitehead duality map. 
\end{corollary}
\begin{proof}
From the assumption of the statement, we have
\[
    H^1(W_0 \cup_Y W_1, \R)^{-\iota^*}\cong H^1(W_0, \R)^{-\iota_0^*} \oplus H^1(W_1, \R)^{-\iota_1^*} = 0
\]
because if we we take $-\iota^*$ fixed point part of the Mayer--Vietoris exact sequence for $(W_0 \cup_Y W_1, W_0, W_1, Y)$, $H^0(Y, \R)^{-\iota^*}=0$ and $H^1(\partial W_i, \Z)^{-\iota_i^*}=0$ for $i=0, 1$ implies $H^1(W_0 \cup_Y W_1, \R)^{-\iota^*}=0$. 
\end{proof}    

\begin{corollary}\label{gluing_disconn}
    Let $Y$ be a rational homology $3$-sphere and $X$ be a compact oriented $4$-manifold with $\partial X=Y$ and $H^1(X, \Z)=0$. 
    Let $W$ be a compact oriented $4$-manifold with $\partial W=-Y \sqcup Y$. Let $\iota$ be an involution on $W$ which switches the two connected components of the boundary. We assume $H^1(W, \Z)^{-\iota^*}=0$. 
    Let $(\mathfrak{s}, I)$ be a Real $\mathrm{spin}^c$ structure on $W$ which covers $\iota$. Let $\mathfrak{t}$ be a $\mathrm{spin}^c$ structure on $Y$ induced by $\mathfrak{s}$. 
    Note that the $\mathrm{spin}^c$ structure on another component of the boundary induced by $\mathfrak{s}$ is $\overline{\mathfrak{t}}$. 
    We assume that there is a $\mathrm{spin}^c$ structure $\mathfrak{s}_X$ on $X$ such that $\mathfrak{s}_X |_{Y}\cong \mathfrak{t}$. 
    Then we can extend the Real structure $(\mathfrak{s}, I)$ on $W$ to a Real structure on $X \cup _Y W \cup_Y X$ and we have the following gluing formula:
    \begin{align}\label{gluing_formula_disconn}
        \Psi(X \cup _Y W \cup_Y X, \overline{\mathfrak{s}_X} \cup_{\overline{\mathfrak{t}}} \mathfrak{s} \cup_{\mathfrak{t}} \mathfrak{s}_X, I)=\epsilon (\Psi_{\Z/2}(X, \mathfrak{s}_X) \wedge \Psi(W, \mathfrak{s}, I))
    \end{align}
   , where $\Psi_{\Z/2}(X, \mathfrak{s}_X)$ be the image of the forgetful map \[\pi_{*}^{U(1)}(SWF(Y, \mathfrak{t})) \to \pi_{*}^{\Z/2}(SWF(Y, \mathfrak{t}))\] of the usual relative Bauer--Furuta invariant of $(X, \mathfrak{s}_X)$. 
\end{corollary}
\begin{proof}
    We take $-\iota^*$ invariant part of the Mayer--Vietoris exact sequence for $X \cup _Y W \cup_Y X, W, X \sqcup X, Y \sqcup Y$ implies 
    \[ 
        H^1(X \cup _Y W \cup_Y X, \R)^{-\iota^*}\cong H^1(W, \R)^{-\iota_0^*} \oplus H^1(X \sqcup X, \R)^{-\iota_1^*}
        \]
        because $H^0(Y \sqcup Y, \R)^{-\iota^*} \cong H^0(X \sqcup X, \R)^{-\iota^*} \cong \R$. 
        It is easy to check that $\Psi(X \sqcup X, \mathfrak{s}_X \sqcup \overline{\mathfrak{s}_X}, I)=\Psi_{\Z/2}(X, \mathfrak{s}_X)$. 
\end{proof}

\subsubsection{Proof of \cref{gluing_general}}
We prove \cref{gluing_general} imitating the proof of gluing theorem of the relative Bauer--Furuta invariant presented in \cite[Section.6]{MR4602726}. Before starting the proof, we introduce some notations and recall some propositions which is essentially proved in \cite[Section.6]{MR4602726}. 

Let $X_0$ and $X_1$ be a $4$-manifolds and $Y$ be a three manifold that satisfies the assumptions in \cref{gluing_general}. 
We use the following notations. 

\underline{Notations}

\begin{itemize}
    \item $V_{X_i}^I:=L^2_{k+1}(i\Omega^1(X_i)^{-\iota_i^*} \oplus \Gamma(S_i^+)^{I_i}).$
    \item $U_{X_i}^I:=L^2_k(i\Omega^0(X_i)^{-\iota_i^*} \oplus i\Omega^+(X_i)^{-\iota_i^*} \oplus \Gamma(S^-_i)^{I_i}).$
    \item $V_X^I:=L^2_{k+1}(i\Omega^1(X)^{-\iota^*}\oplus \Gamma(S^+_X)^I).$
    \item $U_X^I:=L^2_{k}(i\Omega^0(X)^{-\iota^*} \oplus i\Omega^+(X)^{-\iota^*} \oplus \Gamma(S^-_X)^I).$
    \item $H^I:=L^2_{k+1/2}(i(d\Omega^0(Y))^{-\iota^*}\oplus \Omega^0(Y)^{-\iota^*})$.
    \item $V_Y^I:=L^2_{k+1/2}(Coul(Y)^I) \oplus H^I$.
\end{itemize}

We have describe the Spaniel--Whitehead duality map $\epsilon$. Imitating the proof of \cite[Proposition 6.10]{MR4602726}, we have the following proposition. 
\begin{proposition}\label{representative_e}\cite[Proposition 6.10]{MR4602726}
    Let $\{V_n^{i, I}\}_n, \{U_n^{i, I}\}_n$ and $\{V_n^I\}_n$ be finite dimensional approximations of $V_{X_i}^I, U_{X_i}^I$ and $V_Y^I$ respectively. We assume that for all $n$ there exist $\nu_n$ such that $V_n^I=(V^{\nu_n}_{-\nu_n})^I$ and $(SW, p^{\nu_n}_{\infty} \circ r)^{-1}(U^i_n \oplus V_n)=V_n^i$, where $p^{\nu_n}_{\infty}$ is the orthogonal projection on $V_Y^I$ to $(V^{\nu_n}_{\infty})^I$ with respect to the $L^2$ inner product. 
    The morphism $\epsilon(\Psi(X_0, \mathfrak{s}_0, I_0)\wedge \Psi(X_1, \mathfrak{s}_1, I_1))$ is given by the following representative:
    \begin{align}\label{glued_BF}
 \begin{array}{r@{\,\,}c@{\,\,}c@{\,\,}c}
 &B(V_n^{0, I}, R_0)^+\wedge B(V_n^{1, I}, R_1)^+&\longrightarrow&B(U^{0, ,I}_n, \epsilon)^+ \wedge B(U_n^{1, I}, \epsilon)^+ \wedge B(V^I_n, \bar{\epsilon})^+ \\
 &\rotatebox{90}{$\in$}&&\rotatebox{90}{$\in$}\\
 &[(x_0, x_1)]&\longmapsto&[(\tilde{SW}_n(x_0), \tilde{SW}_n(x_1), r_2(x_0)-r_2(x_1))]
 \end{array}
\end{align}
where $R_0$ and $R_1$ are sufficiently large real number and $\epsilon$ and $\bar{\epsilon}$ are sufficiently small positive real number. 
\end{proposition}
\begin{proof}
    Let $(\tilde{N}, \tilde{L})$ be an index pair of a dynamical system on a finite dimensional vector space $V$ that is given by the isolating block (\cite[Definition 3.20]{MR4602726}). 
    In this case, from \cite[Section 4.4]{MR4602726}, the duality map $\epsilon \colon \tilde{N}/\tilde{L} \wedge \tilde{N}/\tilde{\bar{L}} \to B_\delta(V)^+$ is given as follows:
    \[
        \epsilon([x]\wedge[y])=\begin{cases}
            &[a(x)-b(y)] \;\text{if} \; \lVert a(x)-b(y) \rVert \le \delta,\\
            & * \;\; \text{otherwise}
        \end{cases}
    \] 
   , where $a, b \colon \tilde{N} \to \tilde{N}$ be a map that satisfies $\lVert a(x)-x \rVert\le 2\delta $, $\lVert b(y)-y \rVert \le 2\delta$, $a \colon \tilde{N} \to \tilde{N}\setminus \nu_{\delta}(\tilde{\bar{L}})$, $b \colon \tilde{N} \to \tilde{N}\setminus \nu_{\delta}(\tilde{{L}})$, $a(x)=b(x)=x$ for $x \in \tilde{N}\setminus \nu_{3\delta}(\partial \tilde{N})$, and homotopic to the identity. 
    In ordinary case \cite[Section 6.6]{MR4602726}, the composition of the map $\epsilon$ with the smash product of relative Bauer--Furuta maps are deformed into the map that corresponds to \eqref{glued_BF} in the statement. One can check that all the homotopy and construction of the pre-index pairs and index pairs work in our context. 
\end{proof}

\begin{definition}[\cite{MR2312050}]\label{admissible_pair}
    Let $E_1, E_2$ be Hilbert spaces and let $\lVert \cdot \rVert _i$ be the norm of $E_i$ for $i=1, 2$. Let $\lvert \cdot \rvert _1$ be a weaker norm on $E_1$ and we denote the completion of $E_1$ of the norm 
 $\lvert \cdot \rvert _1$ by $\bar{E}_1$. We assume that if $\{x_n\}_n$ be a weakly convergent sequence in $E_1$, then it strongly converges in $\bar{E}_1$. 
 \begin{itemize}
    \item Let $l \colon E_1 \to E_2$ be a bounded linear map.
    \item Let $c \colon E_1 \to E_2$ be a continuous bounded map that extends to a continuous map $\bar{c} \colon \bar{E}_1 \to E_2$. 
\end{itemize}
Then we call $(l, c)$ is an admissible pair. 
\end{definition}
\begin{definition}\label{SWC_triple}
    Let $(l, c)$ be an admissible pair. If $(l, c)$ satisfies the following properties, we call $(l, c)$ SWC-pair:
    \begin{itemize}
        \item $l \colon E_1 \to E_2$ is a Fredholm map. 
        \item There exist $M' >0$ such that for all $x \in E_1$ that satisfies $lx+c(x)=0$, then we have $\lVert x \rVert _1 \le M'$.  
    \end{itemize}
    We call two SWC-pairs $(l_i, c_i)$ for $i=0, 1$ is c-homotopic if there exists a continuous $1$-parameter family $(l_t, c_t)$ of SWC-pairs that connects $(l_0, c_0)$ and $(l_1, c_1)$. 
    We call two SWC-pairs $(l_i, c_i)$ for $i=0, 1$ with possibly different domains and targets are called stably c-homotopic if there are Hilbert spaces $F_1$ and $F_2$ such that $(l_0\oplus id_{F_1}, c_0 \oplus 0_{F_1})$ is c-homotopic to $(l_1\oplus id_{F_2}, c_1 \oplus 0_{F_2})$. 
\end{definition}
\begin{remark}
    In the original paper \cite{MR2312050}, Manolescu introduced the notion of \textit{Conley type} and Lin--Sasahira--Khandhawit call it \textit{SWC-triple} in \cite{MR4602726}. Here we assume that the resulting $4$-manifold $X$ is closed so that we state the definition of the notion of the SWC-triple only in the closed case and we drop the data of the restriction map to the boundary of $X$ from the definition of the SWC-triple. Then we call it SWC-pair since there are only two data in our setting. 
\end{remark}
The following lemma is proved in \cite[Observation 1.]{MR2312050}. However, in our context, this is obvious. 
\begin{lemma}\label{stable_SWC}
    Let $(l, c)$ be an admissible pair and $g \colon E_1 \to E_3$ be a bounded surjective linear map. Then $(l\oplus g, c\oplus 0_{E_3})$ is SWC-pair if and only if $(l |_{\ker{g}}, c|_{\ker{g}})$ is SWC-pair. Moreover, if they are SWC-pair, they are stably c-homotopic. 
\end{lemma}

The reason for introducing the notion of SWC-pair is to describe the stable homotopy of the Bauer--Furuta invariant. The following lemma is useful. The proof is same as the construction of the Bauer--Furuta invariant. 
\begin{lemma}\label{same_BF}
    Let $(l, c)$ be a SWC-pair. Let $\{E_2^n\}_n$ be an increasing sequence of finite dimensional subspaces of $E_2$ such that the orthogonal projection $p_n \colon E_2 \to E_2^n$ weakly converges to the identity. Let us denote $l^{-1}(E^n_2)$ by $E_1^n$. Then there exist a sufficiently small real number $\epsilon_0>0$, a sufficiently large real number $R_0>>1$, and a sufficiently large natural number $N_0>>1$ such that if $n >N_0$, $\epsilon < \epsilon_0$, and $R>R_0$, then 
    \[
        f_n = p_n \circ (l+c)\colon E_1^n \to E_2^n
    \]
    defines a stable homotopy class 
    \[
      h_{(l, c)}= \{ f'_n \colon B(E^n_1, R)/S(E^n_1, R) \to B(E^n_2, \epsilon)/S(E^n_2, \epsilon) \}_n \in \text{colim}_{U  \in \mathcal{U}}[(\ind_{\R}(l) \oplus U)^+ \to U^+]_{\Z/2}
    \]
   , where $\mathcal{U}$ is the universe of $\Z/2$-representation spaces we defined in \cref{definition_of_the_relative_bauer_furuta}. (We regard $\ind_{\R}l$ as an virtual vector space given by $\ind_{\R}l=V_0-V_1$ for some $V_0, V_1\mathcal{U}$.)
Then $h_{(l, c)}$ does not depends on the choice of $N, \epsilon_0$ and $R_0$. Moreover, if $(l, c)$ and $(l', c')$ are stably c-homotopic, then $h_{(l, c)}=h_{(l', c')}$. 
\end{lemma}

Now we can prove the gluing theorem. 
\begin{proof}[Proof of \cref{gluing_general}]
    Recall that the Real Seiberg--Witten map on $X$ is given by 
    \[
        \begin{array}{r@{\,\,}c@{\,\,}c@{\,\,}c}
            SW_X \colon & L^2_{k+1}(i(\ker(d^*\colon \Omega^1(X) \to \Omega^0(X)))^{-\iota^*} \oplus \Gamma(S^+_X)^I)&\longrightarrow&L^2_k(i\Omega^+(X)^{-\iota^*}\oplus \Gamma(S^-_X)^I) \\
            &\rotatebox{90}{$\in$}&&\rotatebox{90}{$\in$}\\
            &(a, \phi)&\longmapsto&(d^+a - \rho^{-1}((\phi \phi^*)_0), \Slash{D}_{A_0}\phi + \rho(a)\phi)
        \end{array} 
    \]
    and if we set $l_X=(d^+, \Slash{D}_{A_0}), c_X(a, \phi)=(- \rho^{-1}((\phi \phi^*)_0), \rho(a)\phi)$, $E_1=L^2_{k+1}(i(\ker(d^*\colon \Omega^1(X) \to \Omega^0(X)))^{-\iota^*} \oplus \Gamma(S^+_X)^I)$, and $E_2=L^2_k(i\Omega^+(X)^{-\iota^*}\oplus \Gamma(S^-_X)^I)$ we can easily check that $(l_X, c_X)$ is a SWC-pair. 
    From \cref{representative_e} and \cref{same_BF}, what we only have to do is to prove the map 
    \begin{align*}
            & L^2_{k+1}(i\Omega^1_{CC}(X_0)^{-\iota^*} \times (\Gamma(S^+_{X_0}))^I) \times  L^2_{k+1}(i\Omega^1_{CC}(X_0)^{-\iota^*} \times (\Gamma(S^+_{X_0}))^I) \to 
            \\ &L^2_k(i\Omega^+(X_0)^{-\iota^*} \times \Gamma(S^-_{X_0})^I) \times L^2_k(i\Omega^+(X_1)^{-\iota^*} \times \Gamma(S^-_{X_1})^I) \times L^2_{k+1/2}(Coul(Y)^I) , 
        \end{align*}
        \[
            (x_0, x_1) \mapsto (SW_{X_0}(x_0), SW_{X_1}(x_1), r_0(x_0) - r_1(x_1))
        \]
    is a SWC-pair and is stably c-homotopic to the Real relative Bauer--Furuta map of $W$, where ${SW}_{X_i}(a, \phi)=(d^+a - \rho^{-1}((\phi \phi^*)_0), \Slash{D}_{A_0}\phi + \rho(a)\phi)$ is the Seiberg--Witten map and $r_i(a, \phi)=(a|_Y, \phi|_Y) \in Coul(Y)^I$ be a restriction map. We denote $(d^+, \Slash{D}_{A_0}) \colon L^2_{k+1}(i\Omega^1(X_i)^{-\iota^*} \oplus \Gamma(S^+_{X_i})^I) \to L^2_{k}(i\Omega^+(X_i)^{-\iota^*} \oplus \Gamma(S^-_{X_i})^I)$ by $l_{X_i}$ for $i=0, 1$. 

    We imitate the method of \cite[Section 6.4]{MR4602726}. There are $7$ steps to deform the relative Bauer--Furuta map. Let us check that each step works in our situation. The main differences between the ordinary case \cite[Section 6.4]{MR4602726} and the our case only are in Step.4 and Step.5. In other cases, we can follow the argument in \cite[Section 6.4]{MR4602726} except for changing notations. 

        \underline{Step.1}\; We move the coulomb gauge condition $d^\ast=0$ to stably $c$-homotopic maps. In our context, 
        \[
            d^\ast \colon i\Omega^1(X)^{-\iota^\ast} \to i\Omega^0(X)^{-\iota^*} 
        \]
        is surjective since we assume that $H^0(X, \Z)^{-\iota^*}=0$. Therefore the following lemma follows from \cref{stable_SWC}. 
        \begin{lemma}\label{coulomb_gauge}
            The SWC-pair $(l_X, c_X)$ is stably c-homotopic to the SWC-pair $(l_X \oplus d^*, c_X \oplus 0)$. 
        \end{lemma}
        
        \underline{Step.2}\; In this step, we prove the following lemma which is a counter part of the \cite[Lemma 6.16]{MR4602726}. 
        \begin{lemma}
            Let \[V^{k+1/2-m, I}_{Y}:=L^2_{k+1/2-m}(i\Omega^0(Y)^{-\iota^*} \oplus i\Omega^1(Y)^{-\iota^*} \oplus \Gamma(S_Y)^I)\] 
            and we denote an operator 
            \[
                V_{X_0}^I \times V_{X_1}^I \to V^{k+1/2-m, I}_{Y};\; (x_0, x_1) \mapsto r(\partial_{\vec{n}}^m(x_0))-r(\partial_{\vec{n}}^m(x_1))
            \]
            by $D^{(m)}$, where $\vec{n}$ is an outward normal vector of $\partial X_0=-\partial X_1$.             
            Then we have 
            \begin{align}
               (proj,D^{(m)}, \dots, D^{(0)})\colon V_{X_0}^I \times V_{X_1}^I \to V_X^I \oplus \bigoplus_{m=0}^{k}V^{k+1/2-m, I}_{Y}
            \end{align}
            is isomorphism. 
            Here $proj$ is given as follows: Let $p$ be the orthogonal projection to the kernel of $(D^{(m)}, \dots, D^{(0)})$ in $L^2_{k+1}$ inner product and let $f$ be a map from the kernel of $(D^{(m)}, \dots, D^{(0)})$ to $V_X^I$ given by the gluing theorem of the Sobolev space. 
            In particular, $((proj \circ ((l_{X_0}\oplus d^*) \times (l_{X_1}\oplus d^*)), D^{(m)}, \dots, D^{(0)}), (proj \circ (c_{X_0} \times c_{X_1}), 0,\dots,0,0))$ is stably c-homotopic to $(l_X \oplus d^*, c_X \oplus 0)$. 
        \end{lemma}
        \begin{proof}
            Just using the Sobolev gluing theorem \cite[Section 4.6, Lemma 3.]{MR2312050}. 
        \end{proof}
        
        \underline{Step.3}\; We denote the operator $((l_{X_0}\oplus d^*) \times (l_{X_1}\oplus d^*))$ by $l_{X_0, X_1}$. 
 In this step, we deform the operator $D^{(i)}$ in order to change the target of the Real Seiberg--Witten map and we can prove a counter part result of \cite[Lemma 6.18]{MR4602726} in our context in the same way. 
        \begin{lemma}
            Let $E^{(m)}$ be an operator that is given by 
            \[
                U_{X_0}^I \times U_{X_1}^I \to V^{k+1/2-m, I}_{Y};\;(y_0, y_1) \mapsto r(\partial_{\vec{n}}^m(y_0))-r(\partial_{\vec{n}}^m(y_1)). 
            \]
            For $m=1, \dots, k$ and $t \in [0, 1]$, let $D^{(m)}_t$ be an operator given by 
            \[
                D^{(m)}_t=(1-t)D^{(m)} + t E^{(m-1)} \circ l_{X_0, X_1}
            \]
            and let $l_t$ and $c_t$ be maps given by
            \begin{align*}
                l_t&:=(proj \circ l_{X_0, X_1}, D^{(k)}_t, \dots, D^{(1)}_t, D^{(0)}) , \\
            c_t&:=(proj \circ (c_{X_0} \times c_{X_1}), tE^{(k-1)}\circ (c_{X_0} \times c_{X_1}), t E^{(0)} \circ (c_{X_0} \times c_{X_1}), 0)\colon \\
            &V_{X_0}^I \times V_{X_0}^I \to U_{X}^I \oplus \bigoplus_{m=0}^{k}V^{k+1/2-m, I}_{Y}. 
        \end{align*}
            Then $l_t$ is Fredholm. Moreover, the zero set $\{ x 
             \in V_{X_0}^I \times V_{X_0}^I  \mid (l_t + c_t)(x)=0\}$ is independent of $t$. 
        \end{lemma}
        \begin{proof}
            The proof of this lemma is same way that of \cite[Lemma 6.18]{MR4602726}. The key point is that $E^{(m-1)} \circ l_{X_0, X_1}-D^{(m)}$ contains at most $m-1$ derivative in the normal vector of $Y$. 
        \end{proof}
        Then we have counter part of the \cite[Proposition 6.17]{MR4602726}.
        \begin{proposition}
            The SWC-pair 
            \begin{align*}
                (&(proj \circ l_{X_0, X_1}, E^{(k-1)} \circ l_{X_0, X_1}, \dots, E^{(0)} \circ l_{X_0, X_1}, D^{(0)}), \\ &(proj \circ (c_{X_0} \times c_{X_1}), E^{(k-1)}\circ (c_{X_0} \times c_{X_1}), E^{(0)} \circ (c_{X_0} \times c_{X_1}), 0))
            \end{align*}
            is c-homotopic to the SWC-pair $((proj \circ l_{X_0, X_1}, D^{(k)}, \dots, D^{(0)}), (proj \circ (c_{X_0} \times c_{X_1}), 0, \dots,0,0))$. 
        \end{proposition}
        
        \underline{Step.4}\; In this step, we use gluing theorem of Sobolev space again. In our context, the counter part of \cite[Lemma 6.19]{MR4602726} is different from original one. There are no constant function that is preserved by $-\iota^*$ on $X$ since we assume that $X$ is connected. 
        Here is the counter part of \cite[Lemma 6.19]{MR4602726}. 
\begin{lemma}
    We have 
    \begin{align*}
        (proj,E^{(k-1)}, \dots, E^{(0)}) \colon U_{X_0}^I \times U_{X_1}^I \to U_{X}^I \oplus \bigoplus_{m=0}^{k}V^{k+1/2-m, I}_{Y}
    \end{align*}
    is an isomorphism. 
\end{lemma}
\begin{proof}
    This is consequence of the gluing theorem of the Sobolev space \cite[Section 4.6, Lemma 3.]{MR2312050}. Note that there is no $\R$-factor that appears in the target of the map in \cite[Lemma 6.19]{MR4602726} because there is no nontrivial constant function in $i\Omega^0(X)^{-\iota^*}$. 
\end{proof}
Therefore we have the counter part of \cite[Lemma 6.20]{MR4602726} immediately. 
\begin{lemma}\label{step4}
    The SWC-pair $((l_{X_0, X_1}, D^{(0)}), (c_{X_0} \times c_{X_1}, 0))$ is c-homotopic to the SWC-pair
    \begin{align*}
        (&(proj \circ l_{X_0, X_1}, E^{(k-1)} \circ l_{X_0, X_1}, \dots, E^{(0)} \circ l_{X_0, X_1}, D^{(0)}), \\ &(proj \circ (c_{X_0} \times c_{X_1}), E^{(k-1)}\circ (c_{X_0} \times c_{X_1}), E^{(0)} \circ (c_{X_0} \times c_{X_1}), 0)). 
    \end{align*}
\end{lemma}
    Let us explain the difference between \cref{step4} and \cite[Lemma 6.20]{MR4602726}. In \cite[Lemma 6.20]{MR4602726}, the $D^{(0)}$-factor in our case is replaced by $(D_{Y}, D_H)$. These operators are given by a decomposition of $D^{(0)}$. We explain the corresponding operators in our context. 

    Let us denote \[\Omega^0_0(Y)^{-\iota^*}=\{ f \in i\Omega^0(Y)^{-\iota^*} \mid \int_{Y_j}f=0\; \text{for all connected component}Y_j \text{of}\;Y \}.\]
    If $Y$ is connected this coincides with $i\Omega^0(Y)^{-\iota^*}$. 
    We have a decomposition  
    \[
        V^{k+1/2, I}_Y = L^2_{k+1/2}(Coul(Y)^I) \oplus L^2_{k+1/2}(i(d\Omega^0(Y))^{-\iota^*} \oplus \Omega^0_0(Y)^{-\iota^*}) \oplus \R^{b_0(Y)-1}. 
    \] 
    
    Let us define $(D_Y, D_H, D_{\R})$ be the each corresponding factor of $D^{(0)}$. Here is a difference between the Real case and the ordinary case. In \cite[Lemma 6.20]{MR4602726}, $D_{\R}$ component disappears because that factor absorbed by the $\R$-factor that appears in the target of the map in \cite[Lemma 6.19]{MR4602726}. 
    After we deform $D_H$ by a homotopy the kernel of $(d^*_{X_0}\times d^*_{X_1}, D_H)$ coincides with the double coulomb condition since $Y_2$, which corresponds to $Y$ in our context, is connected and $d^*a=0$ implies $\boldsymbol{t} \ast a=0$. 
    
    However, in our situation, $Y$ can be disconnected. If the connected component of $Y$ is $2$, $D_{\R}$ does not disappears in our setting. 
        
        \underline{Step.5}\; In this step, we deform $D_H$ component that corresponds to the boundary conditions for gauge fixing. We check that the kernel of $(d^*_{X_0}\times d^*_{X_1}, D_H, D_{\R})$ corresponds the double coulomb condition after we deform $D_H$. 
        For $x_i=(\hat{a}_i, \hat{\phi}_i) \in V^I_{X_i}$, we have a decomposition $\boldsymbol{t}\ast \hat{a}_i=a_i + b_i$, where $a_i \in \ker d^* \subset i\Omega^1(Y)^{-\iota^*}$ and $b_i \in \text{im} d \subset i\Omega^1(Y)^{-\iota^*}$. Then we have $\hat{a}_i|_{Y}=a_i + b_i +e_i dt+e'_i dt$, where $e_i \in i\Omega^0_0(Y)^{-\iota^*}$ and $e'_i$ is a locally constant function that is fixed by $-\iota^*$.    
        Let $D_{H, t}$ be the operator
        \[
            D_{H, t}\colon V_{X_0}^I \times V_{X_1}^I \to H^I ;\; (\hat{a}_0, \hat{a}_1)=(b_0-b_1, t\bar{d}^*(b_0+b_1)+(1-t)\cdot(e_0-e_1)) \]
       , where $\bar{d}^*=d^*(d^* d)^{-1/2}$.  
        To prove $((l_{X_0, X_1}, (D_Y, D_{H,t}, D_{\R})), (c_{X_0} \times c_{X_1}, 0))$ be a continuous family of SWC-pair, we give a family of gauge transformation $u_i^t$ that satisfies the following significance: if $b_0=b_1$, we have $D_H(\hat{a}_0, \hat{a}_1)=0$ if and only if $D_{H, t}(u_0^t \cdot \hat{a}_0, u_1^t \cdot \hat{a}_1)=0$. Moreover, we can take $u_i^t$ from the identity component of the gauge group. The following lemma is a counter part of \cite[Lemma 6.21]{MR4602726}. There are two cases that corresponds that $Y$ is connected or not. 
        
    Firstly we state the case that $Y$ is connected.   We can prove this lemma in the same way as \cite[Lemma 6.21]{MR4602726}. Moreover, in this case, the proof is easier than that of \cite[Lemma 6.21]{MR4602726} since there are no nontrivial constant functions. To compare with the case that $Y$ is disconnected, we describe the proof in this case carefully. 

        \begin{lemma}[The case that $Y$ is connected. ]\label{boundary_connected}
            Let $W^I \subset L^2_{k+2}(X_0, \R)^{-\iota^*} \times L^2_{k+2}(X_1, \R)^{-\iota^*}$ be the subspace of functions $(f_0, f_1)$ that satisfies the following conditions:
            \begin{enumerate}
                \item $d^* d f_i=\Delta f_i=0$ for $i=0, 1$;
                \item $f_0|_Y=f_1|_Y$.  
            \end{enumerate}
            Then the map $\rho_t \colon W^I \to L^2_{k+1/2}(\Omega^0(Y))^{-\iota^*}$ defined by 
            \[
                \rho_t(f_0, f_1)=2t\bar{d^*} d(f_0|_Y)+(1-t)(\partial_{\vec{n}}f_0|_Y - \partial_{\vec{n}}f_1|_Y)
            \]
            is an isomorphism. Note that $\Omega^0(Y)^{-\iota^*}=\Omega^0_0(Y)^{-\iota^*}$ in this case. 
        \end{lemma}
        \begin{proof}
            The key fact is the following:
            Let $W$ be a compact connected Riemannian manifold with boundary. Let $s \ge 2$ be a real number. Then the map
            \begin{align}\label{key1}
                L^2_{s}(W) \to L^2_{s-2}(W) \oplus L^2_{s-1/2}(\partial W);\; f \mapsto (\Delta f, f|_{\partial W})
            \end{align}
            is an isomorphism. Moreover, the map 
            \begin{align}\label{key2}
                L^2_{s}(W) \to L^2_{s-2}(W) \oplus L^2_{s-3/2}(\partial W);\; f \mapsto (\Delta f, (\partial_{\vec{n}}f)|_{\partial W})
            \end{align}
            is a Fredholm map. The kernel is constant functions on $W$ and the image is the set of the pairs $(g, h) \in L^2_{s-2}(W) \oplus L^2_{s-3/2}(\partial W)$ with $\int_W g=\int_{\partial W} h$. In particular, the kernel is one dimensional and the cokernel is isomorphic to $\R^{b_0(\partial W)}$. If the boundary $\partial W$ is connected and there is a non-trivial involution $\iota$ on $W$, the map 
            \begin{align}\label{d*d}
                \bar{d}^* d \colon L^2_{s-1/2}(\partial W)^{-\iota^*} \to L^2_{s-3/2}(\partial W)^{-\iota^*}
            \end{align}
            is isomorphism and the two maps 
            \begin{align}\label{key1*}
                L^2_{s}(W)^{-\iota^*} \to L^2_{s-2}(W)^{-\iota^*} \oplus L^2_{s-3/2}(\partial W)^{-\iota^*};\; f \mapsto (\Delta f, (\bar{d}^* df)|_{\partial W})
            \end{align}
            and 
            \begin{align}\label{key2*}
                L^2_{s}(W)^{-\iota^*} \to L^2_{s-2}(W)^{-\iota^*} \oplus L^2_{s-3/2}(\partial W)^{-\iota^*};\; f \mapsto (\Delta f, (\partial_{\vec{n}}f)|_{\partial W})
            \end{align}
            are isomorphisms. Let $\mathcal{H}(X_i) \subset L^2_{k+2}(X_i)$ be the space of harmonic functions on $X_i$ for $i=0, 1$. 
            Then the following maps are isomorphisms:
            \begin{align*}
                \varphi_{1, i} \colon &\mathcal{H}(X_i)^{-\iota^*} \to L^2_{k+1/2}(Y)^{-\iota^*} ;\; f \mapsto (\bar{d}^* df)|_{\partial W},\\
                \varphi_{2, i} \colon &\mathcal{H}(X_i)^{-\iota^*} \to L^2_{k+1/2}(Y)^{-\iota^*} ;\; f \mapsto (\partial_{\vec{n}}f)|_{\partial W},\\
                \varphi_3 \colon  \mathcal{H}(X_0)^{-\iota^*} \times \mathcal{H}(X_1)^{-\iota^*}& \to L^2_{k+1/2}(Y)^{-\iota^*} \times L^2_{k+1/2}(Y)^{-\iota^*};\; (f_0, f_1) \mapsto (f_0|_Y, f_1|_Y),\\
                \varphi_4^t \colon  L^2_{k+1/2}(Y)^{-\iota^*} \times L^2_{k+1/2}(Y)^{-\iota^*}& \to L^2_{k+1/2}(Y)^{-\iota^*} \times L^2_{k-1/2}(Y)^{-\iota^*};\\ & (g_0, g_1) \mapsto (g_0-g_1, A_0^t(g_0)+A_1^t(g_1) ),
            \end{align*}

           , where $A_i^t(g_i)=t \bar{d}^*d g_i +(-1)^i(1-t)\varphi_{2, i} \circ \varphi_{1, i}^{-1}(\bar{d}^*d g_i)$. For $\varphi_4^t$, 
            from the divergence formula, one can easily check that $A_i^t$ are bounded positive self-adjoint operator and $\varphi_4^t$ is an isomorphism.  The composition of the map 
            \begin{align*}
                \mathcal{H}(X_0)^{-\iota^*} \times \mathcal{H}(X_1)^{-\iota^*} \to L^2_{k+1/2}(Y) \times L^2_{k+1/2}(Y) ;\; (f_0, f_1) \mapsto (f_0|_Y, f_1|_Y) 
            \end{align*}
            and $\varphi_4^t$ is an isomorphism and we have that $\rho_t$ is isomorphic. 
        \end{proof}

    Next, we give the counter part of \cite[Lemma 6.21]{MR4602726} when $Y$ is disconnected. Note that from the assumption that $H^1(X_1, \Z)^{-\iota^*}=0$ and $\#\pi_0(Y/\iota)=\#\pi_0(X_1/\iota)=1$, we have the number of connected components of $Y$ and $X_1$ are both $2$.  
        \begin{lemma}[The case that $Y$ is disconnected.]\label{boundary_disconn}
            We assume that $Y$ is disconnected. Let us denote $Y=Y_0 \sqcup Y_1$ and $X_1=X_{1,0} \sqcup X_{1,1}$ by its connected components. Note that $\iota(Y_0)=Y_1$ and $\iota(X_0)=X_1$. 
            Let $W^I \subset L^2_{k+2}(X_0, \R)^{-\iota^*} \times L^2_{k+2}(X_1, \R)^{-\iota^*}$ be the subspace of functions $(f_0, f_1)$ that satisfies the following conditions:
            \begin{enumerate}
                \item $d^* d f_i=\Delta f_i=0$ for $i=0, 1$;
                \item $f_0|_Y=f_1|_Y$;
                \item \label{orthogonal}Let $u_0 \in \mathcal{H}(X_0)^{-\iota^*}$ be the unique function that satisfies $u_0|_Y=(1, -1) \in L^2(Y_0) \oplus L^2(Y_1)$. We assume $\int_{X_0}\langle du_0, df_0 \rangle =0$. 
            \end{enumerate}
        Then the map $\rho_t \colon W^I \to L^2_{k+1/2}(\Omega^0_0(Y))^{-\iota^*}$ defined by 
            \[
                \rho_t(f_0, f_1)=2t\bar{d^*} d(f_0|_Y)
                +(1-t)(\partial_{\vec{n}}f_0|_Y - \partial_{\vec{n}}f_1|_Y)
            \]
            is well-defined and is an isomorphism. 
        \end{lemma}
        \begin{proof}
            Firstly, we check that the image of $\rho_t$ is in $L^2_{k+1/2}(\Omega^0_0(Y)) \subset L^2_{k+1/2}(\Omega^0(Y))$. 
    We only have to check that $\partial_{\vec{n}}f_i \in L^2_{k+1/2}(\Omega^0_0(Y))$. For $f_1$, both $X_{1, 0}$ and $X_{1, 1}$ are connected so that we have $\partial_{\vec{n}}f_1 \in L^2_{k+1/2}(\Omega^0_0(Y_i))$ for $i=0, 1$ from the divergence formula. 
    For $f_0$, we have to prove that 
    \[                                                  
\int_{Y_0}\partial_{\vec{n}}f_0|_{Y_0}=\int_{Y_1}\partial_{\vec{n}}f_0|_{Y_1}=0. 
    \] 
From the condition that $\iota^*f_0=-f_0$, we have
    \[
        \int_{Y_0}\partial_{\vec{n}}f_0|_{Y_0}+\int_{Y_1}\partial_{\vec{n}}f_0|_{Y_1}=0. 
    \]
Let us show that $\int_{Y_0}\partial_{\vec{n}}f_0|_{Y_0}-\int_{Y_1}\partial_{\vec{n}}f_0|_{Y_1}=0$. This is a $L^2$ inner product with $\partial_{\vec{n}}f_0|_{Y}$ and $(-1, 1) \in L^2(Y_0) \oplus L^2(Y_1)$. Therefore, from the assumption i\cref{orthogonal} in the statement, we have 
\begin{align*}
    \int_{Y_0}\partial_{\vec{n}}f_0|_{Y_0}-\int_{Y_1}\partial_{\vec{n}}f_0|_{Y_1}&=\int_{Y} \langle u_0|_{Y}, \partial_{\vec{n}}f_0|_Y \rangle \\
    &=-\int_{X_0}u_0 \Delta f_0+\int_{X_0} \langle du_0, df_0 \rangle \\
    &=\int_{X_0} \langle du_0, df_0 \rangle \\
    &=0. 
\end{align*}
Thus we have $\rho_t$ is well-defined. 

            Next we prove that $\rho_t$ is an isomorphism.  In this situation, the map $\varphi_{j, 1}$ in the proof of \cref{boundary_connected} have one dimensional kernel and one dimensional cokernel for each $j=1, 2$. The kernel of these maps are the set of locally constant functions on $X_1$ that preserved by $-\iota^*$. The cokernel of these maps are the set of the  locally constant functions on $Y$ preserved by $-\iota^*$. 
            Let $\mathcal{H}_0(X_1)^{-\iota^*} \subset L^2_{k+2}(X_1)^{\iota^*}$ be the set of harmonic functions on $X_1$ that satisfies $\int_{X_{1,0}}f_1=\int_{X_{1,1}}f_1=0$ and $\iota^*f_1=-f_1$. Then the maps
            \begin{align*}
                \varphi'_{1, 1} \colon &\mathcal{H}_0(X_1)^{-\iota^*} \to L^2_{k+1/2}(\Omega^0_0(Y))^{-\iota^*} ;\; f \mapsto (\bar{d}^* df)|_{Y},\\
                \varphi'_{2, 1} \colon &\mathcal{H}_0(X_1)^{-\iota^*} \to L^2_{k+1/2}(\Omega^0_0(Y))^{-\iota^*} ;\; f \mapsto (\partial_{\vec{n}}f)|_{Y}
            \end{align*}
            are isomorphism. 

            The map $\varphi_{1, 0}$ has one dimensional kernel and one dimensional cokernel in this situation. The kernel is generated by $u_0$. 
            The cokernel of $\varphi_{1, 0}$ is the set of locally constant functions on $Y$ that is preserved by $-\iota^*$. 
            
            The map $\varphi_{2, 0}$ is isomorphism and if $f_0$ satisfies $\int_{X_0}\langle du_0, df_0 \rangle =0$, then we have $\varphi_{2, 0}(f_0) \in L^2_{k+1/2}(\Omega^0_0(Y))$. 

            Therefore, let $\mathcal{H}_0(X_0)^{-\iota^*} \subset L^2_{k+2}(X_0)$ be a set of harmonic functions on $X_0$ that satisfies $\int_{X_0}\langle du_0, df_0 \rangle =0$ and $\iota^*f_0=-f_0$. Then the following maps are isomorphism:
            \begin{align*}
                \varphi'_{1, 0} \colon &\mathcal{H}_0(X_i)^{-\iota^*} \to L^2_{k+1/2}(\Omega^0_0(Y))^{-\iota^*} ;\; f \mapsto (\bar{d}^* df)|_{Y},\\
                \varphi'_{2, 0} \colon &\mathcal{H}_0(X_i)^{-\iota^*} \to L^2_{k+1/2}(\Omega^0_0(Y))^{-\iota^*} ;\; f \mapsto (\partial_{\vec{n}}f)|_{Y}. 
            \end{align*}
        Thus we have that 
        \begin{align*}
            \varphi_3' \colon  \mathcal{H}_0(X_0)^{-\iota^*} \times \mathcal{H}_0(X_1)^{-\iota^*}& \to L^2_{k+1/2}(\Omega^0_0(Y))^{-\iota^*} \times L^2_{k+1/2}(\Omega^0_0(Y))^{-\iota^*}; \\ (f_0, f_1) \mapsto &(p(f_0|_Y), p(f_1|_Y)) \end{align*}
            and 
            \begin{align*}
            {\varphi'}_4^t \colon  L^2_{k+1/2}(\Omega^0_0(Y))^{-\iota^*} \times &L^2_{k+1/2}(\Omega^0_0(Y))^{-\iota^*} \to L^2_{k+1/2}(\Omega^0_0(Y))^{-\iota^*} \times L^2_{k-1/2}(\Omega^0_0(Y))^{-\iota^*};\\ & (g_0, g_1) \mapsto (g_0-g_1, A_0^t(g_0)+A_1^t(g_1) )
        \end{align*}
        are isomorphism, where $p \colon \Omega^0(Y) \to \Omega^0_0(Y)$ is the $L^2$ orthogonal projection. Therefore, if we set $W^I_0=\{(f_0, f_1) \in \mathcal{H}_0(X_0)^{-\iota^*} \times \mathcal{H}_0(X_1)^{-\iota^*} \mid p(f_0|_Y-f_1|_Y)=0\}$, we can prove the map
        \[
            \rho'_t \colon W^I_0 \to L^2_{k+1/2}(\Omega^0_0(Y))^{-\iota^*};\; (f_0, f_1) \mapsto 2t\bar{d^*} d(f_0|_Y)
            +(1-t)(\partial_{\vec{n}}f_0|_Y - \partial_{\vec{n}}f_1|_Y)
        \]
        be an isomorphism in the same way in the proof of \cref{boundary_connected}. 
        For $(f_0, f_1) \in W^I_0$ let $c=f_0|_Y-f_1|_Y-p(f_0|_Y-f_1|_Y)$ and let $c' \in L^2(X_1)^{-\iota^*}$ be a locally constant function with $c'|_Y=c$. 
        Then $(f_0, f_1-c') \in W^I$ and $\rho_t(f_0, f_1-c')=\rho'_t(f_0, f_1)$. 
        The isomorphism between $W_0^I$ and $W^I$ is given as above. 
        Therefore we finished the proof.  
        \end{proof}

        Then we have the counter part of \cite[Proposition 6.22]{MR4602726}. 
        \begin{proposition}
            For any $t \in [0, 1]$, the pair $((l_{X_0, X_1},D_Y, D_{H, t} D_{\R}), (c_{X_0} \times c_{X_1}, 0, 0, 0))$ is a SWC-pair. 
            Therefore, this provides the c-homotopy between $t=0$ and $t=1$. 
            Note that \[((l_{X_o, X_1}, D^{(0)}), (c_{X_0} \times c_{X_1}, 0))=((l_{X_0, X_1},D_Y, D_{H, 0} D_{\R}), (c_{X_0} \times c_{X_1}, 0, 0, 0)).\]
            Moreover, the kernel of the operator \begin{align*}&d^*_{X_0, X_1} \oplus D_{H, 1} \oplus D_{\R} \colon L^2_{k+2}(\Omega^1(X_0) \times \Omega^1(X_1))^{-\iota^*} \to \\ & (L^2_{k+1}(\Omega^0(X_0) \oplus \Omega^0(X_0)))^{-\iota^*}  \oplus L^2_{k+1/2}(d\Omega^0(Y) \oplus \Omega^0_0(Y))^{-\iota^*} \oplus \R \end{align*} coincides with $\Omega^1_{CC}(X_0)^{-\iota^*} \times \Omega^1_{CC}(X_1)^{-\iota^*}$. 
        \end{proposition}
        \begin{proof}
            From \cref{boundary_connected} and \cref{boundary_disconn}, we have that there are $1$-parameter family of the gauge transformations $u_i^t$ that satisfies that if $b_0=b_1$, we have $D_H(\hat{a}_0, \hat{a}_1)=0$ if and only if $D_{H, t}(u_0^t \cdot \hat{a}_0, u_1^t \cdot \hat{a}_1)=0$. Then we have that $((l_{X_0, X_1},D_Y, D_{H, t} D_{\R}), (c_{X_0} \times c_{X_1}, 0, 0, 0))$ is a SWC-pair. 

            To prove latter half of the proposition, let $(\hat{a}_0, \hat{a}_1) \in \Omega^1(X_0)^{-\iota^*} \times \Omega^1(X_1)^{-\iota^*}$ be a kernel of the operator $d^*_{X_0, X_1} \oplus D_{H, 1} \oplus D_{\R}$. Let us decompose $\hat{a}_i|_{Y}=a_i + b_i +e_i dt+e'_i dt$, where $a_i \in \ker d^* \subset i\Omega^1(Y)^{-\iota^*}$, $b_i \in \text{im} d \subset i\Omega^1(Y)^{-\iota^*}$, $e_i \in i\Omega^0_0(Y)^{-\iota^*}$, and $e'_i$ is a locally constant function that is fixed by $-\iota^*$ in the collar neighborhood of $Y$.  Then we have 
            \[
                D_{H, 1} \oplus D_{\R}(\hat{a}_0, \hat{a}_1) = (b_0-b_1, \bar{d}^*(b_0+b_1), e_0'-e_1')=0. 
            \]
            If we have $b_0=b_1$, from the second term $\bar{d}^*(b_0+b_1)=2\bar{d}^* b_0=2\bar{d}^* b_0=0$. We assume that $H^1(Y)^{-\iota^*}=0$, we have $b_0=b_1=0$ since they are harmonic $1$-form preserved by $-\iota^*$. From the condition of $d^* \hat{a}_i=0$, we have $\int_{Y}e'_i=0$. Since the number of connected components of $X_1$ and that of $\partial {X_1}=-Y$ are the same, we have $e'_1=0$. Therefore $e_0'=0$. Then we proved this proposition.  
        \end{proof}
        
        \underline{Step.6}\; In this step we deform the action of the harmonic gauge transformations. However in our setting only non-trivial harmonic gauge transformation is $-1$ for any boundary conditions which appears in this proof. Therefore we can skip this step. 
        \underline{Step.7}\; In this step we discuss about the choice of the sequence of the finite dimensional subspace of the Hilbert space of the domain of the Real relative Seiberg--Witten map. 
        We can prove a counter part of \cite[Lemma 6.26]{MR4602726} in the same way.  
        \begin{lemma}
            Let us define 
            \[
                W^{n, t, I}_{X_0, X_1}\colon = \{ (x_0, x_1) \in Coul(X_0)^I  \times Coul(X_1)^I \mid  p^{\infty}_{\mu_n}r(x_0)=tp^{\infty}_{\mu_n}r(X_1), \; p^{-\mu_n}_{-\infty}r(x_1)=tp^{-\mu_n}_{-\infty}r(X_0) \}. 
            \]
            Therefore, for any $R >0$, there exist $N$ and $\epsilon$ with the following significance: For all $n>N$, $t \in [0,1]$, and $(x_0, x_1) \in B^+(W^{n, t, I}_{X_0, X_1}, R)$ that satisfies
            \begin{itemize}
                \item $\lVert p^{\mu_n}_{-\mu_n}(x_0-x_1) \rVert _{L^2_{k+1/2}(Y)} \le \epsilon$, 
                \item $\lVert p_{U_i}\circ SW_{X_i}(x_i) \rVert \le \epsilon$,
            \end{itemize}
            we have $\lVert x_i \rVert _{L^2_{k+1}(X_i)} \le R_i+1$, where $R_i$ for $i=0, 1$ are the constants that appears in \cref{representative_e}. 
        \end{lemma}
Therefore \cref{gluing_general} follows. 
\end{proof}

\section{Definition of the invariant}\label{def_inv_surfaces}
\subsection{Definition of the invariant of the embedded surfaces}
In this section, we will define invariants for surfaces in a $4$-manifold using the invariants of a Real structure discussed in the previous section. Furthermore, in the case of a $2$-knot, we define a more manageable numerical invariant, and provide an alternative definition that does not rely on finite-dimensional approximations of that invariant.

\begin{definition}\label{nice_emb_closed}
    Let $X$ be a closed oriented $4$-manifold and $S\neq \emptyset$ be a connected surface embedded in $X$. If $(X, S)$ satisfies $H_1(X, \Z/2)=0$ and $0=[S]_2 \in H_2(X, \Z/2)$, we call $(X, S)$ is a nice pair. 
\end{definition}

From \cite[Corollary 2.10, Remark 2.11]{nagami2000existence}, if $(X, S)$ is a nice pair, there is a unique double branched cover of $X$ along $S$. We denote $\Sigma(X, S)$ by the double branched cover. 

\begin{lemma}\label{unique_Real_involution}
    Let $(X, S)$ be a nice pair. Let $R \Spinc(\iota)$ be a set of isomorphism classes of Real $\mathrm{spin}^c$ structures on $\Sigma(X, S)$ that covers the covering involution $\iota$. Then 
    \[
 \begin{array}{r@{\,\,}c@{\,\,}c@{\,\,}c}
 \varphi\colon&R \Spinc(\iota)&\longrightarrow&\{ c \in H^2(\Sigma(X, S), \Z)^{-\iota^*}\mid c = w_2(\Sigma(X, S)) \mod2 \}\\
 &\rotatebox{90}{$\in$}&&\rotatebox{90}{$\in$}\\
 &[(\mathfrak{s}, I)]&\longmapsto&c_1(\mathfrak{s})
 \end{array}
 \]
 is a bijection. 
\end{lemma}
\begin{proof}
    From \cite[Lemma 4.5]{konno2023involutions}, there is no $2$-torsion in $H^2(\Sigma(X, S), \Z)$ and $b_1(\Sigma(X, S))=0$. 
    Since $H^1(\Sigma(X, S), \Z)=0$ and \cref{classification_of_real_spinc}, the number of isomorphism class of Real structures on a $\mathrm{spin}^c$ structure $\mathfrak{s}$ on $\Sigma(X, S)$ is at most $1$. 
    Since there are no $2$-torsion in $H^2(\Sigma(X, S), \Z)$, isomorphism class of $\mathrm{spin}^c$ structures on $\Sigma(X, S)$ corresponds to characteristic elements in $H^2(\Sigma(X, S))$. We denote $c_1(\mathfrak{s})$ by the first Chern class of the determinant line bundle of the $\mathrm{spin}^c$ structure $\mathfrak{s}$. 
    If $(\mathfrak{s}, I)$ be a Real $\mathrm{spin}^c$ structure, we have $\iota^*c_1(\mathfrak{s})=-c_1(\mathfrak{s})$. Thus the map $\varphi$ is well defined and injective. 

    To prove $\varphi$ be a surjection, we have to construct a Real structure $I$ on a $\mathrm{spin}^c$ structure with $\iota^*c_1(\mathfrak{s})=-c_1(\mathfrak{s})$. This is shown in \cite[Section.2]{konno2023involutions}. 
\end{proof}

\begin{definition}\label{invariant_closed}
    Let $(X, S)$ be a nice pair and $c \in H^2(\Sigma(X, S), \Z)^{-\iota^*}$. We define a $\Z/2$-equivariant stable cohomotopy invariant of the triple $(X, S, c)$ by 
    \[
        \Psi(X, S, c) := \Psi(\Sigma(X, S), \mathfrak{s}, I)
    \]
   , where $\mathfrak{s}$ is the $\mathrm{spin}^c$ structure on $\Sigma(X, S)$ with $c_1(\mathfrak{s})=c$ and $I$ is its unique Real structure. 
\end{definition}

\subsection{Definition of the numerical invariant of $2$-knots and $P^2$-knots}

\subsubsection{Invariant of $2$-knots}

Let us define the numerical invariant of $2$-knots in the $4$-sphere. 

Let $S \subset S^4$ be a $2$-knot. 
It is clear that $(S^4, S)$ is a nice pair. 
From \cite[Lemma 4.5.]{konno2023involutions}, we see that $b_2(\Sigma(S^4, S))=0$ and $b_1(\Sigma(S^4, S))=0$. 
Therefore $\ind_{\R}(\Slash{D}^I)=0$ and $H^+(\Sigma(S^4, S))^{-\iota^*}=0$. 
Thus the $\Z/2$-equivariant stable cohomotopy invariant $\Psi(S^4, S, c)$ of $(S^4, S)$ lives in
\[  
    \text{colim}_{U\in \mathcal{U}}[(V \oplus U)^+ \to U^+]
\] 
for all Real $\mathrm{spin}^c$ structure $c$, where $V=V_0-V_1$ is a virtual vector space with $V_0, V_1 \subset \mathcal{U}$ and $\dim V=\dim V_0-\dim V_1=0$. Note that we have no canonical choice to trivialize this virtual vector space. In particular, there is no canonical orientation of $V$. 
\begin{definition}\label{def_degree}
    Let $S$ be a $2$-knot. We define $\lvert deg(S, c) \rvert \in \mathbb{N}$ by the absolute value of the mapping degree of $\Psi(S^4, S, c)$. Here we forget the group action of $\Z/2$ or $\Z/4$. If the Real $\mathrm{spin}^c$ structure $c$ comes from the spin structure on $\Sigma(S^4, S)$ and the Real involution $I$ that is given by the composition of a lift $\tilde{\iota}$ of the involution $\iota$ to the spinor bundle and $j$, we denote $\lvert deg(S) \rvert$ by $\lvert deg(S, c) \rvert$. 
\end{definition}
\begin{remark}
    Since $(S^4, S)$ is a nice pair, there is no $2$-torsion in $H^2(\Sigma_2(S^4, S), \Z)$ so that $H^1(\Sigma_2(S^4, S), \Z/2) \cong \hom (H_1(\Sigma_2(S^4, S), \Z), \Z/2)=0$. Thus the spin structure on $\Sigma_2(S^4, S)$ is unique. 
\end{remark}

\begin{proposition}\label{multiplicative}
    Let $S_0$ and $S_1$ be $2$-knots. Let $c \in H^2(\Sigma(S^4, S_0 \# S_1))^{-\iota^*}$ be a Real $\mathrm{spin}^c$ structure. Then there is Real structure $c_i \in H^2(\Sigma(S^4, S_i))^{-\iota_i^*}$ for $i=0, 1$ such that $c=c_0+c_1$ and we have $\lvert deg(S_0 \# S_1, c) \rvert =\lvert \deg(S_0, c_0) \rvert \lvert \deg(S_1, c_1) \rvert$. In particular, $\lvert deg(S_0 \# S_1) \rvert = \lvert deg(S_0)\rvert \lvert deg(S_1) \rvert$. 
\end{proposition}
\begin{proof}
    Let $D_i$ be a $4$-disk embedded in $S^4$ such that $D_i \cap S_i$ be a standard disk embedding for $i=0, 1$. Note that $\partial D_i \cap S_i$ is the unknot. We see that $S_0 \# S_1$ is given by $(S_0 \setminus D_0 \cup S_1 \setminus D_1) \subset (S^4 \setminus D_0 \cup_{-\partial D_0=\partial D_1} S^4\setminus D_1)$. From \cref{gluing_conn}, we have that $\Psi(S^4, S_0 \# S_1, c)$ is given by $\Psi(S^4, S_0, c_0) \wedge \Psi(S^4, S_1, c_1)$. The mapping degree is multiplicative and the proof finished. 
\end{proof}

\subsubsection{Invariant of $P^2$-knots}
Next we define the invariant of $P^2$-knots. 
Let $P \subset S^4$ be an embedding of $\R P^2$. 
Firstly, let us recall the notion of the normal Euler number of $P^2$-knots:
\begin{definition}\label{def_of_normal_euler_number}
    Let $P^2$ be a $P^2$-knot and let $l$ be the local coefficient given by the determinant bundle of $T \R P^2$. Then we have the Euler class $e(N)$ of the normal bundle $N$ of $P$ that values $l$. We have a pairing $\langle e(N), [\R P^2] \rangle$, where $[\R P^2] \in H_2(\R P^2, l)$ is the fundamental class. This is the normal Euler number of $P$. 
\end{definition}

We assume that the normal Euler number $P \circ P$ of $P$ is $\pm 2$. 

\begin{lemma}\label{formal_dim_P^2}
    Let $\Sigma_2(S^4, P)$ be a double branched cover along $P$. Let $\iota \colon \Sigma_2(S^4, P) \to \Sigma_2(S^4, P)$ be a covering involution. 
    Then we have following properties:
    \begin{enumerate}
        \item $b_1(\Sigma_2(S^4, P))=0$.
        \item $b_2(\Sigma_2(S^4, P))=1$. 
        \item If $P \circ P=2$ we have $b^+(\Sigma_2(S^4, P))=0$ and if $P \circ P=-2$ we have $b^+(\Sigma_2(S^4, P))=1$. 
        \item The covering involution acts $H^2(\Sigma_2(S^4, P), \Z)/\text{Tor}$ by $-1$. 
        \item There are no $2$-torsion in $H^2(\Sigma_2(S^4, P), \Z)$.
        \item $\Sigma_2(S^4, P)$ is not a spin $4$-manifold.
    \end{enumerate}
\end{lemma}
\begin{proof}
    From (1) to (5) are shown in \cite[Lemma 4.5]{konno2023involutions}. Our situation corresponds to the case that $L=L'=\text{unknot}$ and $S$ is a cobordism between unknot given by $P\setminus D^2\sqcup D^2$ in \cite[Lemma 5.5]{konno2023involutions}. 
    Nagami \cite[Theorem 1.1]{nagami2000existence} tells us that the double branched cover along non-orientable surface is not spin. Thus we have (6). 
\end{proof}
\begin{proposition}\label{invariant_of_P2_knot}
    Let $P$ be a $P^2$-knot in $S^4$ with $P \circ P= \pm 2$. Then we have following properties:
    \begin{itemize}
        \item $(S^4, P)$ be a nice pair in the meaning of \cref{nice_emb_closed}. 
        \item The isomorphism class of Real $\mathrm{spin}^c$ structure on $\Sigma_2(S^4, P)$ such that the formal dimension of the Real Seiberg--Witten equation is $0$ corresponds to $\text{Ext}(H^2(\Sigma_2(S^4, P), \Z)^{-\iota^*}) \cong \text{Tor}(H_1(\Sigma_2(S^4, P),\Z)^{-\iota^*})$. 
    In particular, if $\Sigma_2(S^4, P)$ be an integer homology $\C P^2$ or $\bar{\C P^2}$ we have unique Real $\mathrm{spin}^c$ structure that satisfies formal dimension of the Real Seiberg--Witten equations is $0$. 
    \end{itemize}
\end{proposition}
\begin{proof}
    If $X=S^4$, $H_1(X, \Z/2)=0$ and $0=[P]_2 \in H_2(X, \Z/2)$. Therefore $(S^4, P)$ is a nice pair. 
From \cref{unique_Real_involution}, Real $\mathrm{spin}^c$ structures on $\Sigma_2(S^4, P)$ corresponds to characteristic elements in $H^2(\Sigma_2(S^4, P))$ that satisfies $\iota^*c=-c$. 
From (4) in \cref{formal_dim_P^2}, it coincides to the elements of the set $(\hom(H_2(\Sigma_2(S^4, P), \Z), \Z) \oplus \text{Ext}(H^2(\Sigma_2(S^4, P), \Z), \Z)) \setminus 2\hom(H_2(\Sigma_2(S^4, P), \Z), \Z)$. Therefore a characteristic element is determined by its torsion part and the value $\langle c^2, [\Sigma_2(S^4, P)] \rangle$. 
Let $c$ be a characteristic element. The formal dimension $d(c)$ of the Real Seiberg--Witten moduli space is given by 
\[
    d(c)=\frac{c^2-\sigma(\Sigma_2(S^4, P))}{8}-b^+_{-\iota^*}(\Sigma_2(S^4, P)). 
\]
If $P \circ P=2$, we have $b^+_{-\iota^*}(\Sigma_2(S^4, P))=0$ and $\sigma(\Sigma_2(S^4, P))=-1$. Let $c$ be a characteristic element with $c^2=-1$, we have $d(c)=0$. 

If $P \circ P=2$, we have $b^+_{-\iota^*}(\Sigma_2(S^4, P))=1$ and $\sigma(\Sigma_2(S^4, P))=1$. Let $c$ be a characteristic element with $c^2=9$, we have $d(c)=0$. Thus we prove the second bullet of the proposition. 
\end{proof}
\begin{definition}\label{def_degree_P2}
    Let $P \in S^4$ be a $P^2$-knot with $P \circ P=\pm 2$ and $c$ be a characteristic element such that $d(c)=0$. 
    Then the non-equivariant Real Bauer--Furuta invariant of the Real $\mathrm{spin}^c$ structure that corresponds to $c$ is determined by the degree of $\Psi(S^4, P, c)$. 
    Then we define the degree $\lvert deg(P, c) \rvert $ by the absolute value of the mapping degree of $\Psi(S^4, P, c)$. 
    In particular, if $\Sigma_2(S^4, P)$ is an integer homology $\C P^2$ or $-\C P^2$, we denote $\lvert deg(P, c) \rvert $ by $\lvert deg(P) \rvert$. 
\end{definition} 

Using \cref{gluing_general}, we have following connected sum formula for the degree of $2$-knots and $P^2$-knots. 
\begin{theorem}\label{conn_sum_S_P}
    Let $S \subset S^4$ be a $2$-knot and 
    let $P$ be a $P^2$-knot. 
    We denote a $P^2$-knot $P \# S$ by $P'$. Let $c_P$ and $c_S$ be a Real $\mathrm{spin}^c$ structure on $\Sigma_2(S^4, P)$ and $\Sigma_2(S^4, S)$ respectively. We assume that the formal dimension of the moduli space of the Real Seiberg--Witten equations of these Real $\mathrm{spin}^c$ structure. Then we have the Real $\mathrm{spin}^c$ structure $c'=c_P \# c_S$ on $\Sigma_2(S^4, P')=\Sigma_2(S^4, P) \# \Sigma_2(S^4, S)$ satisfies the formal dimension of the moduli space of the Real Seiberg--Witten equations is $0$ and we have
    \begin{align}
       \lvert deg(P \# S, c') \rvert= \lvert deg(P, c_P) \rvert \lvert deg(S, c_S) \rvert. 
    \end{align}
\end{theorem}

\subsection{Alternative definition of the degree}
To calculate the degrees of $2$-knots and $P^2$-knots, we define an alternative definition of the degree that does not use the finite dimensional approximation. 

Let $X$ be a closed $4$-manifold with involution $\iota$ and $H^1(X, \R)^{-\iota^*}=0$. Suppose $(\mathfrak{s}, I)$ is a Real $\mathrm{spin}^c$ structure which covers $\iota$ and the formal dimension of the Real seiberg--Witten moduli space $d(\mathfrak{s})$ is $0$. 
Let $\mathcal{C}=L^2_{k+1}(i\Omega^1(X) \oplus \Gamma(S^+))$ and let $\mathcal{G}_0=\{ e^{if} \in L^2_{k+2}(\mathcal{G}) \mid f \in L^2_{k+2}(X), \int_X f=0 \}$. 
Note that from the condition $H^1(X, \R)^{-\iota^*}=0$, the Real geuge group $\mathcal{G}^I$ is isomorphic to $\{\pm 1\} \times \mathcal{G}_0^I$. 
We denote $\tilde{\mathcal{B}}^I$ by $\mathcal{C}^I/\mathcal{G}_0^I$. Note that $\tilde{\mathcal{B}}^I$ is contractible. 
    Then we have that there exist $\eta \in i\Omega^+(X)^{-\iota^*}$ such that the framed moduli space $\tilde{\mathcal{M}}_{\eta}^I(X, \mathfrak{s}, I) \subset \tilde{\mathcal{B}}^I$ of the Real Seiberg--Witten equations perturbed by $\eta$ is transverse. The proof of the existence of $\eta$ is the same as the ordinary Seiberg--Witten theory. Note that $\tilde{\mathcal{M}}_{\eta}^I(X, \mathfrak{s}, I)$ is a set of finite points. 
    Since $\tilde{\mathcal{B}}^I$ contractible, the determinant bundle of the Real $\mathrm{spin}^c$ Dirac operator on $\tilde{\mathcal{B}}^I$ is isomorphic to trivial real line bundle. However, we have no canonical trivialization on this line bundle. 

    Therefore, the signed count of the elements of $\tilde{\mathcal{M}}_{\eta}^I(X, \mathfrak{s}, I)$ is well defined up to over all sign. Note that our counting including reducible solutions so that we have no assumptions on $b^+(X)^{-\iota^*}$. 

\begin{definition}\label{alternative_def_deg}
    Let $X$ be a closed $4$-manifold with involution $\iota$ and $H^1(X, \R)^{-\iota^*}=0$. 
    Suppose $(\mathfrak{s}, I)$ is a Real $\mathrm{spin}^c$ structure which covers $\iota$ and the formal dimension of the formal dimension of the Real seiberg--Witten moduli space $d(\mathfrak{s})$ is $0$. 
    We define $\lvert deg'(X, \mathfrak{s}, I) \rvert$ by the absolute value of the signed count of the elements of $\tilde{\mathcal{M}}_{\eta}^I(X, \mathfrak{s}, I)$, where $\eta$ be a generic perturbation such that the framed moduli space is transversal. 
\end{definition}
\begin{proposition}\label{alternative}
    Let $\Sigma$ be a $2$-knot or $P^2$-knot.
    Let $c=(\mathfrak{s}, I)$ be a Real $\mathrm{spin}^c$ structure on $\Sigma_2(S^4, \Sigma)$ that satisfies the formal dimension of the Real Seiberg--Witten moduli is $0$. 
    We have that the degree $\lvert deg(\Sigma, c)\rvert$ coincides with $\lvert deg'(\Sigma_2(S^4, \Sigma), \mathfrak{s}, I) \rvert$. 
\end{proposition}
\begin{proof}
    From the compactness of $\tilde{\mathcal{M}}_{\eta}^I(\Sigma_2(S^4, \Sigma), \mathfrak{s}, I)$, if we take the finite dimensional approximation of the Real Seiberg--Witten equations perturbed by $\eta$, the the $0$ locus of the finite dimensional approximated Real Seiberg--Witten equations convergent to $\tilde{\mathcal{M}}_{\eta}^I(\Sigma_2(S^4, \Sigma), \mathfrak{s}, I)$. If we take the sufficiently approximated representative of $\Psi(S^4, \mathfrak{s}, I)$, we have natutal bijention to the $0$ locus of the approximated solutions and true solutions. Then one can easily check that the signed count of the true solutions and finite dimensional approximated solutions are the same. Thus we finished the proof.  
\end{proof}

Using this alternative description of the degree, let us compute the invariant for the unknot $U \subset S^4$. 
\begin{proposition}\label{invariant_Unknot}
    Let $U$ is the unknotting embedding of $S^2$ to $S^4$. Then we have $\lvert deg(U) \rvert =1$.      
\end{proposition}
\begin{proof}
    We regard $U$ is the embedding $(\R^2)^+ \to (\R^4)^+$. 
    One can easily check that $\Sigma_2(S^4, U)$ is diffeomorphis to $S^4$ and the covering transformation $\iota \colon (\R^4)^+ \to (\R^4)^+$ is induced by  
    \[
        \R^4 \ni (x_1, x_2, x_3, x_4) \mapsto (x_1, x_2, -x_3, -x_4).
    \]
    Therefore, the positive scalar curvature metric on $S^4$ is preserved by the involution. So that all the solutions of the Real Seiberg--Witten equations are reducible. We have $H^+(S^4, \R)^{-\iota^*}=0$ and $H^1(S^4, \R)^{-\iota^*}=0$, then the reducible solution is unique up to contractible gauge transformations and this is transverse. Then we have $\lvert deg(U) \rvert =1$. 
\end{proof}

In \cref{calculations_subsubsection}, we use $3$-dimensional invariant to compute the degree of the $2$-knot. 
Let $M$ be a closed oriented $3$-manifold with an involtuion $\iota$. 
Suppose $H^1(M, \Z)^{-\iota^*}=0$. 
Let $(\mathfrak{t}, I)$ be a Real $\mathrm{spin}^c$ structures on $M$ which covers $\iota$. 
The formal dimension of the moduli space of $3$-dimensional Seiberg--Witten equations is $0$ and one can check there is a $1$-form $\nu$ with $\iota^* \nu=-\nu$ such that all the solutions of the Real Seiberg--Witten equations perturbed by $\nu$ is transverse. 

Let $\mathcal{C}(Y)$ be a $L^2_{k+1}$ completeion of the $3$-dimensioanl configuation space and $\mathcal{G}_0\colon = \{ e^{if} \mid f \in L^2_{k+2}(Y, \R), \int_Y f=0\}$ be a subgroup of the gauge group. 
Then $\tilde{\mathcal{B}}^I=\mathcal{C}^I(Y)/\mathcal{G}^I_0$ be a contractible and the determinant line bundle of the index of the linearized Real Seiberg--Witten equations is isomorphic to the trivial real line bundle however we have no canonical trivialization. 
Then the signed count of the framend moduli space $\tilde{\mathcal{M}}^I_{\nu}(M, \mathfrak{t}, I) \subset \tilde{\mathcal{B}}^I$ is well defined up to over all sign. 

\begin{definition}\label{3d_deg}
    Let $M$ be a closed oriented $3$-manifold with an involtuion $\iota$. Suppose $H^1(M, \Z)^{-\iota^*}=0$. 
    Let $(\mathfrak{t}, I)$ be a Real $\mathrm{spin}^c$ structures on $M$ which covers $\iota$. 
    We define $\lvert deg(M, \mathfrak{t}, I) \rvert$ by the absolute value of the signed count of the framed moduli space $\tilde{\mathcal{M}}^I_{\nu}(M, \mathfrak{t}, I)$. 
\end{definition}

To compute the signed count of $4$- or $3$-dimensional framed moduli space of the Real Seiberg--Witten equations $\tilde{\mathcal{M}}^I_{\nu}$, we prove some propositions.
\begin{proposition}\label{coincides_with_two_quotients}
    Let $M$ be a $3$-or $4$-manifold with involtuion $\iota$ and $(\mathfrak{s}, I)$ be a Real $\mathrm{spin}^c$ structure. Suppose $\mathcal{C}(M)$ be the configuation space of the Seiberg--Witten equations and there is a $I$-equivariant decomposition $\mathcal{G}\cong \mathrm{U}(1) \times \mathcal{G}^{\perp}$. Note that the action $I$ on $\mathcal{G}$ is given by $I\cdot u\colon = I \circ u \circ I^{-1}=\iota^*\bar{u}$. Then the natural map 
    \[
        \mathcal{C}^I/(\mathcal{G}^{\perp})^I \to (\mathcal{C}/\mathcal{G}^{\perp})^I
    \]
    is injective. 
    If the fixed point part of $\iota^*$ on $H^1(M, \Z)$ is trivial i.e. $H^1(M, \Z)^{\iota^*}=0$, then the map above is bijective. 
    In particular, if $H^1(M, \Z)^{\iota^*}=0$, the moduli space $\tilde{\mathcal{M}}^I(M, \mathfrak{s}, I) \subset \mathcal{C}^I/(\mathcal{G}^{\perp})^I$ is homeomorphic to the fixed point set of $I$ in the framed moduli space of the ordinary Seiberg--Witten equations on $M$. 
\end{proposition}
\begin{proof}
    Firstly we prove the injectivity of the map $\mathcal{C}^I/\mathcal{G}_0^I \to (\mathcal{C}/\mathcal{G}_0)^I$. 
    If $(A, \phi)$ and $(A', \phi')$ are in $C^I$ and there exist a gauge transformation $u \in \mathcal{G}^{\perp}$ such that $u \cdot (A, \phi)=(A', \phi')$. Then we have $A-A'=-u^{-1}du \in i\Omega^1(M)^{-\iota^*}$ and the assumption that the decomposition $\mathcal{G} \cong \mathrm{U}(1) \times \mathcal{G}^{\perp}$ is equivariant under the involtuion $I$ implies $u \in (\mathcal{G}^{\perp})^I$. 

    Secondly, we prove $H^1(M, \Z)^{\iota^*}=0$ implies the map is surjective. 
    From the condition $H^1(M, \Z)^{\iota^*}=0$, we have $(\mathcal{G}^{\perp})^{\iota^*}$ is contractible and for all $u \in (\mathcal{G}^{\perp})^{\iota^*}$ there is a real valued function $f$ on $M$ such that $f=\iota^*f$  and $u=e^{if}$. 
    Let us take $(A, \phi) \in \mathcal{C}$ with $I(A, \phi)=u \cdot (A, \phi)$ for some $u \in \mathcal{G}^{\perp}$. 
    Then $u$ satisfies $u =\iota^*u$ from $I^2=1$. Therefore $u \in (\mathcal{G}^{\perp})^{\iota^*}$ and there is a real function $f$ such that $u=e^{if}$ and $\iota^*f=f$. 
    Let us take $v=e^{if/2}$ and let us consider $v\cdot (A, \phi)$. One can check that $v\cdot (A, \phi) \in C^I$. 
    This proves $\mathcal{C}^I/(\mathcal{G}^{\perp})^I=(\mathcal{C}/\mathcal{G}^{\perp})^I$.
\end{proof}

Next proposition is used in \cref{counting_real_sw}. 

\begin{proposition}\label{MxS1}
    Let $M$ be a closed oriented $3$-manifold with an involtuion $\iota$. Suppose $H^1(M, \Z)=0$ and $(\mathfrak{t}, I)$ be a Real $\mathrm{spin}^c$ structures on $M$ which covers $\iota$. 
    Let us denote $X$ by the $4$-manifold $M \times S^1$. Let $\mathfrak{u}$ be the spin structure on $S^1$ that is preserved by the Lie group action of $S^1$. 
    Then $\mathfrak{t} \times \mathfrak{u}$ is a $\mathrm{spin}^c$ structure on $X$ with the Real structure $I'$ induced by $I$ which covers the involtuion $\iota'=\iota \times id_{S^1}$. 
    Then we have $\lvert deg'(X, \mathfrak{s}, I') \rvert=\lvert deg(M, \mathfrak{t}, I) \rvert$. 
\end{proposition}
\begin{remark}
    In the ordinary case, the Seiberg--Witten invariant of the $4$-manifold with circle action is obtained by Baldridge \cite{MR1946410}. \cref{MxS1} is the Real version of \cite[Theorem B.]{MR1946410} in the case that a $4$-manifold $X$ is direct product of $S^1$ and a $3$-manifold. 
    Thus in our situation, we can regard the solution on a $4$-manifold $M \times S^1$ as the periodic  solution of the Seinberg--Witten flow on $M$. Actually, it is a solution of the $3$-dimensional equations. 
    So that the proof of \cref{MxS1} which we give below is much easier than \cite[Theorem B.]{MR1946410}. 
    However, one can check that the argument in \cite{MR1946410} works in our situation. 
\end{remark}
\begin{proof}[Proof of \cref{MxS1}]
    Since $H^1(M, \Z)^{-\iota^*}=0$, we have $H^1(M \times S^1, \Z)^{-\iota'^*}=0$ and the formal dimension of the Real Seiberg--Witten moduli space is $0$. Thus we can define $\lvert deg'(X, \mathfrak{s}, I') \rvert$ and $\lvert deg(M, \mathfrak{t}, I) \rvert$. Let $\nu \in i \Omega^1(M)^{-\iota^*}$ be a $1$-form such that the framed moduli space $\tilde{\mathcal{M}}^I_{\nu}(M, \mathfrak{t}, I)$ is transversal. 
    Let $\eta \in i\Omega^+(X)^{-\iota'^*}$ be a self-dual $2$-form induced by $\nu$. 
    
    We write $\mathcal{M}_{\eta}^I(X, \mathfrak{s}, I')_0 \subset (\mathcal{C}(X)/\mathcal{G}(X)^{\perp})^I$ for the fixed point set of the $I$-action on the framed moduli space of the ordinary Seiberg--Witten equations. 
    From \cref{coincides_with_two_quotients}, the map 
    \[
        \tilde{\mathcal{M}}^I_{\eta}(X, \mathfrak{s}, I') \to \mathcal{M}^I_{\eta}(X, \mathfrak{s}, I')_0
    \]
    is injective. We also define $\mathcal{M}^I_{\nu}(M, \mathfrak{t}, I)_0 \subset (\mathcal{C}(M)/\mathcal{G}(M)^{\perp})^I$ in the same manner. 
    From \cref{coincides_with_two_quotients}, the map 
    \[
        \tilde{\mathcal{M}}^I_{\nu}(M, \mathfrak{t}, I) \to \mathcal{M}^I_{\nu}(M, \mathfrak{t}, I)_0
    \]
    is bijective since $H^1(M, \Z)=0$. 

From \cref{bijectivity} below, we have the map 
\[
        \pi^* \colon  \mathcal{M}^I_{\nu}(M, \mathfrak{t}, I)_0 \to \mathcal{M}^I_{\eta}(X, \mathfrak{s}, I')_0
    \]
    is bijective, where $\pi \colon X \to M$ is the projection. 

Then in the commutative diagram
\[
\begin{CD}
\tilde{\mathcal{M}}^I_{\nu}(M, \mathfrak{t}, I) @>f>> \mathcal{M}^I_{\nu}(M, \mathfrak{t}, I)_0 \\
@V\pi^*VV @V\pi^*VV \\
\tilde{\mathcal{M}}^I_{\eta}(X, \mathfrak{s}, I') @>g>> \mathcal{M}^I_{\eta}(X, \mathfrak{s}, I')_0
\end{CD}
\]

    $f$ and $\pi^*$ in the right side is bijective and $g$ is injective. Therefore $\pi^*$ in the left side is bijective. 

    Let us prove that $\tilde{\mathcal{M}}^I_{\eta}(X, \mathfrak{s}, I')$ is transversal and the signed count of the elements in $\tilde{\mathcal{M}}^I_{\eta}(X, \mathfrak{s}, I')$ and that of $\tilde{\mathcal{M}}^I_{\nu}(M, \mathfrak{t}, I)$ are the same. 
    If $(A, 0)$ is a reducible solution on $M$, then from the asuumption of the perturbation $1$-form $\nu$, we have that the kernel of the $3$-dimensional Dirac oerator $D^3_A$ in $\Gamma(S)^I$ is $0$. 
    One can check that the reducible soltuion $\pi^*(A, 0)$ on $X$ is transversal directly using $\ker D^3_A=0$. 

    To prove the irreducible solutions on $X$ are transversal and to prove $\lvert deg'(X, \mathfrak{s}, I') \rvert=\lvert deg(M, \mathfrak{t}, I) \rvert$, let us consider the linealization of the Seiberg--Witten equations with gauge fixing condition. 

    We identify $T^*X \cong \R\cdot dt \oplus \pi^* T^*M \cong \R \oplus \wedge^+$ in the usual way. 
    Then both $i\Omega(X)^1\oplus \Gamma(S^+_X)$ and $i\Omega^0(X) \oplus i\Omega^+(X) \oplus \Gamma(S^-_X)$ are identified with the space of the closed path in $i\Omega^0(M) \oplus i\Omega^1(M) \oplus \Gamma(S_M)$.  
    Let $\gamma \in i\Omega(X)^1\oplus \Gamma(S^+_X)$ be an arbitrary configuation given by $\pi^*\gamma'$ for some $\gamma'=(A', \phi') \in i\Omega^1(M) \oplus \Gamma(S_M)$. 
    Suppose $(\xi(t), b(t), \psi(t)) \in i\Omega^0(M) \oplus i\Omega^1(M) \oplus \Gamma(S_M)$ for $t \in \R/2\pi\Z$ is in the kernel of the $4$-dimensional linearized Seiberg--Witten equations at $\gamma$ with gauge fixing condition, i.e.:
    \begin{align*}
        \frac{\partial}{\partial t}\xi(t)&=-d^*b-i\text{Re}(\langle i\phi, \psi \rangle), \\
        \frac{\partial}{\partial t}b(t)&=-d\xi-*db+(\phi \psi^*)_0+(\psi  \phi^*)_0, \\
        \frac{\partial}{\partial t}\psi(t)&=\xi \phi-\Slash{D}_{A'}\psi-\rho(b)\phi. 
    \end{align*}
It is easily seen that the right hand side of the equation is linearization $L$ of the $3$-dimensional Seiberg--Witten equations with the gauge fixing condition.  Since $L$ is a self-adjoint operator on $i\Omega^0(M) \oplus i\Omega^1(M) \oplus \Gamma(S_M)$, the above equation can be written as 
\begin{equation}\label{linearlization}
    \frac{d}{dt} y=Ly
\end{equation}
where $y(t)=(\xi(t), b(t), \psi(t))$. Since $y(0)=y(2\pi)$, the set of the solutions of \eqref{linearlization} coincides with the kernel of $L$. It is clear that the above argument runs for the Real setting. 
 
Therefore if the $S^1$-invariant irreducible solution $\gamma$ on $M \times S^1$ is a transversal in $\mathcal{C}^I(M)/(\mathcal{G}(M)^{\perp})^I$, then $\pi^* \gamma$ is transversal in $\mathcal{C}^I(X)/(\mathcal{G}(X)^{\perp})^I$. 
    
    We are left to prove that the signed counting of $3$-dimensional solutions and $4$-dimensional solutions coincides. 
Let $[\gamma_0]$ be the unique reducible solution in $\tilde{\mathcal{M}}_{\nu}^I(M, \mathfrak{t}, I)$. Let $[\gamma]$ be an arbitrary element in $\tilde{\mathcal{M}}_{\nu}^I(M, \mathfrak{s}, I)$ and $\{[\gamma(t)]\}_{t \in [0, 1]} \subset \mathcal{C}^I/(\mathcal{G}^{\perp})^I$ be a smooth path with $\gamma(0)=\gamma_0$ and $\gamma(1)=\gamma$. 
Suppose the spectral flow of the linealization of the $3$-dimensioanl Real Seiberg--Witten equations with gauge fixed condition is transversal on the path $\gamma(t)$. 

From the argument above, the kernel of the linealization of the $4$-dimensional equations at $\pi^*\gamma(t)$ coincides with the kernel of the linealization of the $3$-dimensional equations. 
This implies the sign of $[\gamma] \in \tilde{\mathcal{M}}_{\nu}^I(M, \mathfrak{t}, I)$ coincides with that of $\pi^* [\gamma] \in \tilde{\mathcal{M}}^I_{\eta}(X, \mathfrak{s}, I')$. 
\end{proof}
\begin{lemma}\label{bijectivity}
    Under the same condition in \cref{MxS1}, we have the map 
    \[
        \pi^* \colon  \mathcal{M}^I_{\nu}(M, \mathfrak{t}, I)_0 \to \mathcal{M}^I_{\eta}(X, \mathfrak{s}, I')_0
    \]
    is bijective, where $\pi \colon X \to M$ is the projection. 
\end{lemma}
\begin{proof}
    One can easily check that reducible solution is unique in $\mathcal{M}^I_{\nu}(M, \mathfrak{t}, I)_0$ and $\mathcal{M}^I_{\eta}(X, \mathfrak{s}, I')_0$ since $H^1(M, \R)^{-\iota^*}=H^1(X, \R)^{-\iota^*}=0$. We will focus on irreducible solutions. 
    Let us show that for each irreducible solution $(A, \Phi)$ of the ordinary Seiberg--Witten equations on $(M \times S^1, \mathfrak{s})$ there is a gauge thransformation $u \in \mathcal{G}(X)^{\perp}$ such that $(A, \Phi)$ is $S^1$-invariant solution. 
    Since $X=M \times S^1$ and $\mathfrak{s}$ is the pullback of the $\mathrm{spin}^c$ structure $\mathfrak{t}$ and the perturbation $\eta$ is induced by the pullback of the $1$-form $\nu$ on $M$, we have the topological energy $\mathcal{E}^{\text{top}}_{\eta}(\cdot)$ of the perturbed Seiberg--Witten equations is identically $0$. 
    So that if $(A, \Phi)$ be a soltuion on $X$, the analytic energy $\mathcal{E}^{\text{an}}_{\eta}(A, \Phi)$ of the Seiberg--Witten equations perturbed by $\eta$ is $0$. (One can fined the definition of these energy functionals in \cite{kronheimer2007monopoles}.) 
    We write $t$ is the periodic coodinate of $S^1\cong \R/2\pi\Z$. From the equation \cite[Section 10.6, (10.12)]{kronheimer2007monopoles}, for the $4$-dimensional Seiberg--Witten solution $(A, \Phi)=(\hat{A}(t)+c dt, \hat{\Phi}(t))$, where $(\hat{A}(t), \hat{\Phi}(t))=(A, \Phi)|_{Y\times \{t\}}$ and $c \colon X \to i \R$, we have 
    \[
     \frac{d}{dt}\hat{\Phi}(t)+c\hat{\Phi}(t)=0
    \]
    and $(\hat{A}(t),\hat{\Phi}(t))$ satisfies $3$-dimensioanl Seiberg--Witten equations perturbed by $\nu$ for all $t \in S^1$. 
    Then $\hat{\Phi}(t)$ is given by
    \[
        \hat{\Phi}(t)=e^{-\int_{0}^{t} c(\tau)d\tau}\hat{\Phi}(0)
    \]
    and $\hat{\Phi}(2\pi)=\hat{\Phi}(0)$. From the unique continuation of $\hat{\Phi}(t)$, we see that $\int_{0}^{2\pi} c(\tau)d\tau \in 2\pi\Z$. Thus $u= e^{\int_{0}^{t} c(\tau)d\tau}$ is a gauge transformation on $X$. One can easily check that $u \cdot (A, \Phi)$ is a temporal gauge and is pullback of a solution of the $3$-dimensional equations. 
    
    Suppose $[(A, \Phi)]=\pi^*[(B, \phi)] \in \mathcal{M}^I_{\eta}(X, \mathfrak{s}, I')_0$. Then there exist gauge transformation $u \in \mathcal{G}(X)^{\perp}$ on $X$ such that $I'(A, \Phi)=u \cdot (A, \Phi)$. Since $(A, \Phi)$ is $S^1$-invariant and $I$ commutes with the $S^1$-action, $u$ is $S^1$-invariant gauge transformation. Then $[(B, \phi)] \in \mathcal{M}^I_{\nu}(M, \mathfrak{t}, I)_0$ and the map
    \[
       \pi^* \colon  \mathcal{M}^I_{\nu}(M, \mathfrak{t}, I)_0 \to \mathcal{M}^I_{\eta}(X, \mathfrak{s}, I')_0
    \]
    is surjective. 
    One can easily check that the map $\pi^*$ is injective because if there exist a gauge transformation $u \in \mathcal{G}^{\perp}$ on $X$ such that $\pi^*(B, \phi)=u \cdot \pi^*(B', \phi)$ for some $[(B, \phi)]$ and $[(B', \phi)]$ in  $\mathcal{M}^I_{\nu}(M, \mathfrak{t}, I)_0$, then $u$ is $S^1$-invariant. 
\end{proof}

\section{Examples and Calculations}\label{Examples_of_surfaces}
\subsection{Surfaces in $K3\#\overline{\C P^2}\# N \C P^2$}
In this subsection we prove \cref{main1}. 
\begin{proposition}\label{Freedman_Perron_Quinn}
    Let $f \colon K3\#\overline{\C P^2} \to 3\C P^2 \# 20 \overline{\C P^2}$ be a homeomorphism.
    If $N\ge 1$,   
    there is a diffeomorphism 
    \[\Phi \colon K3\#\overline{\C P^2}\# N \C P^2 \to 3\C P^2 \# 20 \overline{\C P^2} \# N \C P^2\] 
    such that $\Phi$ is $C^0$-isotopic to the homomorphism $f \# N id_{\C P^2}$. 
\end{proposition}
\begin{proof}
    It is shown that if $N\ge 1$, $K3\#\overline{\C P^2}\# N \C P^2 $ and $ 3\C P^2 \# 20 \overline{\C P^2} \# N \C P^2$ are diffeomorphic. Let us take a diffeomorphism     \[\Phi' \colon K3\#\overline{\C P^2}\# N \C P^2 \to 3\C P^2 \# 20 \overline{\C P^2} \# N \C P^2. \] 
Let us denote $\psi$ be the automorphism on $H^2(3\C P^2 \# 20 \overline{\C P^2} \# N \C P^2, \Z)$ given by $(f \# N id_{\C P^2} \circ {\Phi'}^{-1})^*$. From the Wall's theorem for the diffeomorphisms \cite[Theorem.2]{MR0163323}, we see that there is a self-diffeomorphism $\Psi$ on $ 3\C P^2 \# 20 \overline{\C P^2} \# N \C P^2$ which satisfies $\Psi^*=\psi$ because the intersection form of $3\C P^2 \# 19 \overline{\C P^2} \# (N-1) \C P^2$ is indefinite and $3\C P^2 \# 19 \overline{\C P^2} \# (N-1) \C P^2 \# S^2 \times S^2$ is diffeomorphic to $ 3\C P^2 \# 20 \overline{\C P^2} \# N \C P^2$. 
We shall define the diffeonmorphism $\Phi$ using the composition $\Psi \circ \Phi'$. This is $C^0$ isotopic to $f \# N id_{\C P^2}$ because $(f \# N id_{\C P^2} \circ \Phi^{-1})^*=(f \# N id_{\C P^2} \circ {\Phi'}^{-1} \circ \Psi^{-1})^*$ is the identity on $H^2(3\C P^2 \# 20 \overline{\C P^2} \# N \C P^2, \Z)$ and from the theorem of Quinn \cite[Theorem 1.1]{MR0868975}, this is $C^0$ isotopic to the identity. 
\end{proof}
\begin{definition}\label{exotic_surfaces}
    Let $S_0$ be an embedded $\C P^1$ in $\C P^2$ that homology class is the twice of the generator of $H_2(\C P^2, \Z)$. Let $S'$ be an embedded $2$-sphere in $\# N \C P^2$ given by taking inner connected sums to $S_0$ in each connected component of $\# N \C P^2$. 
    Let $S$ be an embedded $2$-sphere in $3\C P^2 \# 20 \overline{\C P^2} \# N \C P^2$ given by the image of the composition of the embedding 
    \[
        S' \to \#N \C P^2 \setminus D^4 \to 3\C P^2 \# 20 \overline{\C P^2} \# N \C P^2. 
    \]
    Let $\tilde S$ be an embedded $2$-sphere in $3\C P^2 \# 20 \overline{\C P^2} \# N \C P^2$ given by the image of the following embedding:
    \[
    \begin{CD}
        S' \to \#N \C P^2 \setminus D^4 \to K3 \# \overline{\C P^2} \# N \C P^2  @>{\Phi}>> 
    \end{CD} 3\C P^2 \# 20 \overline{\C P^2} \# N \C P^2. 
    \]
\end{definition}
\begin{lemma}\label{homeo}
    The embedded surfaces $S$ and $\tilde S$ are $C^0$-isotopic. 
\end{lemma}
\begin{proof}
    From \cref{Freedman_Perron_Quinn}, we have that $\Phi$ is $C^0$-isotopic to the homeomorphism $f \# N id_{\C P^2}$. We see that two embeddings
    \[
        S' \to \#N \C P^2 \setminus D^4 \to 3\C P^2 \# 20 \overline{\C P^2} \# N \C P^2
    \]
    and 
    \[
        \begin{CD}
            S' \to \#N \C P^2 \setminus D^4 \to K3 \# \overline{\C P^2} \# N \C P^2  @>{f \# N id_{\C P^2}}>> 
        \end{CD} 3\C P^2 \# 20 \overline{\C P^2} \# N \C P^2
    \]
    are the same embeddings. This proves the lemma. 
\end{proof}

\begin{theorem}\label{exoticness}
    The pairs $(3\C P^2 \# 20 \overline{\C P^2}, S)$ and $(3\C P^2 \# 20 \overline{\C P^2}, \tilde S)$ are not diffeomorphic.
\end{theorem}
\begin{proof}
    Let us compare the invariant of the pairs $(3\C P^2 \# 20 \overline{\C P^2}\#N \C P^2 , S)$ and $(3\C P^2 \# 20 \overline{\C P^2}\#N \C P^2 , \tilde S)$. 

\underline{Calculations of the invariant of $S$}
We see that the double branched cover along $S$ is 
\[
    3\C P^2 \# 20 \overline{\C P^2}\#N S^2 \times S^2 \# 3\C P^2 \# 20 \overline{\C P^2}. 
\]
The covering involution $\iota$ is given by the $\Z/2$-equivariant connected sums of 
\[
    S^2 \times S^2 \ni (x, y) \mapsto (y, x) \in S^2 \times S^2
\]
 and the switching the connected components $3\C P^2 \# 20 \overline{\C P^2}\setminus D^4 \sqcup 3\C P^2 \# 20 \overline{\C P^2}\setminus D^4$. 
 Since the $4$-manifold $3\C P^2 \# 20 \overline{\C P^2}\#N S^2 \times S^2 \# 3\C P^2 \# 20 \overline{\C P^2}$ is simply connected, if we take a $\mathrm{spin}^c$ structure $\mathfrak{s}$ with $\iota^*c_1(\mathfrak{s})=-c_1(\mathfrak{s})$ we have unique Real structure on $\mathfrak{s}$ that covers $\iota$. 
 From the definition of $\iota$, all Real $\mathrm{spin}^c$ structures are given by the following way: Let $\mathfrak{s}_0$ be the spin structure on $S^2 \times S^2$ and $\mathfrak{s}_1$ is an arbitrary $\mathrm{spin}^c$ structure on $3\C P^2 \# 20 \overline{\C P^2}$. 
 We denote $\mathfrak{s}$ by a $\mathrm{spin}^c$ structure on $3\C P^2 \# 20 \overline{\C P^2}\#N S^2 \times S^2 \# 3\C P^2 \# 20 \overline{\C P^2}$ that is given by $\mathfrak{s}_1 \cup \mathfrak{s}_0 \cup \mathfrak{s}_1$. 
 
 We see that $\iota^*c_1(\mathfrak{s})=-c_1(\mathfrak{s})$ and we have the unique Real structure $I$. 
 From \cref{gluing_disconn}, we see that the non-equivariant Real Bauer--Furuta invariant is given by 
 \begin{align*}
    &\Psi(3\C P^2 \# 20 \overline{\C P^2}\#N S^2 \times S^2 \# 3\C P^2 \# 20 \overline{\C P^2}, \mathfrak{s}, I) \\& =\epsilon(\Psi(3\C P^2 \# 20 \overline{\C P^2}\setminus B^4, \mathfrak{s}_1) \wedge \Psi(\# N S^2 \times S^2 \setminus B^4 \sqcup B^4, \mathfrak{s}_0, I))\\
    & = \Psi(3\C P^2 \# 20 \overline{\C P^2}, \mathfrak{s}_1) \wedge \Psi(\# N S^2 \times S^2 , \mathfrak{s}_0, I). 
 \end{align*}
 Since $3\C P^2 \# 20 \overline{\C P^2}$ has positive scalar curvature, we have that $\Psi(3\C P^2 \# 20 \overline{\C P^2}, \mathfrak{s}_1)=0$ for all $\mathfrak{s}_1$. Thus we have $\Psi(3\C P^2 \# 20 \overline{\C P^2}\#N S^2 \times S^2 \# 3\C P^2 \# 20 \overline{\C P^2}, \mathfrak{s}, I)=0$ for all Real $\mathrm{spin}^c$ structure $(\mathfrak{s}, I)$. 

\underline{Calculations of the invariant of $\tilde S$}
Let $\tilde S'$ be an embedded $2$-sphere given by the following embedding:
    \[
        S' \to \#N \C P^2 \setminus D^4 \to K3 \# \overline{\C P^2} \# N \C P^2 .
    \]
    We see that $(K3 \# \overline{\C P^2} \# N \C P^2, \tilde S')$ and $(3\C P^2 \# 20 \overline{\C P^2}, \tilde S)$ are diffeomorphic and the invariant of $(K3 \# \overline{\C P^2} \# N \C P^2, \tilde S')$ and $(3\C P^2 \# 20 \overline{\C P^2}, \tilde S)$ are the same. 
We can easily see that the double branched cover of $K3 \# \overline{\C P^2} \# N \C P^2$ along $\tilde{S}'$ is 
\[
    K3 \# N S^2 \times S^2 \# K3
\]
and the covering involution is similar in the case \underline{Calculations of the invariant of $S$}. 
Form the same argument in  
\underline{Calculations of the invariant of $S$}, we see that all Real $\mathrm{spin}^c$ structure $(\mathfrak{s}, I)$ are given by $\mathfrak{s}=\mathfrak{s}_1 \cup \mathfrak{s}_0 \cup \overline{\mathfrak{s}_1}$ and the unique Real involution $I$. Using \cref{gluing_disconn}, we see that 
\begin{align*}
    &\Psi(K3\#N S^2 \times S^2 \# K3, \mathfrak{s}, I) \\& =\epsilon(\Psi(K3 \setminus B^4, \mathfrak{s}_1) \wedge \Psi(\# N S^2 \times S^2 \setminus B^4 \sqcup B^4, \mathfrak{s}_0, I))\\
    & = \Psi(K3 , \mathfrak{s}_1) \wedge \Psi(\# N S^2 \times S^2 , \mathfrak{s}_0, I) . 
 \end{align*}
 If we take the $\mathrm{spin}^c$ structure $\mathfrak{s}_1$ on $K3$ that comes from the spin structure on $K3$, we have $\Psi(K3 , \mathfrak{s}_1)$ is the generator of $\pi_1^{st}(S^0) \cong \Z/2$. From \cref{the_invariant_of_S2xS2} below, we have that 
 \[
    \Psi(K3\#N S^2 \times S^2 \# K3, \mathfrak{s}, I)= 1 \in \pi_1^{st}(S^0) \cong \Z/2=\{0, 1\}. 
 \]
 This proves the theorem. 
\end{proof}

\begin{lemma}\label{the_invariant_of_S2xS2}
    Let $\iota $ be the involution on $\# N S^2 \times S^2$ that is given by the equivariant connected sums of the involution
    \[
        S^2 \times S^2 \ni (x, y) \mapsto (y, x) \in S^2 \times S^2. 
    \]
    Then we have that the unique Real $\mathrm{spin}^c$ structure $(\mathfrak{s}_0, I)$ is the pair that comes from the spin structure with the unique Real involution. The non-equivariant Real Bauer--Furuta invariant $\Psi(\# N S^2 \times S^2 , \mathfrak{s}_0, I)$ is a generator of $\pi_0^{st}(S^0) \cong \Z$. 
\end{lemma}
\begin{proof}
    It is easy to check that the involution $\iota$ acts trivially on the cohomology ring of $\# N S^2 \times S^2$. 
    If $(\mathfrak{s},I)$ be a Real structure that covers $\iota$, we have $\iota^* c_1(\mathfrak{s})=-c_1(\mathfrak{s})$. This implies that $\mathfrak s$ should be $\mathfrak{s}_0$. Since $\# N S^2 \times S^2$ is simply connected, there is a unique Real involution $I$ on $\mathfrak{s}$. 

    Let us calculate the invariant $\Psi(\# N S^2 \times S^2 , \mathfrak{s}_0, I)$. The formal dimension of the moduli space is $0$. 
We see that $b^+(\# N S^2 \times S^2)=b^+_{\iota}(\# N S^2 \times S^2)$ and $\# N S^2 \times S^2$ have a positive scalar curvature. 
Therefore we have that the only solution of the Real Seiberg--Witten equation is reducible and this is a non-degenerate solution. 
The element $\pi_0^{st}(S^0)$ is characterized by its mapping degree. 
We see that the mapping degree of the Real Bauer--Furuta invariant $\Psi(\# N S^2 \times S^2 , \mathfrak{s}_0, I)$ is $\pm 1$. This proves the lemma. 
\end{proof}
\subsection{Twist and rolling spun knots}\label{Twist and rolling spun knots}
    In this section, we see that the calculations of the invariant of the $2$-knot, using the construction of the twist and rolling spun knot given by Plotnick \cite{MR0780592}. Firstly, we recall the constion given in \cite[Section 1]{MR0780592}. 
\subsubsection{Construction of $k$-twist $\alpha$-roll spun knot}

    Let $K$ be a $1$-knot in $S^3$. Let $\nu(K)$ be a tubular neighborhood of the knot $K$ and $M_K :=S^3\setminus \nu(K)$. Denote $X_K:=M_K \times S^1$. Note that $\partial X_K=-\partial \nu(K) \times S^1 \cong T^3$. 
    
\begin{definition}\label{basis}
    Let $m$ and $l$ be embedded circles in $\partial M_k$ that represents meridian and longitude respectively. We denote $\{p\} \times S^1 \subset M_K \times S^1=X_K$ by $h$, where $p \in \partial M_K$. 
\end{definition}
One can easily prove the following lemma:
\begin{lemma}
    The first homology of $X_K$ is $\Z \langle m \rangle \oplus \Z \langle h \rangle$ for all knot $K$. 
\end{lemma}

\begin{definition}\label{twin}
    Let $P$ be a $4$-manifold with boundary given by the following way:
\begin{itemize}
    \item Let $\{(\psi, \theta) \mid -\pi \le \psi \le \pi,\; 0 \le \theta \le 2\pi\}$ be a spherical coordinate of $S^2$ and $\{(r, \varphi) \mid 0 \le r \le 1, \; 0 \le \varphi \le 2\pi \}$ be a polar coordinate of $D^2$. 
    \item Let $(\varphi, t,(s, \theta)), \; -1\le t \le 1, 0 \le \varphi\le 2\pi$ be a coordinate of $S^1 \times [0,1] \times D^2$, where $(s, \theta)$ is the polar coordinate of $D^2$. 
    \item Let us denote $D_0^3 \subset S^2 \times D^2$ by the set $\{ ((\psi, \theta), (1, \varphi)) \in S^2 \times D^2\mid 0\le \psi \le \pi/4\}$. We denote $D_1^3=\{ ((\psi, \theta), (1, \varphi)) \in S^2 \times D^2\mid 3\pi/4 \le \psi \le \pi\}$. 
    \item Let $\alpha_i$ for $i=0, 1$ be a homeomorphism that are given as follows:
    \[
        \alpha \colon D_i^3 \to S^1 \times \{i\} \times D^2 ;\; ((\psi, \theta), (1, \varphi)) \mapsto (\varphi, (-1)^i, (\sqrt{2}\sin \psi, \theta)). 
    \]
    \item Let $P$ be a $4$ manifold given by the smoothing of the corner of 
    \[
        S^2 \times D^2 \cup_{\alpha_0, \alpha_1} S^1 \times [0, 1] \times D^2. 
    \]
\end{itemize}
\end{definition}
\begin{definition}
    Let $P$ be the $4$-manifold we defined in \cref{twin}. Let $S \subset P$ the embedded 
    $2$-sphere in $P$ given by $\{ ((\psi,0),(r, \varphi))\mid 0\le r \le 1, \; 0\le \varphi \le 2\pi, \; \psi=0, \pi\} \cup S^1 \times I \times \{0\}$. We call this $2$-sphere $S$ by the core of $P$. 
\end{definition}

We can easily check that the boundary of $P$ is $T^3$. We define  generators of $\pi_1(\partial P)\cong H_1(\partial P, \Z) \cong \Z^3$. 
\begin{definition}
Let us take a base point $(0, 0, (1, 0)) \in \partial S^1 \times [-1, 1] \times D^2$. we define $3$ loops in the boundary of $P$ as follows:
    \begin{align*}
        e_1 & =\{ (0, 0, (1, \theta)) \in S^1 \times [-1,1] \times D^2 \mid 0\le \theta \le 2\pi \}, \\
        e_2 & =\{ (\varphi, 0, (1, 0)) \in S^1 \times [-1,1] \times D^2 \mid 0\le \varphi \le 2\pi \}, \\
        e_3 & =\{ (0, t, (1, 0)) \in S^1 \times [-1,1] \times D^2 \mid -1 \le t \le 1 \} \cup \{ ((\psi, 0), (1, 0)) \in S^2 \times D^2 \mid \frac{\pi}{4} \le \psi \le \frac{3\pi }{4}\}. \\
    \end{align*}

\end{definition}
\begin{remark}
    Note that $e_1$ is the meridian of the core $S$ of $P$. 
\end{remark}
One can easily check that the following lemma holds. 
\begin{lemma}
    Let $P$ be the $4$-manifold given as above. Then we have
    the map
    \[
        H_1(\partial P) \to H_1(P)
    \] 
    induced by the inclusion is surjective and the kernel is $\Z\langle e_1 \rangle \oplus \Z\langle e_2 \rangle$. Moreover, $e_3$ is an generator of $H_1(P) \cong \Z$. 
\end{lemma}

Let us define the construction of a $2$-knot that is called $k$-twist $\alpha$-roll spun of the knot $K$. One can find the proof of the following lemma in \cite[Section 1, Section 2]{MR0780592}.  
\begin{lemma}\label{well_def_k_twist_a_roll}
    Let $k$ and $\alpha$ be integers. Let $A \in SL_3(\Z)$ be a matrix
    \[A=
        \begin{pmatrix}
            1&k&0\\
            0&1&0\\
            0&\alpha&1
        \end{pmatrix}.
    \]    
    Let $g(A) \colon \partial P \to \partial X_K$ be a difeomorphism 
    that $g(A)_* \colon H_1(\partial P) \to H_1(\partial X_K)$ is $A$ if we take the basis of $H_1(\partial P)$ by $\langle e_1, e_2, e_3 \rangle$ and the basis of $H_1(\partial X_K)$ by $\langle m, h, l \rangle$. 
Then we have that the $4$-manifold $\Sigma_A=P\cup_{g(A)} X_K$ is diffeomorphic to $S^4$. 
\end{lemma}
From \cref{well_def_k_twist_a_roll}, we can define a $2$-knot in $S^4$. 
\begin{definition}\label{k_twist_a_roll}
    Let $S \subset P$ be the core in $P$. Let $i \colon P \to \Sigma_A \cong S^4$ be the natural inclusion.  We call the image $i(S)$ the $k$-twist $\alpha$-roll spun of the knot $K$. We denote this $2$-knot by $\tau_{k, \alpha}(K)$. 
\end{definition}

\subsubsection{Double branched cover and Real spin$^c$ structures}
    Next we see that if we put an additional assumption on $k$ and $\alpha$, we can construct the double branched cover $\Sigma_2(S^4, \tau_{k, \alpha}(K))$ from the double branched cover of $P$ along $S$ and the double cover of $X_K$. Let us firstly consider that the double (branched) cover of $P$ and $X_K$.

\begin{lemma}
    We have a double branched cover of $P$ along $S$. 
    Let us denote $\tilde{P}$ by the double branched cover along $S$. Then $\partial \tilde{P}$ be a double cover of $\partial P$ that is given by the map
    \[
        \pi_1(\partial P) \cong H_1(\partial P) \to \Z/2;\;e_1 \mapsto 1, \; e_2, e_3 \mapsto 0. 
    \]
\end{lemma}
\begin{proof}
We have that the meridian of $S=\{ ((\psi,0),(r, \varphi))\mid 0\le r \le 1, \; 0\le \varphi \le 2\pi, \; \psi=0, \pi\} \cup S^1 \times I \times \{0\}$ is $e_1=\{ (0, 0, (1, \theta)) \in S^1 \times [0,1] \times D^2 \mid 0\le \theta \le 2\pi \}$. Note that $P \setminus S$ is given by $(S^2 \times D^2 \setminus \{ p_N, p_S \} \times D^2) \cup_{S^0 \times (D^2\setminus \{0\} \times S^1)} S^1 \times I \times (D^2 \setminus \{0\})$. From the Mayer--Vietoris exact sequence, we have 
\[
    H_1(P \setminus S) = \Z \langle e_1 \rangle \oplus \Z \langle e_3 \rangle. 
\]
Therefore we can take the double cover of $P\setminus S$ that is given by the map 
\[
    H_1(P \setminus S) \mapsto \Z/2;\;e_1 \mapsto 1, \; e_3 \mapsto 0.
\]
The normal bundle of $S$ is trivial in $P$, we can fill the double cover of $P\setminus S$ and this is the double branched cover $\tilde{P}$. 
The latter half of the statement follows from the definition of the double cover. 
\end{proof}
    
    We see that $M_k$ is a homology $S^1 \times D^2$ so that there is the unique double covering of $M_k$. 
    That double covering is denoted by $\tilde{M}_K$ and the covering involution is denoted by $\iota'$. 
    We consider the manifold $\tilde{X}_K :=\tilde{M}_K \times S^1$ with the involution $\iota:=\iota' \times id_{S^1}$. 
    The proof of the following lemma is obvious from the definition of $\tilde{X}_K$. 

\begin{lemma}
    The double cover $\partial \tilde{X}_K \to \partial X_K$ is given by the map 
    \[
        H_1(\partial X_K) \to \Z/2;\:m \mapsto 1, \; h, l \mapsto 0. 
    \]
\end{lemma}

We give a sufficient and necessary condition of $A$ that there is a lift of $\tilde{A}$ such that $g(\tilde{A}) \colon \partial \tilde{P} \to \partial \tilde{X}_K$ is a lift of $g(A)$. 
\begin{proposition}\label{k/2+a_even}
    Let $A \in SL_3(\Z)$ be the matrix given in \cref{well_def_k_twist_a_roll}. If $k$ is even integer, we have a lift $\tilde{A} \colon H_1(\partial \tilde{P}) \to H_1(\partial \tilde{X}_K)$ of $A$ and a lift $g(\tilde{A})$ of $g(A)$. The lift $g(\tilde{A})$ can be taken to be commutative with the covering involutions. 
\end{proposition}
\begin{proof}
    Let us assume there is a commutative diagram
    \[
\begin{CD}
H_1(\partial \tilde{P}) @>\tilde{A}>> H_1(\partial \tilde{X}_K) \\
@V\pi VV @V\pi VV \\
H_1(\partial P) @>A>> H_1(\partial X_K)
\end{CD},
\]
    where $\pi$ are the projection of the double covers. 
    Note that all three manifolds in the above diagram is $T^3$. 
    Let us take a basis $\langle \tilde{e}_1, e_2, e_3 \rangle$ of $H_1(\partial \tilde{P})$, where $\tilde{e}_1$ is a two fold cover of $e_1 \subset \partial P$. Also we take a basis $\langle \tilde{m}, h, l \rangle$ of $H_1(\partial \tilde{X}_K)$, where $\tilde{m}$ is the two fold cover of the meridian $m \subset \partial X_K$. 
    Then the matrix representation of $\pi$ is 
    \[
        \begin{pmatrix}
            2&0&0\\
            0&1&0\\
            0&0&1
        \end{pmatrix}. 
    \]
    Therefore, the matrix $\tilde{A}$ is given by 
    \[
        \begin{pmatrix}
            1&\frac{k}{2}&0\\
            0&1&0\\
            0&\alpha&1
        \end{pmatrix}. 
    \]
    and $k$ must be in $2\Z$. Conversely, if we give $\tilde{A}$ as above, the diagram above commutes. The latter half of the statement follows easily. 
\end{proof}
From \cref{k/2+a_even}, we have the following corollary;
\begin{corollary}\label{gluing_double_branched_cover}
    Let $k$ be an even integer. Then the double branched cover of $S^4$ along $\tau_{k, \alpha}(K)$ is given by $\tilde{P} \cup_{g(\tilde{A})} \tilde{X}_K$. 
\end{corollary}
Let us denote by $\Sigma_2(S^4, \tau_{k, \alpha}(K))$ that branched cover. 

From now on, we assume $k=2k'$ be an even integer. 

Next we consider Real Spin$^c$ structures on $\tilde{P}$, $\tilde{X}_K$, and $\Sigma_2(S^4, \tau_{k, \alpha}(K))$. 
For simplicity, we consider the case that $(S^4, \tau_{k, \alpha}(K))$ is a nice pair as defined in \cref{nice_emb_closed}. 
\begin{lemma}\label{H_1}
    If the Alexander--polynomial $\Delta_K(t)$ of the knot $K$ satisfies $\lvert \Delta_K(-1) \rvert=1$ and $k$ be an even integer, then $(S^4, \tau_{k, \alpha}(K))$ is a nice pair. 
\end{lemma}
\begin{proof}
    Since $H^2(S^4;\Z/2)=0$, the condition that $[\tau_{k, \alpha}(K)]_2=0$ is always satisfied. We show $H^1(\Sigma_2(S^4, \tau_{k, \alpha}(K));\Z/2)=0$. 

    Note that $\lvert \Delta_K(-1) \rvert=1$ implies that the double branched cover of $S^3$ along the knot $K$ is integer homology three sphere. Then $H_1(\tilde{X}_K, \Z) \cong \Z \langle \tilde{m} \rangle \oplus \Z \langle h \rangle $. We have $H_1(\tilde{P}) = \Z \langle e_3 \rangle$. Using the Mayer--Vietoris exact sequence
    \[
       \dots \to H_1(\partial \tilde{P};\Z) \to H_1(\tilde{P};\Z) \oplus H_1(\tilde{X}_K;\Z) \to H_1(\Sigma_2(S^4, \tau_{k, \alpha}(K));\Z)\to 0
    \]
and from the fact that, we have that the map $H_1(\partial \tilde{P};\Z) \to H_1(\tilde{P};\Z) \oplus H_1(\tilde{X}_K;\Z)$ is surjection. Therefore $H_1(\Sigma_2(S^4, \tau_{k, \alpha}(K));\Z)=0$.
From the universal coefficient theorem, we have $H^1(\Sigma_2(S^4, \tau_{k, \alpha}(K));\Z/2)=0$.  
\end{proof}

\begin{lemma}\label{real_spinc_2_knot}
    Let $K$ be a knot with $\lvert \Delta_K(-1)\rvert=1$ and $k$ be an even integer. Then the unique Real Spin$^c$ structure on $\Sigma_2(S^4, \tau_{k, \alpha}(K))$ is the unique Real spin structure. 
\end{lemma}
\begin{proof}
    We see that $b_2(\Sigma_2(S^4, \tau_{k, \alpha}(K)))=0$ and from the proof of \cref{H_1}, we have $H_1(\Sigma_2(S^4, \tau_{k, \alpha}(K));\Z)=0$ then $H^2(\Sigma_2(S^4, \tau_{k, \alpha}(K));\Z)=0$. Therefore the only spin$^c$ structure on $\Sigma_2(S^4, \tau_{k, \alpha}(K))$ is the spin structure. 
    Note that the spin structure is unique since $H^1(\Sigma_2(S^4, \tau_{k, \alpha}(K));\Z/2)=0$. 
    Let $\mathfrak{s}$ be the spin structure. This is obviously preseved by the involtuion $\iota$ so that there is a Real structure $I=\tilde{\iota}\circ j$, where $\tilde{\iota}$ is a lift of the involtuion $\iota$ to the spinor bundle. 
    Note that there is two lift $\pm \tilde{\iota}$ however $iI(-i)=-I$ implies these Real spin$^c$ structures are isomorphic. 
    Since $H^1(\Sigma_2(S^4, \tau_{k, \alpha}(K));\Z)^{\iota^*}=0$, so that the isomorphism class on $\mathfrak{s}$ is unique. 
\end{proof}

Let us consider that Real spin structures on $\tilde{P}$ and $\tilde{X}_K$. 
One can easily see that the involutions on $\tilde{P}$, $\tilde{X}_K$, and $\partial \tilde{P}$ act trivially on $H_1(\cdot, \Z/2)$. 
Therefore we have all the spin structure have Real structure. 
We only have to consider that which spin structure on $\tilde{P}$ and $\tilde{X}_K$ are glued by the map $g(\tilde{A})$. 

\begin{definition}
    Let $\mathfrak{s}'_0$ be the spin structure on $\tilde{M}_K$, which is the double cover of $S^3\setminus \nu(K)$, that can be extended to $\Sigma_2(S^3, K)$. 
    Let us define $\mathfrak{s}_{\tilde{X}_K}$ is the product spin structure $\mathfrak{s}'_0 \times \mathfrak{t}_1$ on $\tilde{X}_K = \tilde{M}_K \times S^1$, where $\mathfrak{t}_1$ is the spin structure invariant under $S^1$ action. 
    Note that $\mathfrak{t}_1$ is a generator of $1$-st spin bordism group. 
\end{definition}
We prove later that $\mathfrak{s}_{\tilde{X}_K}$ is the only Real spin structure on $\tilde{X_K}$ whose Real relative Bauer--Furuta invariant can be nontrivial.  

\begin{proposition}\label{condition_of_alpha}
     Let $k$ be an even integer and $K$ be a knot with $\lvert \Delta_K(-1) \rvert=1$. If $\frac{k}{2}+\alpha$ is an odd integer, the spin structure $\mathfrak{s}_{\tilde{X}_K}$ can be extended to $\Sigma_2(S^4, \tau_{k, \alpha}(K))$. 
\end{proposition}
\begin{proof}
    Let us denote $\mathfrak{s}_{T^3}$ be the spin structure on \[\partial \tilde{P} \cong (\R\langle \tilde{e}_1 \rangle \oplus\R\langle e_2 \rangle \oplus \R\langle e_3 \rangle)/(\Z \langle \tilde{e}_1 \rangle \oplus\Z\langle e_2 \rangle \oplus \Z \langle e_3 \rangle) \cong  T^3\] that is invariant under $T^3$ action. 

    Let us take a basis $\epsilon_1, \epsilon_2$, and $\epsilon_3$ of $H^1(T^3; \Z/2)$ by the dual basis of $\tilde{e}_1, e_2$, and $e_3$. 
    Then if we take the identification $\partial \tilde{X}_K \cong T^3$ by $\tilde{m} \mapsto \tilde{e}_1, h \mapsto e_2$, and $l \mapsto e_3$, we have $\mathfrak{s}_{\tilde{X}_K}|_{\partial \tilde{X}_K}=\mathfrak{s}_{T^3} \otimes_{\R} l_{\epsilon_1+\epsilon_3}$, where $l_{\gamma}$ is a real line bundle on $T^3$ whose first Stiefel--Whitney class is $\gamma \in H^1(T^3; \Z/2)$. 
    Spin structures on $\partial \tilde{P}$ which we can extend it to $\tilde{P}$ are isomorphic to a spin structure given by  $\mathfrak{s}_{T^3} \otimes_{\R} l_{\epsilon_1+\epsilon_2+a\epsilon_3}$, where $a \in \Z/2$. 

    The matrix representation of $g(\tilde{A})^* \colon H^1(\partial \tilde{X}_K;\Z/2) \to H^1(\partial \tilde{P};\Z/2)$ is given as follows: If we identify $\partial \tilde{X}_K \cong T^3$ as above and we use the basis of $H^1(T^3;\Z/2)$ by $\langle \epsilon_1, \epsilon_2, \epsilon_3 \rangle$ as above, then the matrix is 
    \[
        \begin{pmatrix}
            1&0&0\\
            \frac{k}{2}&1&\alpha\\
            0&0&1
        \end{pmatrix}. 
    \]
    Then $g(\tilde{A})^*(\epsilon_1+\epsilon_3)=\epsilon_1+(\frac{k}{2}+ \alpha)\epsilon_2+\epsilon_3$. If we have $\frac{k}{2}+ \alpha=1 \in \Z/2$, then we $g(\tilde{A})^*l_{\epsilon_1+\epsilon_3}=l_{\epsilon_1+\epsilon_2+\epsilon_3}$. Thus we finish the proof. 
\end{proof}

\subsubsection{Calculations of the Real Bauer--Furuta invariants}\label{calculations_subsubsection}

Then we will calculate the Real Bauer--Furuta invariant using the gluing result \cref{gluing_conn}. Firstly, We compute the Real Seiberg--WItten Floer homotopy type of $T^3$ with free involution. 

\begin{lemma}\label{real_spin_T3}
    Let $\iota \colon T^3 \to T^3$ be an involution that is given by $(x, y, z) \mapsto (x+\pi, y, z) \in T^3=(\R/2 \pi \Z)^3$. Let $\mathfrak{s}$ be a spin structure on $T^3$ that is extendable to $D^2 \times T^2$.  Then the lift of this involution $\tilde{\iota}$ to the spin structure $\mathfrak{s}$ is order $4$. 
    In particular, we have a canonical Real spin$^c$ structure on $\mathfrak{s}$ that is given by $I=\tilde{\iota}\circ j$. 
    Moreover, all Real spin$^c$ structure on $T^3$ with the first Chern class is $0$ is given by $(\mathfrak{s} \otimes_{\R}l, I)$, where $l$ is any real line bundle on $T^3$. 
\end{lemma}
\begin{proof}
    Note that the involution $\iota$ is extended to $D^2 \times T^2$ by $((r, x), y, z) \mapsto ((r, x+\pi), y, z)$, where $\{(r, x)\mid 0\le r \le 1, x \in \R/2\pi\Z \}$ be a polar coordinate of $D^2$. 
    From the assumption of $\mathfrak{s}$, there is a spin structure $\tilde{\mathfrak{s}}$ on $D^2 \times T^2$ with $\tilde{\mathfrak{s}}|_{T^3}=\mathfrak{s}$. 
    Note that all the spin structures on $D^2 \times T^2$ are preserved by the involution. The fixed point set of the involution on $D^2 \times T^2$ is codimension $2$. 
    Then from \cite[Lemma 4]{MR4365044}, then we have that the lift $\tilde{\iota}$ is order $4$. 
    The latter half of the statement is easily follows from the classification result of real spin$^c$ structure \cref{classification_of_real_spinc}. 
\end{proof}

\begin{proposition}\label{real_seiberg_witten_of_T3}
    Let $\mathfrak{s}$ be a spin structure on $T^3$ that can be extended to $D^2 \times T^2$. Then the Real Seiberg--Witten floer homotopy type of $(T^3, (\mathfrak{s} , I))$ is 
    \[
       SWF_G(T^3,\mathfrak{s},I) = (S^0, 0, 0)   
    \] 
    for $G=\Z/4$ or $\Z/2$. 
\end{proposition}
\begin{proof}
    We take the flat metric on $T^3$. Then all solutions of the Seiberg--Witten equations are reducible. Since $H^1(T^3, \R)^{-\iota^*}=0$, we have that there is only one reducible in each real spin$^c$ structure. We have that the kernel of the spin Dirac operator on the spin structure $\mathfrak{s}$ is $0$ since we have $-1$ monodromy along the loop $\{(x, 0, 0)\mid x \in \R/2\pi\Z\} \subset T^3$. Therefore the reducible solution is transverse. Then by definition of the conley index and we have the Real Seiberg--Witten Floer homotopy type is given by 
    \begin{align*}
       & [((V^{0}_{-\nu}((i\ker{d^*})^{-\iota^*}) \oplus V^{0}_{-\nu}(\Gamma(S)^I))^+, \dim(V^{0}_{-\nu}((i\ker{d^*})^{-\iota^*})), \dim_{\mathbb{K}_G}(V^{0}_{-\nu}(\Gamma(S)^I))+n(Y, \mathfrak{t}, g)/2)]\\
       & =[S^0, 0,n(T^3, \mathfrak{s}, g)/2 ]
    \end{align*}
where $g$ is the flat metric on $T^3$ and $n(T^3, \mathfrak{s}, g)$ is the correction term that is given by 
\[
    n(T^3, \mathfrak{s}, g)=\ind_{\C}(D^+)+\frac{\sigma(D^2 \times T^2)}{8} 
\]
where $D^+$ is the spin Dirac operator on $D^2 \times T^2$ with boundary operator $D_{T^3}$ has trivial kernel. There is a diffeomorphism that reverse the orientations of $D^2 \times T^2$ and its boundary $T^3$. Moreover, that diffeomorphism is isometry on $T^3$. Then we have $\ind_{\C}(D^+)=0$ and we have $\sigma(D^2 \times T^2)=0$, then $n(T^3, \mathfrak{s}, g)=0$. 
\end{proof}

Then using the gluing result \cref{gluing_conn}, we have the following product theorem of the degree. 
\begin{proposition}\label{product_degree}
    Let $K$ be a knot in $S^3$ and $k$ be an even integer. Let $(\mathfrak{s}, I)$ is a Real spin structure on $\tilde{X}_K$. 
    We assume that $(\mathfrak{s}, I)$ can be extended to $\Sigma_2(S^3, K) \times S^1$. Then the relative Real Bauer--Furuta invariant is represented by sequence of maps
    \[
       \{ f_n \colon U_n' \to U_n \}_n
    \]
     if we forget $\pm 1$-action to the spinor. 
    We define the degree $deg(\tilde{X}_K, \mathfrak{s},I)$ of $(\tilde{X}_K, \mathfrak{s}, I)$ by mapping degree $deg(f_n)$ of $f_n$ for sufficiently large $n$. Similarly, we define the degree $deg(\tilde{P}, \mathfrak{s}', I')$ of $(\tilde{P}, \mathfrak{s}', I')$, where $(\mathfrak{s}', I)$ are Real spin structure on $\tilde{P}$. 
    If there is a $\Z/2$-equivariant map $g \colon \partial \tilde{P} \to \partial \tilde{X}_K$ such that $g^*((\mathfrak{s}, I)|_{\partial \tilde{X}_K})=(\mathfrak{s}', I')$. 
    Then we have 
    \begin{equation}\label{prod_formula_deg}
        \lvert deg(\tilde{P} \cup_g \tilde{X}_K, \mathfrak{s} \cup_g \mathfrak{s'}, I \cup_g I' )\rvert =\lvert deg(\tilde{X}_K, \mathfrak{s},I) \rvert \lvert deg(\tilde{P}, \mathfrak{s}', I') \rvert. 
    \end{equation}
\end{proposition}
\begin{proof}
    From the assumption that $\mathfrak{s}$ can be extended to $\Sigma_2(S^3, K) \times S^1$, from \cref{real_seiberg_witten_of_T3}, the Real Bauer--Furuta invariant is a map between sphere spectrums. 
    The Spaniel--Whitehead duality of the sphere spectrum is the smash product. 
    From \cref{gluing_conn}, the Real Bauer--Furuta invariant of $(\tilde{P} \cup_g \tilde{X}_K, \mathfrak{s} \cup_g \mathfrak{s'}, I \cup_g I' )$ is the composition of the smash product and the Spaniel--Whitehead duality map of the Real Seiberg--Witten Floer homotopy type of the boundary. 
     Therefore, the mapping degree of the smash product of two maps between sphere spectrum is product of mapping degrees of each map. Then we proved this proposition. 
\end{proof}

Using the above product formula, we can calculate the degree of $\tilde{P}$. 
\begin{proposition}\label{calculations_for_P}
    Let $(\mathfrak{s}', I')$ be the Real spin structure on $\tilde{P}$ that satisfies $\mathfrak{s}'|_{\partial \tilde{P}}=\mathfrak{s}_{T^3} \otimes_{\R} l_{\epsilon_1+\epsilon_2+\epsilon_3}$. 
    Then we have $\lvert deg(\tilde{P}, \mathfrak{s}', I') \rvert =1$. 
\end{proposition}
\begin{proof}
    Let us take the unknot $U$ in $S^3$. From \cref{gluing_double_branched_cover} and \cref{condition_of_alpha}, we have 
    \[
        \Sigma_2(S^4, \tau_{2, 0}(U))=\tilde{P} \cup_{g(\tilde{A})} \tilde{X}_U
    \]
    and the restriction of the unique Real spin structures on $\Sigma_2(S^4, \tau_{2, 0}(U))$ to $\tilde{P}$ is $(\mathfrak{s}', I')$. 
    One can easily check that $\tau_{2, 0}(U)$ is isotopic to $U$ by using another description of the $k$-twist spun knot that is given in \cite[Section 6.2]{MR3588325}. 
    We show that the absolute value of the degree of the unknot is $1$. $\Sigma(S^4, U) \cong S^4$ has positive scalar curvature (PSC) metric and that is preserved by the covering involution. From $H^2(S^4)=0$ and $H^1(S^4)=0$, we have that there is only one reducible solution of Real Seiberg--Witten equations. Moreover, this is transverse since the metric is PSC. 
    From \cref{product_degree}, we have
    \[
        1=\lvert deg(\tilde{X}_U, \mathfrak{s},I) \rvert \lvert deg(\tilde{P}, \mathfrak{s}', I') \rvert
    \]
    and we have $\lvert deg(\tilde{P}, \mathfrak{s}', I') \rvert=1$. 
\end{proof}

Next we calculate the degree $\lvert deg(\tilde{X}_K, \mathfrak{s},I) \rvert$. 
\begin{lemma}
    Let $K$ be a knot in $S^3$. 
    Let $(\mathfrak{s}, I)$ is a Real spin structure on $\tilde{X}_K$. 
    We assume that $(\mathfrak{s}, I)$ can be extended to $\Sigma_2(S^3, K) \times S^1 \cong X_K \cup_{\tilde{S}^1 \times T^2} 
     \tilde{D}^2 \times T^2$, where $\tilde{S}^1$ and $\tilde{D}^2$ are $S^1$ and $D^2$ with the involution $z \mapsto -z \in S^1 \subset D^2 \subset \C$. 
    Let us denote $(\tilde{\mathfrak{s}}, \tilde{I})$ by the extended Real spin structure on $\Sigma_2(S^3, K) \times S^1$. 
    Then we have 
    \[
        \lvert deg(\tilde{X}_K, \mathfrak{s},I) \rvert = \lvert deg(\Sigma_2(S^3, K) \times S^1, \tilde{\mathfrak{s}}, \tilde{I}) \rvert
    \]
   , where $deg(\Sigma_2(S^3, K) \times S^1, \tilde{\mathfrak{s}}, \tilde{I})$ is a mapping degree of the Real Bauer--Furuta invariant of $(\Sigma_2(S^3, K) \times S^1, \tilde{\mathfrak{s}}, \tilde{I})$. 
\end{lemma}
\begin{proof}
    Note that we can define $\lvert deg(\Sigma_2(S^3, K) \times S^1, \tilde{\mathfrak{s}}, \tilde{I}) \rvert$ because $H^1(\Sigma_2(S^3, K) \times S^1)^{-\iota^*}=0$ and $b_2(\Sigma_2(S^3, K) \times S^1)=0$, we have that the Real Bauer--Furuta invariant is in $\pi_0^{st}(pt) \cong \Z$ if we forget $\pm 1$-action on the spinor bundle. 
    From \cref{gluing_conn} and \cref{real_seiberg_witten_of_T3}, we have that 
    \[
        \lvert deg(\tilde{X}_K, \mathfrak{s},I) \rvert \lvert deg(\tilde{D}^2 \times T^2, \tilde{\mathfrak{s}}|_{\tilde{D}^2 \times T^2}, \tilde{I}|_{\tilde{D}^2 \times T^2}) \rvert = \lvert deg(\Sigma_2(S^3, K) \times S^1, \tilde{\mathfrak{s}}, \tilde{I}) \rvert. 
    \]
    By considering the case that the knot $K$ is the unknot, we can prove $\lvert deg(\tilde{D}^2 \times T^2, \tilde{\mathfrak{s}}|_{\tilde{D}^2 \times T^2}, \tilde{I}|_{\tilde{D}^2 \times T^2}) \rvert=1$ in the same way in \cref{calculations_for_P}. 
\end{proof}

We are left to calculate the degree $\lvert deg(\Sigma_2(S^3, K) \times S^1, \tilde{\mathfrak{s}}, \tilde{I}) \rvert$. 

\begin{definition}
    We assume a knot $K$ in $S^3$ satisfies $\lvert \Delta_K(-1)\rvert=1$. Let $(\mathfrak{t}_K, I)$ be the unique Real spin structure on $\Sigma_2(S^3, K)$. 
    We take that a perturbation of the $3$-dimensional Real Seiberg--Witten equations such that all the solutions are transverse. 
    Then we define $\lvert deg(K) \rvert$ by the absolute value of the signed count of the framed moduli space. 
\end{definition}
\begin{remark}
    In the Real setting, the framed moduli space is a finite points and the reducible solution is also a point. Then $\lvert deg(K) \rvert$ coincides with 
    \[
       \lvert 2\#M^*(\Sigma_2(S^3, K), \mathfrak{t}_K, I)\pm 1 \rvert
    \] 
   , where $M^*(\Sigma_2(S^3, K), \mathfrak{t}_K, I)$ is the moduli space of the irreducible Real Seiberg--Witten equations and $\#$ means signed count. The factor $2$ comes from the $\pm1$-action to the spinor bundle and $\pm 1$ comes from the reducible solution. 
    Note that the signed count of all solutions of that perturbed Real Seiberg--Witten equations is independent of the choice of the perturbation. 
\end{remark}

\begin{proposition}\label{counting_real_sw}
    We assume a knot $K$ in $S^3$ satisfies $\lvert \Delta_K(-1)\rvert=1$. Let $(\mathfrak{t}_K, I)$ be the unique Real spin structure on $\Sigma_2(S^3, K)$. Let $\mathfrak{s} =\mathfrak{t}_K \times \mathfrak{t}'$ be the product spin structure on $\Sigma_2(S^3, K) \times S^1$, where $\mathfrak{t}'$ is the spin structure on $S^1$ that is invariant under the $S^1$ action. We denote $l$ be the M{\"o}bius bundle on $S^1$. 
    Then we have
    \begin{itemize}
        \item $\lvert deg(\Sigma_2(S^3, K) \times S^1, \tilde{\mathfrak{s}}\otimes_{\R} l , \tilde{I}) \rvert=1$, and 
        \item $\lvert deg(\Sigma_2(S^3, K) \times S^1, \tilde{\mathfrak{s}}, \tilde{I}) \rvert$ coincides with $\lvert deg(K) \rvert$. 
    \end{itemize}
\end{proposition}
\begin{proof}
    Let us take a perturbation of $3$-dimensional Real Seiberg--Witten equations with all the solution of the perturbed equations are transverse. From \cref{MxS1}, the second bullet of the statement holds. Let us prove the first statement. We show that there is no irreducible soluiton of the Real Seiberg--Witten equations on the Real $\mathrm{spin}^c$ structure $(\tilde{\mathfrak{s}}\otimes_{\R} l , \tilde{I})$. Note that $\tilde{\mathfrak{s}}\otimes_{\R}l$ is isomorphic to $\tilde{\mathfrak{s}}$ as a $\mathrm{spin}^c$ structure. 
    Let $F$ be a map from the configuation space with respect ot the $\mathrm{spin}^c$ structure $\tilde{\mathfrak{s}}$ to that of $\tilde{\mathfrak{s}}\otimes_{\R}l$ given as follows: 
    \[
       F \colon  (A, \phi) \mapsto (A-\frac{idt}{2}, e^{it/2}\phi)
    \]
   , where $t \in [0, 2\pi]$. Note that this map preserves the solution of the Seiberg--Witten equations however does not $I$-equivariant. 
    Let $F(A, \Phi)$ be an irreducible solution of the Real Seiberg--Witten equations on $(\Sigma_2(S^3, K) \times S^1, \tilde{\mathfrak{s}}\otimes l, \tilde{I})$. Then $(A, \Phi)$ is a solution of the ordinary Seiberg--Witten equations. 
    From the argument in the proof of \cref{bijectivity}, there exist the gauge transformation $u$ on on $\Sigma_2(S^3, K) \times S^1$ such that $u \cdot (A, \Phi)$ is pull back of the $3$-dimensional solution $(B, \phi)$, i.e. 
    $u \cdot (A, \Phi)=\pi^*(B, \phi)$. 

    Then $u^{-1}F(\pi^*(B, \phi))$ is an irreducible solution of the Real Seiberg--Witten equation. Thus we have 
    \[
        Ie^{it/2}u^{-1} (B, \phi)= e^{it/2}u^{-1}(B, \phi) \;\text{i.e.} I(B, \phi)=u^{-2} e^{it} (B,\phi) .
    \] 
 However, since $I(B, \phi)$ and $(B,\phi)$ are independent of $t$ and $\phi$ is not identically $0$, $u^{-2}e^{it}$ is independent $t$ and 
 \[
   \frac{1}{2\pi i} \int_{S^1} 2u^{-1}du=\frac{1}{2\pi i}\int_{S^1} idt=1.
 \] 
 This is a contradiction since $\int_{S^1} 2u^{-1}du/2\pi i$ is an even integer. 
\end{proof}

We give an explicit example of the knot $K$ in $S^3$. 
\begin{proposition}\label{237}
    Let $K$ be the $(-2, 3, 7)$ pretzel knot. Then the double branched cover of $K$ is an integer homology $3$-sphere and $\lvert deg(\Sigma_2(S^3, K)) \rvert=3$. 
\end{proposition}
\begin{proof}
    The double branched cover of $S^3$ along $K$ is the orientation reversed Brieskorn sphere $-\Sigma(2, 3, 7)$. The covering involution is the complex conjugate. Thus we have $\Sigma_2(S^3, K)$ is a homology $3$-sphere. 
    Let us take a metric $g$ on $-\Sigma(2, 3, 7)$ and a connection $\nabla^{\circ}$ on the spinor bundle of the spin structure $-\Sigma(2, 3, 7)$ as defined by Mrowka--Ozsv\'{a}th--Yu in \cite[Section 5.2]{MR1611061}. 
    From the result of Mrowka--Ozsv\'{a}th--Yu \cite[Theorem 1.]{MR1611061}, the framed moduli space $\mathcal{M}_0$ of the $3$-dimensional ordinary Seiberg--Witten equations on $-\Sigma(2, 3, 7)$ is diffeomorhic to $\{pt\} \sqcup S^1 \sqcup S^1$. From \cref{coincides_with_two_quotients}, the $I$-fixed point set of $\mathcal{M}_0$ coincides with the framed moduli space of the Real Seiberg--Witten equations $\tilde{\mathcal{M}}^I$ if the metric $g$ and the connection $\nabla^{\circ}$ are preserved by the involtuion $I$. 

    Firstly, we prove that the metric and the connection are preserved by $I$. Recall the definition of $g$ and $\nabla^{\circ}$. 
    The Brieskorn sphere $-\Sigma(2, 3, 7)$ is a sphere bundle of an orbifold line bundle $\pi \colon N \to \Sigma$ on an orbifold $\Sigma$. 
    Let $\eta$ be a $1$-form that is $i\eta$ is a connection form on the sphere bundle of $N$ with constant curvature, i.e. $id\eta=i\xi d\mu_{\Sigma}$, where $d\mu_{\Sigma}$ is a volume form on $\Sigma$ and $\xi$ is for some real number. 
    The metric $g$ is given by $g=\eta^2 + \pi^* g_{\Sigma}$ here $g_{\Sigma}$ is an orbifold Riemannian metric on ${\Sigma}$. 
    Let $\iota \colon -\Sigma(2, 3, 7) \to -\Sigma(2, 3, 7)$ the complex conjugate. 
    Then $\iota$ reverses the $S^1$-action on $-\Sigma(2, 3, 7)$ so that if we take $g_{\Sigma}$ preseved by the involution, we have $\iota^*\eta=-\eta$ and $\iota^*g=g$. 
    From \cite[Lemma 5.5.]{MR1611061} the connection $\nabla^{\circ}$ is given by a linear comnination of the Spin connection induced by the Levi-Civita connection of $g$, Clifford multiplication, and a quadratic term of $\eta$. 
    Thus $\nabla^{\circ}$ is preserved by $I$. 
    
    Secondly, we consider the $I$-action on $\mathcal{M}_0$. It is clear that $I$ preserves the reducible solution. There are two possibility of the action on the irreducible solutions: $I$ preseves the connected components or not. 
    If $I$ preserves the connected components, sinse $I$ reverses the $S^1$-action, one can easily seen that there are two points in each connected component so that we have $4$-irreducible solutions in $\tilde{\mathcal{M}}^I$. 

    We prove that each connected components are preserved by $I$. From \cite[Section 5.5.]{MR1611061}, each components are characterized as follows: Let $\rho(\eta)$ be the Clifford multiplication of $\eta$. Then $(i\rho(\eta))^2=1$ and its eigenvalue is $\pm 1$. The spinor of the irreducible solutions are the section of one of these eigenspaces. 
    On the other hand, we have 
    \[
            I(i\rho(\eta)\phi)=-i\rho(\iota^* \eta)I(\phi)=i\rho(\eta)I(\phi) 
        \]
        for all $\phi \in \Gamma(S)$. 
    Therefore, $I$ is commutative with $i\rho(\eta)$ and each components of the irreducible solutions in $\mathcal{M}_0$ are preserved by $I$-action.

    To calculate $\lvert deg(\Sigma_2(S^3, K)) \rvert$, we have to calculate the spectral flow between the reducible solution and irreducible solutions. We shall prove that spectral flow is $\pm 1$. 

In our situation, the line bundle $\lfloor \frac{k}{2} \rfloor$, that is given in \cite[Definition 10.1]{MR1611061} is the trivial orbifold bundle on $\Sigma$. 
Therefore, from \cite[Theorem 10.12]{MR1611061}, we have the dimension of the ordinary moduli moduli space of the flow between an irreducible solution and the reducible solution is $1$ since $\dim \mathcal{D}(\underline{\C}, \underline{\C})=0$. 

In the ordinary case, from \cite[Theorem 10.6.]{MR1611061}, we can identify with the ordinary framed moduli space $\mathcal{M}_0([\gamma], [\gamma_0])$ of the flow between an irreducible solutions $\gamma$ and the reducible solutions $\gamma_0$ is identified with $\mathcal{D}^{\circ}(\underline{\C}, \underline{\C})$, which is the based moduli space of the divisors on the rules surface that does not intersect $\Sigma^{\pm}$, defined in \cite[Definition 10.4.]{MR1611061}. 
In our setting, we can easily seen that all elements of $\mathcal{D}^{\circ}(\underline{\C}, \underline{\C})$ is represented by the constant section on the ruled surface so that $\mathcal{D}^{\circ}(\underline{\C}, \underline{\C}) \cong \C^*$. 
Then the connected components of $\mathcal{M}_0([\gamma], [\gamma_0])$ is $1$ and the $I$ fixed point set of this framed moduli space is not empty and $1$-dimensional since $I$ kills the $\mathrm{U}(1)$-action. 
This moduli space is transverse so that the Atiyah--Patodi--Singer index of the linealization of the Seiberg--Witten flow is $1$. This implies that the spectral flow is $1$. 
\end{proof}

Using \cref{237} and \cref{counting_real_sw}, we have the following theorem:

\begin{theorem}\label{237ka}
    Let $K$ be a $(-2, 3, 7)$-pretzel knot and $k$ be an even integer. If $\frac{k}{2}+\alpha$ be an even integer, $\lvert deg(\tau_{k, \alpha}(K))\rvert =1$. If $\frac{k}{2}+\alpha$ be an odd integer, we have $\lvert deg(\tau_{k, \alpha}(K))\rvert =3$
\end{theorem}
\begin{proof}
    Since the double branched cover of $S^3$ along $K$ is the Brieskorn sphere $-\Sigma(2, 3, 7)$ and this is integer homology sphere, this knot satisfies the assumption of \cref{counting_real_sw}. 
\end{proof}

\subsection{$P^2$-knots}

In this section, we calculate examples of $P^2$-knot. 
\begin{proposition}\label{standard_P2_knot}
    Let $P_{\pm}$ be a standard $P^2$-knot with $P_{\pm} \circ P_{\pm}=\pm 2$. Note that $\Sigma_2(S^4, P_+)=\bar{\C P^2}$, $\Sigma_2(S^4, P_-)=\C P^2$, and the covering involution is the complex conjugate. 
Then we have $deg(P_+)=1$ and $deg(P_-)=0$. 
\end{proposition}
\begin{proof}
    We see that the PSC metric of $\C P^2$ is preserved by the complex conjugate. Therefore, in both cases, all the solutions of the Real Seiberg--Witten equations are reducible. 
    
    First, we prove $deg(P_-)=0$. If $(A, 0)$ be a solution of the Real Seiberg--Witten equations, $F_{A^{\tau}}^+=0$. However, the de Rham cohomology class of $F_{A^{\tau}}^+$ is $c$ and this is not $0$. Therefore there is no solutions in this case. 

    Second, we prove $deg(P_+)=1$. This case we have $b^+_{-\iota^*}(\bar{\C P^2})=0$ so that we have the unique reducible solution up to gauge transformation. Moreover, we can easily check that this is transverse solution since the Real Dirac index is $0$ and we have PSC metric. 
\end{proof}

From the computation above and \cref{conn_sum_S_P}, we have following corollaries. 
\begin{corollary}\label{kinoshita}
    If a $P^2$-knot $P$ is given by the connected sum of $P_-$ and for some $2$-knot $S$. Then its degree is $0$ for all Real $\mathrm{spin}^c$ structure on $\Sigma_2(S^4, P_- \# S)$. 
    In particular, if there is a $P^2$-knot $P$ with $P \circ P=-2$ that satisfies $\lvert deg(P, c) \rvert \neq 0$ for some Real $\mathrm{spin}^c$ structure, the $P^2$-knot $P$ is a counter example of the Kinoshita conjecture. 
\end{corollary}
\begin{corollary}\label{computation_of_P}
    If a $P^2$-knot $P$ is given by the connected sum of $P_+$ and for some $2$-knot $S$. For simplicity, we assume that $\Sigma_2(S^4, S)$
is homology $S^4$. 
Then we have 
\[
    \lvert deg(P) \rvert =\lvert deg(S) \rvert. 
\]
\end{corollary}

From \cref{computation_of_P}, \cref{multiplicative}, and \cref{237ka}, we have the following theorem:

\begin{theorem}\label{main_theorem_1}
    Let $K$ be the $(-2, 3, 7)$-pretzel knot, $k$ be an even integer, $n$ be a non-negative integer, and $\alpha$ be an integer. Let $P_{k, \alpha}^n$ be a $P^2$-knot that is given by the connected sum of $P_+$ and $n$ times connected sum of $\tau_{k, \alpha}(K)$: $P_+ \# n\tau_{k, \alpha}(K)=P_+ \# \tau_{k, \alpha}(K)\# \tau_{k, \alpha}(K)\dots \# \tau_{k, \alpha}(K)$. 
    Then if $\frac{k}{2}+\alpha$ be an even integer, $\lvert deg(P_{k, \alpha}^n) \rvert=1$ for all $n$. 
    If $\frac{k}{2}+\alpha$ is odd, $\lvert deg(P_{k, \alpha}^n) \rvert=3^n$ and then we have $(S^4, P_{k, \alpha}^n)$ and $(S^4, P_{k, \alpha}^{n'})$ are not pairwise diffeomorphic if $n\neq n'$.   
\end{theorem}

From \cref{main_theorem_1}, we can partially give the answer of the question that is given by Kamada \cite[Epilogue]{MR3588325}. 
\begin{corollary}\label{even_twist}
    Let $K$ be the $(-2, 3, 7)$-pretzel knot and $k$ and $k'$ be an even integer. If $k - k'$ is not divided by $4$, then $P_+ \# \tau_{k, 0}(K)$ and $P_+ \# \tau_{k', 0}(K)$ are not smoothly isotopic. 
\end{corollary}

Next, we prove that the existence of the exotic $P^2$-knot. 
\begin{theorem}\label{main_theorem_2}
    Let $n>0$ be a natural number.  Let $K$ be the $(-2, 3, 7)$-pretzel knot. Then $P_{0, 1}^n\colon =P_+ \# n \tau_{0, 1}(K)$ is topologically isotopic to $P_+$. Here we denote $n \tau_{0, 1}(K)$ by the connected sum of $n$ copy of $\tau_{0, 1}(K)$. However, $P_{0, 1}^n$ is not smoothly isotopic to $P_+$. Moreover, $(S^4, P_{0, 1}^n)$ and $(S^4, P_+)$ are not pairwise diffeomorphic. 
\end{theorem}
The latter half of the statement immediately follows from \cref{main_theorem_1}. Let us prove that $P_{0, 1}^n$ and $P_+$ are topologically isotopic. 

We use the following result obtained by Conway--Orson--Powell \cite[Theorem A.]{conway2023unknotting}. 
\begin{theorem}[Conway--Orson--Powell\cite{conway2023unknotting}]
    Let $F \subset S^4$ be a $\Z/2$-surface of non-orientable genus $h$ and normal Euler number $e$. If $h \le 3$ or $\lvert e \rvert \neq 2h$, then $F$ is unknotted.  
\end{theorem}
In our case, the non-orientable genus of $\R P^2$ is $1$ so that we can apply this theorem. 
We prove the following theorem.  
\begin{theorem}\label{fundamental_group_general}
    Let $K$ be a $(-2, 3, 7)$ pretzel knot. Let $X$ be a $4$-manifold and $S \subset X$ be an embedded surface. Suppose $\pi_1(X\setminus S)\cong \Z/2$ and its generator is the meridian of hte surface $S$. Then $\pi_1(X \# S^4 \setminus S \# \tau_{0, 1}(K))\cong \Z/2$ and its generator is the meridean of the embedded surface $S \# \tau_{0, 1}(K)$. 
\end{theorem}
From \cref{fundamental_group_general}, we have the following propositions. 
\begin{proposition}\label{fundamental_group}
    The fundamental group of $\pi_1(S^4\setminus P_{0, 1}^n)$ is $\Z/2$. 
\end{proposition}
\begin{proof}
    We see the fundamental group $\pi_1(S^4, P_+) \cong \Z/2$ and the generator is the meridean of $P_+$. Then we can apply \cref{fundamental_group_general} $n$ times. Then we have this result. 
\end{proof}
\begin{theorem}\label{conic}
    Let $C \subset \C P^2$ be a conic, which is a holomorphic embedding of $S^2$ whose homology class is $2[\C P^1] \in H_2(\C P^2, \Z)$. 
    Let $C \# n \tau_{0, 1}(K)$ be an embedding of $S^2$ in $\C P^2 \cong \C P^2 \# S^4$ given by inner connected sum of $C$ and $n \tau_{0, 1}(K)$, the $n$ times inner connected sum of $\tau_{0, 1}(K)$. 
    Then $\pi_1(\C P^2 \setminus C \# n\tau_{0, 1}(K))\cong \pi_1(\C P^2\setminus C) \cong \Z/2$ however $(\C P^2, C)$ and $(\C P^2, C \# n\tau_{0, 1}(K))$ are not pairwise diffeomorphic if $n \neq 0$. 
\end{theorem}
\begin{proof}
    The fundamental group $\pi_1(\C P^2 \setminus C \# n \tau_{0, 1}(K))$ is calculated by \cref{fundamental_group_general} since $\pi_1(\C P^2 \setminus C)\cong \Z/2$. 
    Let us prove $(\C P^2, C)$ and $(\C P^2, C \# n \tau_{0, 1}(K))$ are not pairwise diffeomorphic. 
    Note that the double branched cover of $\C P^2$ along $C$ is $S^2 \times S^2$ and the covering involtuions is $(x, y) \mapsto (y, x)$. 
    By \cref{the_invariant_of_S2xS2}, the formal dimension of the Real spin structure $(\mathfrak{s}, I)$ on $S^2 \times S^2$ is $0$ and the absolute value of the degree of the Real Bauer Furuta map $\lvert deg(S^2 \times S^2,\mathfrak{s}, I) \rvert$ is well-defined. 
    Moreover, \cref{the_invariant_of_S2xS2} tells us that $\lvert deg(S^2 \times S^2,\mathfrak{s}, I) \rvert=1$. 
    On the other hand, the double branched cover $\Sigma_2(\C P^2, C \# \tau_{0, 1}(K))$ have 
    an unique Real $\mathrm{spin}^c$ structure $(\mathfrak{s}', I')$ on $\Sigma_2(\C P^2, C \# n \tau_{0, 1}(K))$ since $H^2(\Sigma_2(\C P^2, C \# n \tau_{0, 1}(K)), \Z)^{-\iota^*}=0$, $H^1(\Sigma_2(\C P^2, C \# n \tau_{0, 1}(K)), \Z)^{-\iota^*}=0$, and $\Sigma_2(\C P^2, C \# n \tau_{0, 1}(K))$ has an unique spin structure, which preseved by $\iota$. From \cref{gluing_conn} and \cref{237ka}, we have that 
    \[\lvert deg(\Sigma_2(\C P^2, C \# n \tau_{0, 1}(K)), \mathfrak{s}', I') \rvert=3^n.\] Then we finished the proof. 
\end{proof}

From now, we prove \cref{fundamental_group_general}. 
\begin{proof}[Proof of \cref{fundamental_group_general}]
    From \cite[Section 4.]{MR0780592}, the fundamental group of $S^4\setminus \tau_{0, 1}(K)$ is given by 
    \[
        \langle \pi_1(X) \mid lxl^{-1}x^{-1}, x \in \pi_1(X), \; l \;\text{is longitude.} \rangle. 
    \]
    From the Van Kampen theorem, we obtain 
    \[
        \pi_1(X \# S^4 \setminus S \# \tau_{0, 1}(K))\cong \pi_1((X \setminus S))\ast_{\pi_1(D^3 \setminus U)}\pi_1(S^4\setminus \tau_{0, 1}(K))=\Z/2 \ast_{\Z} \pi_1(S^4\setminus \tau_{0, 1}(K)). 
    \]
    Therefore we have 
    \[
        \pi_1(X \# S^4 \setminus S \# \tau_{0, 1}(K))=\langle \pi_1(X) \mid lxl^{-1}x^{-1}, m^2, x \in \pi_1(X), \; l \;\text{is longitude},\;m \;\text{is meridian}\rangle.
    \]
    We now use the result of Nakae \cite[Proposition 2.1]{MR3024026}:
    \begin{theorem}[Nakae\cite{MR3024026}]\label{nakae}
        Let $K_s$ be a $(-2, 3, 2s+1)$-pretzel knot with $s \ge 3$. 
        Then the knot group has the representation
        \[
            \langle c, l \mid clcl^{-1}c^{-1}l^{-s}c^{-1}l^{-1}clcl^{s-1} \rangle
        \]
        and an element which represents the meridean $m$ is $c$ and an element of the longitude $L$ is $c^{-2s+2}lcl^{s}cl^{s}clc^{-2s-9}$. 
    \end{theorem}
    Now we show that if we add relations $c^2=1$ and $L$ is in center of the group, we have $l=1$. We know that we have a non trivial double cover of $X \# S^4 \setminus S \# \tau_{0, 1}(K)$, then the knot group is non trivial. Then we have $\pi_1(X \# S^4 \setminus S \# \tau_{0, 1}(K))$ is $\Z/2$ and the generator is $c$, the meridian of $S\# \tau_{0, 1}(K)$. 
    
From \cref{nakae} and $c^2=1$, we have that the knot group of $P_{0, 1}$ is
\[
    \langle c, l \mid clcl^{-1}cl^{-s}cl^{-1}clcl^{s-1}=1, c^2=1, L=lcl^{s}cl^{s}clc \;\text{is in center.} \rangle. 
\]
   , where $s=3$. We prove $l=1$ by 3 steps. Until \underline{Step 2.}, we can carry our caliculation out in general $s$. So that we caliculate in general $s$ until the final step. 
    \begin{itemize}
        \item \underline{Step 1.} Let us calculate $cl^{-1}L$ in two ways. Firstly, from the definition of $L$, we have 
        \[
            cl^{-1}L=cl^{-1}lcl^{s}cl^{s}clc=l^{s}cl^{s}clc. 
        \]
        Since $L$ commutes $cl^{-1}$, then 
        \[
            cl^{-1}L=Lcl^{-1}=lcl^{s}cl^{s}clccl^{-1}=lcl^{s}cl^{s}c. 
        \]
        Then we have 
        \[
            l^{s}cl^{s}clc=lcl^{s}cl^{s}c.
        \]
        Thus we have $cl^sc=l^{1-s}cl^scl^{s-1}$. 
        \item \underline{Step 2.} From the previous step, we have 
        \[
            cl^{-s}c=l^{1-s}cl^{-s}cl^{s-1}. 
        \]
        From $clcl^{-1} \cdot cl^{-s}c \cdot l^{-1}clcl^{s-1}=1$, we have 
        \begin{align*}
            1&=clcl^{-1} \cdot cl^{-s}c \cdot l^{-1}clcl^{s-1}\\
            &=clcl^{-1} \cdot l^{1-s}cl^{-s}cl^{s-1} \cdot l^{-1}clcl^{s-1}\\
            &=clc\cdot l^{-s}cl^{-s}cl^{s-2}clcl^{s-1}\\
            &=cl^2c\cdot cl^{-1}c l^{-s}cl^{-s}cl^{-1}\cdot l^{s-1}clcl^{s-1}.\\
        \end{align*}
        Now, we have $L^{-1}=cl^{-1}c l^{-s}cl^{-s}cl^{-1}$ since $L=lcl^{s}cl^{s}clc$ and $c^2=1$. Then
        \begin{align*}
            &cl^2c\cdot cl^{-1}c l^{-s}cl^{-s}cl^{-1}\cdot l^{s-1}c \cdot lcl^{s-1}\\
            &=cl^2c\cdot L^{-1} l^{s-1}c \cdot lcl^{s-1}\\
            &=cl^2cl^{s-1}c \cdot L^{-1} \cdot lcl^{s-1}\\
            &=cl^2cl^{s-1}c \cdot cl^{-1}c l^{-s}cl^{-s}cl^{-1} \cdot lcl^{s-1}\\
            &=cl^2c\cdot l^{s-2} \cdot c l^{-s}cl^{-1}\\
            &=cl^2c\cdot l^{2s-2}l^{-s}\cdot c l^{-s}cl^{-1}.\\
        \end{align*}
        We have $l^{-s}\cdot c l^{-s}cl^{-1}=clcL^{-1}$. Then 
        \begin{align*}
            &cl^2c\cdot l^{2s-2}\cdot l^{-s} c l^{-s}cl^{-1}\\
            &=cl^2c\cdot l^{2s-2}clcL^{-1}. 
        \end{align*}
        Therefore, we have 
        \[
            L=cl^2c\cdot l^{2s-2}clc=lc L cl^{-1}=lc \cdot  cl^2c l^{2s-2}clc \cdot cl^{-1}=l^3cl^{2s-2}c. 
        \]
        \item \underline{Step 3. }
        Since $cl^{-3}L=Lcl^{-3}$, we have 
        \begin{align*}
            cl^{-3}L&=cl^{-3}l^3cl^{2s-2}c\\
            &=l^{2s-2}c\\
            &=L cl^{-3}\\
            &=l^3cl^{2s-2}c cl^{-3}\\
            &=l^3cl^{2s-5}. 
        \end{align*}
        Therefore, by comparing the second line and fifth line, we have
        \[
            l^{2s-5}c=cl^{2s-5}. 
        \]
        Now, in our case, $s=3$. then we have $l$ is in the center of the group. From $clcl^{-1}cl^{-s}cl^{-1}clcl^{s-1}=1$, we have 
        \[
            clcl^{-1}cl^{-s}cl^{-1}clcl^{s-1}=c^6l^{1-1-s-1+1+s-1}=l^{-1}=1. 
        \]
    \end{itemize}
    Thus we have the desired result. 
\end{proof}

From \cref{fundamental_group}, we have the following theorem:
\begin{theorem}\label{fundamental_group_of_double_branched_cover}
    Let $K$ be the $(-2, 3, 7)$ pretzel knot. Then we have the following results. 
    \begin{itemize}
        \item The fundamental group of the double branched cover $\Sigma_2(S^4, n \tau_{0, 1}(K))$ is trivial. In particular, $\Sigma_2(S^4, n \tau_{0, 1}(K))$ is homotopy $S^4$. 
        \item The fundamental group of the double branched cover $\Sigma_2(S^4, P_+\# n \tau_{0, 1}(K))$ is trivial. In particular, $\Sigma_2(S^4, P_+\# n \tau_{0, 1}(K))$ is homotopy $\overline{\C P^2}$. 
    \end{itemize}
\end{theorem}
\begin{proof}
    Firstly, we prove the second bullet. From \cref{fundamental_group}, we have that the double cover of $S^4\setminus P_+ \# n\tau_{0,1}(K)$ is trivial since the kernel of the nontrivial map $\pi_1(S^4\setminus P_+ \# n\tau_{0,1}(K))\to \Z/2$ is $\{1\}$. From the Van Kampen's theorem, we see that the double branched cover is simply connected. 
    To prove the first bullet, let us denote $G$ by the fundamental group $\pi_1(\Sigma_2(S^4, \tau_{0,1}(K)))$. Then from the Van Kampen's theorem, we have 
    \[
        \pi_1(\Sigma_2(S^4, P_+ \# n\tau_{0,1}(K)))\cong \pi_2(\Sigma_2(S^4, P_+))\ast G \cong G. 
    \]
Since $\pi_1(\Sigma_2(S^4, P_+ \# n\tau_{0,1}(K)))=\{1\}$, we have $G$ is trivial. 
\end{proof}

From \cref{fundamental_group_of_double_branched_cover}, the double branched cover of $P_{0, 1}^n$ and $n \tau_{0, 1}$ are homotopy $\C P^2$ and $S^4$ respectively. However, we do not know wheter or not they are diffeomorphic to $\C P^2$ or $S^4$. We have the following two applications about the group action of $S^4$ and $\C P^2$. 

\begin{theorem}\label{Z2_action_CP2}
        At least one of the following two statements holds:
        \begin{itemize}
            \item There exist an of exotic $\C P^2$. 
            \item There exist an infinite family of exotic $\Z/2$-actions on $\C P^2$. 
        \end{itemize}
\end{theorem}
\begin{proof}
    Suppose the statement in the first bullet is folse. 
    The embeddings of $\R P^2$ to $S^4$ given in \cref{main_theorem_2} are topologically isotopic, so that these involutions on $\C P^2$ are topologically istopic since they are covering involutions of topologically isotopic surfaces. However, our invariant for $(\Sigma_2(S^4, P_{0, 1}^n), \iota)$ are distinct for different $n$ so that $(\Sigma_2(S^4, P_{0, 1}^n), \iota)$ and $(\Sigma_2(S^4, P_{0, 1}^m), \iota)$ are not $\Z/2$-equivariant isotopic if $n \neq m$.  
\end{proof}

\begin{theorem}\label{Z2_action_S4}
    At least one of the following two statements holds:
    \begin{itemize}
        \item There exist a counter example of the smooth Poincar\'{e} conjecture. 
        \item There exist infinite family of exotic $\Z/2$-actions on $S^4$ that do not preserve any positive scalar curvature metric. 
    \end{itemize}
\end{theorem}
\begin{proof}
    Suppose the statement in the first bullet is folse. 
    Then the branched covers $\Sigma_2(S^4, P_{0, 1}^m)$ are diffeomorphic to $S^4$. If the covering involtuion $\iota$ on $\Sigma_2(S^4, P_{0, 1}^m)$ admits a positive scalar curvature $g$ with $\iota^*g=g$, then from the our invariant should be $1$. This is a contradiction.  
\end{proof}

From \cref{main_theorem_1} and \cref{main_theorem_2}, we deduce that the complement of a tubular neighborhood of $P_+$ and $P_+\# n \tau_{0, 1}(K)$ are exotic $4$-manifolds. They are rational homology $4$-ball and the boundaty is a quotient of $S^3$. Note that they are absolutely exotic, not relative to the boundary.

\begin{theorem}\label{complement}
    Let $K$ be the $(-2, 3, 7)$-pretzel knot.
    Let $\nu_n:=\nu(P_+\# \tau_{0, 1}(K))$ be a tubular neighborhood of the $P^2$-knot $\nu(P_+\# \tau_{0, 1}(K))$. Then $X_n:=S^4\setminus \nu_n$ are homeomorphic to $S^4\setminus \nu(P_+)$ however they are not diffeomorphic to each other. 
\end{theorem}
\begin{proof}
    Firstly, they are homeomorphic since $P_+\#\tau_{0, 1}(K)$ are topologically isotopic. Let us prove that they are not diffeomorphic. 
    
    Let $\Sigma_2(\nu_n)$ be the double branched cover of $\nu_n$ along $P_+\# n\tau_{0, 1}(K)$. Let $\tilde{X}_n$ are the double cover of $X_n$. Let $c$ is the canonical Real spin$^c$ structure on $\Sigma_2(S^4, P_+\#\tau_{0, 1}(K))$ as menstioned in \cref{invariant_of_P2_knot}. 
    We calculate $\Psi(\tilde{X}_n, c|_{\tilde{X}_n})$. Note that the boundary $\partial \tilde{X}_n$ is a quotient of $S^3$ so that $\partial \tilde{X}_n$ has a PSC matric and the covering involtion preserves the PSC metric. Hence the Real Seiberg--Witten Floer homotopy type is of the form
    \[
        (S^0, n, m),
    \]
    where $n \in \Z , \; m\in \Q$. We see that $b^+(\tilde{X})=b^+(\Sigma_2(\nu_n))=0$, $(c|_{\tilde{X}})^2-\sigma(\tilde{X})=(c|_{\Sigma_2(\nu_n)})^2-\sigma(\Sigma_2(\nu_n))=0$ and the Real Bauer--Furuta invariants $\Psi(\tilde{X}_n, c|_{\tilde{X}_n})$ and $\Psi(\Sigma_2{\nu_n}, c|_{\Sigma_2{\nu_n}})$ are of the form \[
        f_n \colon U'_n \to U_n ,
    \]
    where $\dim U'_n=\dim U_n$. 
    From \cref{gluing_formula_conn}, we have
    \[
        \Psi(\Sigma_2(S^4, P_+\#\tau_{0, 1}(K)), c)\cong \Psi(\tilde{X}_n, c|_{\tilde{X}_n})\wedge \Psi(\Sigma_2{\nu_n}, c|_{\Sigma_2{\nu_n}}). 
    \]
Therefore, by taking the mapping degree of these Real Bauer--Furuta invariant, we have
\[
   3^n= \lvert deg(P_+\#n\tau_{0, 1}(K)) \rvert=\lvert deg \Psi(\tilde{X}_n, c|_{\tilde{X}_n}) \rvert \lvert deg \Psi(\Sigma_2{\nu_n}, c|_{\Sigma_2{\nu_n}}) \rvert. 
\]
Note that $(\Sigma_2(\nu_n), \iota)$ are mutually $\Z/2$-equivariant diffeomorhic, where $\iota$ is the covering involtuion. Then if we consider the case $n=0$, we deduce that $\lvert deg \Psi(\Sigma_2{\nu_n}, c|_{\Sigma_2{\nu_n}}) \rvert=1$ for all $n$. Therefore we have 
\[
    \lvert deg \Psi(\tilde{X}_n, c|_{\tilde{X}_n}) \rvert=3^n
\]
so that $\{(\tilde{X}_n, \iota)\}_n$ are not $\Z/2$-equivariant diffeomorhic to each other. We show that $X_n=\tilde{X_n}/\iota$ are not diffeomorhic to each other. Suppose that $X_n$ and $X_m$ are diffeomorhic. Since $\tilde{X}_n$ and $\tilde{X}_m$ are simply connected, then the diffeomorphism from $X_n$ to $X_m$ lift to their universal cover $\tilde{X}_n$ and $\tilde{X}_m$ as a $\Z/2$-equivariant diffeomorphism. This is a contradiction. 
\end{proof}


\begin{acknowledgement}
    The author  would like to show his deep appreciation to  Mikio Furuta, Hokuto Konno, Daniel Ruberman, Masaki Taniguchi, and Junpei Yasuda for helpful discusssions. 
    The author is supported by JSPS KAKENHI Grant Number 21J22979 and WINGS-FMSP program at the Graduate school of Mathematical Science, the University of Tokyo.
\end{acknowledgement}

\bibliographystyle{jplain}
\bibliography{bibidatabase}

\end{document}